\documentclass[12pt]{article}

\usepackage[labelfont=bf]{caption}

\usepackage[normalem]{ulem}
\usepackage{amsmath}
\newcommand{\stkout}[1]{\ifmmode\text{\sout{\ensuremath{#1}}}\else\sout{#1}\fi}

\usepackage{algorithm}
\usepackage{algpseudocode}
\usepackage{enumerate}

\usepackage{amsmath}
\DeclareMathOperator*{\argmin}{arg\,min}

\usepackage[margin=1.5cm]{geometry}

\usepackage{amsfonts,amsmath,amssymb}
\usepackage{fancyhdr}
\usepackage{graphicx}
\usepackage{float}
\usepackage[nottoc,notlot,notlof]{tocbibind}
\usepackage{amsmath}
\usepackage{tikz}
\usetikzlibrary{decorations.pathreplacing,calligraphy}

\newcommand{\balpha}{\mbox{\boldmath$\alpha$}}

\newcommand{\bpi}{\mbox{\boldmath$\pi$}}
\newcommand{\bphi}{\mbox{\boldmath$\phi$}}
\newcommand{\btheta}{\mbox{\boldmath$\theta$}}

\usepackage{rotating}

\newtheorem{exa}{Example}

\title{Markov Decision Process and Approximate Dynamic Programming for a Patient Assignment Scheduling problem
	}

\author{
Ma{\l}gorzata M. O'Reilly
~\thanks{email: malgorzata.oreilly@utas.edu.au}
~\thanks{School of Natural Sciences, University of Tasmania, Hobart TAS 7001, Australia.}
\ \
Sebastian Krasnicki
~\footnotemark[2]
\ \
James Montgomery~\thanks{School of Information \& Communication Technology, University of Tasmania, Hobart TAS 7001, Australia.}
\\ \ \ 
Mojtaba Heydar~\thanks{BHP, Perth, Western Australia.}
\ \   
Richard Turner~\thanks{School of Medicine, University of Tasmania, Hobart TAS 7001, Australia.}
\ \ 
Pieter Van Dam~\thanks{School of Nursing, University of Tasmania, Hobart TAS 7001, Australia.}
\ \ 
Peter Maree~\thanks{Strategy and Planning, Department of Health, Hobart TAS 7001, Australia.}
}

\date{\normalsize \today}

\begin{document}
\maketitle

\begin{abstract}
	We study the Patient Assignment Scheduling (PAS) problem in a random environment that arises in the management of patient flow in the hospital systems, due to the stochastic nature of the arrivals as well as the Length of Stay distribution. 
	
	We develop a Markov Decision Process (MDP) which aims to assign the newly arrived patients in an optimal way so as to minimise the total expected long-run cost per unit time over an infinite horizon. We assume Poisson arrival rates that depend on patient types, and Length of Stay distributions that depend on whether patients stay in their primary wards or not.
	
	Since the instances of realistic size of this problem are not easy to solve, we develop numerical methods based on Approximate Dynamic Programming. We illustrate the theory with numerical examples with parameters obtained by fitting to data from a tertiary referral hospital in Australia, and demonstrate the application potential of our methodology under practical considerations.
	
\end{abstract}

\bigskip

\bigskip

{\bf Funding:} \quad This research was supported by funding through the Australian Research Council Linkage Project LP140100152.

\section{Introduction}\label{sec:Intro}

In this paper we consider the Patient Assignment Scheduling (PAS) problem in which, at the start of each day, newly arrived patients are assigned to beds in different wards, considering their needs, priority, and available resources, in a way so as to optimise total expected daily cost over an {\em infinite} time horizon. We use the term `bed' in the sense of the resource that consists of the staff (nurses and clinicians) available to attend to a patient in a physical bed, rather than the physical bed itself.

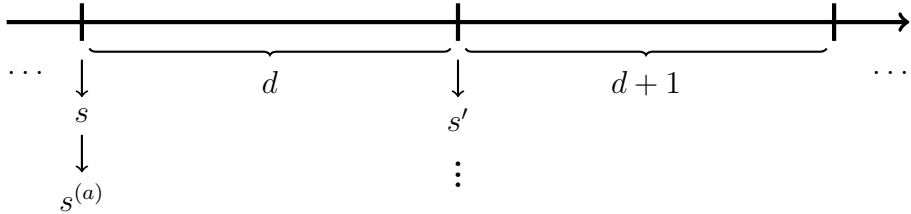
\begin{figure}[H]
	\begin{center}
		\begin{tikzpicture}
		\draw [thick][decorate, decoration = {brace, mirror}] (1.1,-0.35) --  (5.9,-0.35);
		\draw [thick][decorate, decoration = {brace, mirror}] (6.1,-0.35) --  (10.9,-0.35);
		
		\draw [->][thick] (1,-0.5) --  (1,-1);
		\node [below] at (1,-1) {$s$};
		\draw [->][thick] (1,-1.5) --  (1,-2);
		\node [below] at (1,-2) {$s^{(a)}$};
		
		\draw [->][thick] (6,-0.5) --  (6,-1);
		\node [below] at (6,-1) {$s'$};
		\node [below] at (6,-1.5) {\bf $\vdots$};
		
		\draw [->][ultra thick] (0,0) -- (12,0);
		\draw [-][ultra thick] (1,0.25) -- (1,-0.25);
		\draw [-][ultra thick] (6,0.25) -- (6,-0.25);
		\draw [-][ultra thick] (11,0.25) -- (11,-0.25);
		
		\node [below] at (0.25,-0.5) {\bf $\ldots$};
	\node [below] at (3.5,-0.5) {$d$};	
\node [below] at (8.5,-0.5) {$d+1$};				
		\node [below] at (11.75,-0.5) {\bf $\ldots$};
		
		
				
		\end{tikzpicture}
			\caption{Decision epochs of the PAS problem. At the start of time period $d$ we observe some state~$s$ and then make some decision~$a$ about how to allocate/transfer patients. This transforms the system into a post decision state $s^{(a)}$ and then the system evolves in a stochastic manner, until new state $s'$ is observed at the start of the next time period $d+1$.}
			\label{D_Epoch}  
	\end{center}	
\end{figure}

We assume that that patients are assigned to the beds at the start of time period $d=0,1,2,\ldots$, see Figure~\ref{D_Epoch}. The allocation involves a {\em group} of patients, waiting to be admitted to a suitable bed. Upon ceompleting their treatments, patients are discharged. The duration of time between arrival and discharge is a random variable referred to as the {\em Length of Stay} (LoS). We assume an {\em infinite} horizon problem and develop a Markov Decision Process (MDP) to solve this stochastic problem so as to minimise the total expected long-run cost per unit time.

Here, without loss of generality we assume that the index $d=0,1,2,\ldots$ corresponds to {\em day} $d$, but note that this could denote some other time interval of interest such as an $8$-hour time block.

Patient assignment is a key process in the management of hospitals, whereby bed managers allocate inbound patients to appropriate wards, rooms and beds. These decisions are made by teams of ward and bed managers at intermittent times throughout the day, using their experience and understanding of the natural dynamics of the hospital. From 2018 to 2019, the Australian hospital system recorded 11.5 millions hospitalisations, and 8.4 million emergency department presentations, resulting in 365,000 emergency admissions~\cite{AIHW_2018}. This corresponded to 31,500 patient assignment decisions made daily across Australia, with 1000 of these being unplanned~\cite{AIHW_2018}.

These decisions are not trivial, since complications can arise due to variations in hospital policy and the highly stochastic and complex nature of healthcare. Potential patients can arrive at any time of day, with a variety of conditions and specific needs. Hospital policy may affect where patients can be assigned, on the basis of a patient's gender or to comply with prescribed healthcare standards. Managers may also evaluate whether it is best to transfer patients between wards. Because of these complicating factors it can be unclear if the assignments made by bed managers are optimal. 

The importance of maintaining good patient flow cannot be understated. Carter, Puch and Larson in their literature review of emergency department (ED) crowding found that ED overcrowding has a significant positive correlation with patient mortality and with patients leaving the hospital untreated~\cite{Carter2014}. Morley et al. in their literature review also found an increase in patient mortality, as well as a higher exposure to error, poorer patient outcomes and increased patient length of stay, both in the ED, and in the ward to which a patient is eventually assigned~\cite{Morley2018}. They also noted that the inability to quickly assign patients from the ED to an inpatient bed, was a significant contributor to overcrowding. 

Regarding overcrowding in non-emergency wards, there is a long standing debate on what level of hospital occupancy is too high. Bagust et al.~\cite{Bagust1999} found in their simulation-based approach that hospital occupancy above $85\%$ created discernible risk, while occupancy above $90\%$ caused regular bed crises. This figure has since been adopted by many as a standard benchmark for hospital occupancy~\cite{Cummings2012}. However, multiple other authors~\cite{Bain2010,Cummings2012} have emphasized that adopting $85\%$ as the figure at which hospitals operate most efficiently, misinterprets Bagust et al.~\cite{Bagust1999}, as this figure only applies to the specific example in~\cite{Bagust1999}. 
Further, Bain et al.~\cite{Bain2010} state that using any single figure as a target for average occupancy is overly simplistic. Cummings et al.~\cite{Cummings2012} note that while the results of~\cite{Bagust1999} have been misinterpreted, there is a positive correlation between hospital occupancy and emergency department delays. 

Solutions to these problems require improvements in the management of patient flow. Howell et al.~\cite{Howell2008} found that by assigning additional resources to the tasks of bed management and patient assignment, ED length of stay was reduced by an average of 98 minutes for admitted patients. Clearly, tangible benefits are made to patients if they are assigned correctly.

We call the mathematical formulation of patient assignment, the Patient Assignment Scheduling (PAS) problem, where we attempt to find optimal assignments of all patients as they arrive to the hospital. This problem was described by Demeester et al in~\cite{Demeester2010}, Bilgin et al. in~\cite{Bilgin2012}, and Vancroonenburg et al. in~\cite{2016Van}.

Demeester et al.~\cite{Demeester2010} define a useful approach of hard and soft constraints. Patient assignments that violate hard constraints are infeasible, while assignments that violate soft constraints incur some cost, which contribute to the objective function. An example of a hard constraint is that two patients cannot be assigned to the same bed, while an example of a soft constraint is that a patient should be assigned to the ward that corresponds to the their condition.
In their model, the variable for patient length of stay (LoS) used  is deterministic, and so too is the cost of assigning a patient. This approach ignores the inherent stochasticity of a hospital system, but provides a useful framework for building upon.

Bilgin et al.~\cite{Bilgin2012} and Vancroonenburg et al. in~\cite{2016Van} use a very similar approach, defining a set of constraints to form an objective function, but each adding their own novel idea. Bilgin et al.~\cite{Bilgin2012} combine the PAS problem with a nurse scheduling problem, while Vancroonenburg et al. in~\cite{2016Van} introduce parameters to study the effect of uncertainty and so their model is no longer fully deterministic. However, to truly capture the stochasticity of the system, consideration must be given to how the system evolves in time. The model must incorporate current assignment decisions having an impact on future assignments.

Using a similar integer programming approach to~\cite{Bilgin2012,Demeester2010,2016Van}, Abera et al. \cite{Abera2020} demonstrate how the PAS problem can be solved in a stochastic scenario. They consider first that a patient's LoS is a random variable that may take any value on some distribution, and is unknown at the time of decision making. They then consider that patient arrivals to the system are also random, meaning that the cost of making any given assignment is stochastic. That is, the cost of an assignment depends on how long a patient stays, and what types of patients will arrive next. To evaluate the assignment cost, and hence the optimal assignment, Abera et al. \cite{Abera2020} simulate the system several times to determine the expected cost over some planning horizon. The planning horizons used as examples in~\cite{Abera2020} are 14, 28 and 56 days. 

One shared result that is found in~\cite{Abera2020,Bilgin2012,Demeester2010,2016Van}, is that the integer linear program of any non-trivial size problem, cannot be solved explicitly. That is, it was not achievable to find the exact optimal solution. Demeester et al.~\cite{Demeester2010} found in their first example of a hospital with {\raise.17ex\hbox{$\scriptstyle\sim$}}$120$ beds, an optimal solution was not found after a whole week of computation, which is an unacceptable amount of time for assigning patients. To overcome this, various heuristics were required in all of the above work. The authors of the papers listed above all begin by using a neighbourhood search technique. That is, starting with some initial solution, and searching `neighbouring' solutions to see if there is an improvement. The exact methodology of how the neighbourhood is defined and searched, differs between authors. Demeester et al.~\cite{Demeester2010} use a Tabu search, Abera et al.~\cite{Abera2020} use simulated annealing, while Bilgin et al.~\cite{Bilgin2012} use a hyper-heuristic approach (heuristics for selecting heuristics). Regardless of the approach, it is clear that heuristics are a necessity for solving problems of any realistic size.

While the formulation of the PAS problem as an integer program is well researched, alternative formulations exist. Hulshof et al.~\cite{Hulshof2016} and Dai and Shi \cite{Dai2019} build models related to the PAS problem, based on Markovian Decision Processes (MDP). Such an alternative formulation allows for an exploration of possible decisions that a hospital manager can make. Hulshof et al.~\cite{Hulshof2016} explore resource assignment in hospitals and decision making on how many patients should be assigned to different queues. Dai and Shi~\cite{Dai2019} explore the use of patient overflow (allowing wards to `overflow' into one another when full) and take a `higher-level' approach to decision making. Rather than making assignment decisions for each patient, decisions are made on whether to allow patient overflow for each time period.

Similarly to the integer program approaches, the exact solutions of these MDPs cannot be obtained for large-size problems. Approximate Dynamic Programming (ADP) approaches are applied in~\cite{Dai2019,Hulshof2016} to overcome what is called the curse of dimensionality. The ADP approach, described by Powell in~\cite{Powell2009}, allows the cost in each state to be approximated by a set of features, which record some partial information about states, and corresponding weights. In addition to ADP, Dai and Shi \cite{Dai2019} utilize a least-squares temporal difference (LSTD) learning algorithm to solve the problem.

We now discuss the work of Heydar et al.~\cite{Heydar2021}, which is the basis of the model proposed here. Heydar et al.~\cite{Heydar2021} models the PAS as a continuous-time MDP, where the system is observed at every {\em individual} patient arrival or departure. The goal of the model is to identify the best possible patient assignment given information about the arriving patient and the occupancy of the hospital. Heterogeneity is included in the model by defining different patient classes and wards, and assuming that each patient class has an ideal ward, while other wards are unsuitable. In addition to assigning patients directly, Heydar et al.~\cite{Heydar2021} also consider transfers of patients between wards to best fit patients, so as to minimise their appropriately defined objective function. 

However, Heydar et al.~\cite{Heydar2021} note that transfers are positively correlated with patient mortality and LoS, as shown in~\cite{Hall2012}, meaning that transfers are given careful consideration and an additional cost in~\cite{Heydar2021}. We adopt this feature as a key component in our model.

Similarly to Abera et al.~\cite{Abera2020}, Heydar et al.~\cite{Heydar2021} evaluate total expected cost over a finite horizon, meaning limited consideration is given to the state of the system in the long run. In order to address this gap, we generalise the model used in Heydar et al.~\cite{Heydar2021} and extend it to an {\em infinite horizon problem}. In order to address the curse of dimensionality, we propose methodology based on ADP and the approximate policy iteration algorithm used by Dai and Shi in~\cite{Dai2019}.

Furthermore, we assume that patient assignments are made at {\em discrete time points} such as at the start of each day. This is arguably a more realistic approach involving a {\em group of patients} to be assigned, which also means that assignment decisions are far more complex, as multiple assignments must be made, and an order of assignments must be decided. 

Finally, we demonstrate that the algorithm used by Dai and Shi in~\cite{Dai2019} can be adapted to our model to find near-optimal solutions to the PAS problem.

The rest of this paper is structured as follows. In Section~\ref{sec:MDPmodel} we describe the key components of our model, such as state space, transition probabilities, decision variables, constraints, and cost variables. In Section~\ref{esc:ADP} we give the details of the Approximate Dynamic Programming approach, and illustrate the application of the theory through our numerical examples in Section~\ref{sec:numex}. This is followed by concluding remarks in Section~\ref{sec:conclusion}.
	
\section{Markov Decision Process}\label{sec:MDPmodel}
Let $\mathcal{I} = \{1, 2, \dots, I\}$ be the set of all patient types, where type $i\in\mathcal{I} $ may correspond to the medical needs of the patients, their age, gender, and other aspects of inclusive care, see e.g.~\cite{2012RCT}. We assume that type-$i$ patients arrive at a rate $\lambda_i$ {\em per day} on day $d$. Further, we assume that there are $K$ wards in the hospital, each with capacity $m_k$, $k\in\mathcal{K}=\{1,\ldots ,K\}$. The waiting room labelled $(K+1)$ has capacity $m_{K+1}$.

The set of wards that are suitable for patients type $i$ is denoted $\mathcal{K}(i)$, for some $\mathcal{K}(i)\subset \mathcal{K}$. We denote by $w(1,i)\in\mathcal{K}(i)$ the best ward for type-$i$ patient,  by $w(2,i)\in\mathcal{K}(i)$ the second-best ward for type-$i$ patient, and so on, and let $\mathcal{W}(i)=(w(1,i),\ldots, w(K,i))$ be the ordered sequence of wards in $\mathcal{K}(i)$ corresponding to type-$i$ patients. For example, if $\mathcal{K}(i)=\{1,2,3\}\subset \mathcal{K}=\{1,2,3,4,5\}$ and $\mathcal{W}(i)=(3,1,2)$, then this means that ward $3$ is the best ward for type-$i$ patient, ward $1$ is the second best, ward $2$ is the third best, and wards $4$ and $5$ are not suitable for type-$i$ patients.

In Section~\ref{MDP_1} below, we contruct a model  based on a suitable Markov Decision Process to find the optimal policy $\bpi^*$ from the set of all policies $\bpi =(a^{\bpi}(s))_{s\in\mathcal{S}}$ consisting of decisions $a^{\bpi}(s)$ taken whenever state $s$ in state space $\mathcal{S}$ is observed, so as to minimise the long-run expected cost per unit time,
\begin{eqnarray}
E^*
=\min_{\bpi} E^{\bpi}
=\min_{\bpi} \lim_{D\to\infty}
\left(\frac{1}{D}\ 
\mathbb{E}
\left(
\sum_{d=0}^{D-1}
C(S_d,a^{\bpi}(S_d)) 
\right)
\right),
\
\bpi^* 
=\arg\min_{\bpi} E^{\bpi},
\end{eqnarray} 
where  $S_d\in\mathcal{S}$ is a state of the system observed at the start of day $d$, and $C(S_d,a^{\bpi}(s))$ is the cost of a decision then taken assuming policy $\bpi$ is in place.

In practice, we find $\bpi^*$ by solving the Bellman's optimality equation for all states $s\in\mathcal{S}$,
\begin{eqnarray}\label{Bellman}
E^*+V(s) &=&  \min_{a \in \mathcal{A}(s)} \Bigg\{ C(s,a) +  \mathbb{E} \left( V(s' )\  \vert\  (s,a)  \right)   \Bigg\}
= \min_{a \in \mathcal{A}(s)} \Bigg\{ C(s,a) +  
\sum_{s' \in \mathcal{S}} \mathbb{P} \left( s^{\prime}\ \vert\  (s,a) \right) V(s') \Bigg\}, \nonumber\\
\end{eqnarray}
where $V(s)$ is the minimum long-run average cost given current state $s$, and $\mathcal{A}(s)$ is the set of all decisions that are possible when state $s$ is observed. We use Approximate Dynamic Programming methods~\cite{Powell_WhatYouShldKnow,Powell_ADPText,Powell_2016_persp} to solve Equation~\eqref{Bellman} in real-sized problems, in which the number of states and decisions are intractable, as noted in~\cite{Dai2019,gocgun2018simulation,Heydar2021,Hulshof2016}.

\subsection{Model} \label{MDP_1}
	
\noindent Consider a discrete-time Markov chain $\{S_d:d=0,1,\ldots\}$, where $S_d= \left( [N_{k,i}]_{\mathcal{K} \times \mathcal{I}}, [Q_i]_{1 \times \mathcal{I}} \right)\in\mathcal{S}$ is the state of the system at the start of day $d=0,1,2,\ldots$ recording $N_{k,i}$, the number of type-$i$ patients in ward $k$, and $Q_i$, the number of newly arrived type-$i$ patients who are yet to be assigned to the wards. The state space $\mathcal{S}$ of the process is given by,
\begin{eqnarray*}
\mathcal{S}&=&
\{
\left( [n_{k,i}]_{\mathcal{K} \times \mathcal{I}}, [q_i]_{1 \times \mathcal{I}} \right):
n_{k,i}\geq 0,q_i\geq 0,  \sum_{i=1}^{I} n_{k,i}\leq m_k,
\sum_{i=1}^{I}q_i
\leq m_{K+1}
\},
\end{eqnarray*}
with possibly additional constraints on the total number of accepted arrivals, as discussed below.
 
\subsection{Decisions}
Suppose that we observe state $s=\left( [n_{k,i}]_{\mathcal{K} \times \mathcal{I}}, [q_i]_{1 \times \mathcal{I}} \right)\in\mathcal{S}$ at the start of a day and choose a suitable decision $a\in\mathcal{A}(s)$ from the set of available decisions $\mathcal{A}(s)$, such that $a=(x_{k,i},y_{k,l,i})_{i\in\mathcal{I};k,\ell=1,\ldots K}$, where	
\begin{itemize}
	\item $x_{k,i}$ is the number of type-$i$ newly arrived patients to assign to ward $k$, and, 
	\item $y_{k,\ell,i}$ is the number of type-$i$ patients to be transferred from ward $k$ to ward~$\ell$.
\end{itemize}
	
Each decision should satisfy the following sets of constraints. First, all newly arrived type-$i$ patients have to be assigned to some ward, and so
\begin{eqnarray} \label{constraint_1}
	\sum_{k=1}^{K}  x_{k,i}= q_i,\   \forall i \in \mathcal{I}.
\end{eqnarray}
	
Next, the total number of patients in ward $k$, denoted $n_{k,\bullet,t}$, and given by,
\begin{eqnarray}
n_{k,\bullet}&=& \sum_{i=1}^{I}n_{k,i},\\
n_{k,i}
&=&
 n_{k,i} + x_{k,i}  + 
 \sum_{\substack{\ell = 1\\ \ell \ne k}}^{K}   \big( y_{\ell,k,i} - y_{k,\ell,i} \big),
\end{eqnarray}
where $n_{k,i}$ is the total number of type-$i$ patients in ward $k$, cannot exceed the capacity of the ward, and so,
\begin{eqnarray} \label{constraint_2}
	n_{k,\bullet} \le m_k,\ \forall k \in \mathcal{K}.
\end{eqnarray}
	
Further, we may only transfer the available patients, which gives,
\begin{eqnarray} \label{constraint_3}
	\sum_{\substack{\ell = 1 \\ \ell \ne k}}^{K} y_{k,\ell,i} \le n_{k,i},\ \forall k \in \mathcal{K},\ \forall i \in \mathcal{I},
\end{eqnarray}
and if a transfer occurs, it should be in one direction, that is,
\begin{eqnarray}\label{indicator_constraint}
I\{y_{k,\ell,i}\not= 0\}+I\{y_{\ell,k,i}\not= 0\}&\leq & 1,\
\forall k,\ell \in \mathcal{K},\ \forall i \in \mathcal{I},
\end{eqnarray}
where $I\{\cdot \}$ is an indicator function taking value $1$ if the statement in the brackets is true, and $0$ otherwise, which ensures that only one of $y_{k,\ell,i}$ and $y_{\ell,k,i}$ takes a nonzero value.

\subsection{Transition probabilities}
Now, given decision $a$ and current state $s \in \mathcal{S}$, we define the following key random variables:
\begin{itemize}
	\item $Z_{k,i}$, the number of type-$i$ patients who depart from ward $k$ during one day;
	\item $Z=\sum_{k=1}^{K} \sum_{i=1}^{I} Z_{k,i}$, the total number of departures;
	\item $B=\sum_{k=1}^{K}m_k-Z$ , the total number of available beds after departures;
	\item $Q= \sum_{i=1}^{I} Q_i$, the total number of arrivals;
\end{itemize}
which we use below to determine the transition probabilities of the process $\{S_d:d=0,1,\ldots\}$.	The distribution of these variables depends on $(s,a)$, as we discuss below.

Assuming state $s=\left( [n_{k,i}]_{\mathcal{K} \times \mathcal{I}}, [q_i]_{1 \times \mathcal{I}} \right)$ and decision $a = (x_{k,i},y_{k,l,i})_{i\in\mathcal{I};k\ell=1,\ldots K}$, the resulting post-decision state is $s^{(a)} = \left( [n^{(a)}_{k,i}]_{\mathcal{K} \times \mathcal{I}}\right) \equiv \left( [n^{(a)}_{k,i}]_{\mathcal{K} \times \mathcal{I}},[0]_{1 \times \mathcal{I}}\right)$ such that,
\begin{equation}\label{postA}
	n^{(a)}_{k,i}= 
n_{k,i} + x_{k,i}+ 
 \sum_{\substack{\ell = 1 \\ \ell \ne k}}^{K} \big(y_{\ell,k,i}-  y_{k,\ell,i} \big)
	,\ \forall k \in \mathcal{K},\ \forall i \in \mathcal{I}.
\end{equation}

Next, the process transitions from state $s^{(a)}$ on day $d$ to some state $s' = \left( [n'_{k,i}]_{\mathcal{K} \times \mathcal{I}}, [q'_i]_{1 \times \mathcal{I}} \right)$ on day $(d+1)$, for $d=0,1,2,\ldots$, with probability given by
\begin{eqnarray}\label{transition_A}
\mathbb{P} \left( s^{\prime}\ \vert\  (s,a) \right)
&=&
\mathbb{P}
\left(
(Z_{k,i}=z_{k,i})_{k=1,\ldots,K;i=1,\ldots,I}
 \mbox{ and } (Q_i=q'_i)_{i=1,\ldots,I}
\right)
\nonumber\\
&=&
\mathbb{P}
\left(
(Z_{k,i}=z_{k,i})_{k=1,\ldots,K;i=1,\ldots,I}
\right)
\times \mathbb{P}
\left(
(Q_i=q'_i)_{i=1,\ldots,I}\ | \
(Z_{k,i}=z_{k,i})_{k=1,\ldots,K;i=1,\ldots,I}
\right)
\nonumber\\
&=&
\left(
\prod_{k=1}^{K} \prod_{i=1}^{I}{\mathbb{P}(Z_{k,i}=z_{k,i})}
\right)
\times 
\mathbb{P}
\left(
(Q_i=q'_i)_{i=1,\ldots,I}\ | \
(Z_{k,i}=z_{k,i})_{k=1,\ldots,K;i=1,\ldots,I}
\right)
\end{eqnarray}
when the condition
\begin{equation}
n'_{k,i}=n^{(a)}_{k,i}-z_{k,i}
\end{equation} 
is met for all $k\in\mathcal{K}$, $i\in\mathcal{I}$; and $\mathbb{P} \left( s^{\prime}\ \vert\  (s,a) \right)=0$ otherwise.

We apply conditional probabilities 
$\mathbb{P}
\left(
(Q_i=q'_i)_{i=1,\ldots,I}\ | \
(Z_{k,i}=z_{k,i})_{k=1,\ldots,K;i=1,\ldots,I}
\right)$ in~\eqref{transition_A} since the number of accepted type-$i$ arrivals may depend on the number of available beds after departures, and consider the following alternative modelling approaches:
\begin{itemize}
	\item Suppose that there is no restriction on the total number of arrivals, as in Dai and Shi~\cite{Dai2019}. Then $Q_i$ follows Poisson distribution, $Q_i\sim Poi(\lambda_i)$, and so
	\begin{eqnarray}
		\mathbb{P}
		\left(
		(Q_i=q'_i)_{i=1,\ldots,I}\ | \
		(Z_{k,i}=z_{k,i})_{k=1,\ldots,K;i=1,\ldots,I}
		\right)
	&=&
	\prod_{i=1}^I 
	\mathbb{P}(Q_i = q_i)
	=
	\prod_{i=1}^I 
	\dfrac{(\lambda_i)^{q_i}
		e^{-\lambda_i}}{(q_i)!}.
	\end{eqnarray}	
In this approach, arriving patients are not lost to the system, however the number of patients still waiting to be assigned and their total waiting times could be very large (and exceed many days). This may not be a realistic assumption since in practice the capacity of the system is limited (due to the availability of beds and staffing) and there are limits on the accepted maximum waiting times for various patient types (e.g. $24$-hour limit for some types of 
emergency patients).
 	
\item Suppose that the total number of arrivals $Q$ may not exceed the total number of available beds. Then, with $b=\sum_{k=1}^K m_k- \sum_{k=1}^K \sum_{i=1}^I z_{k,i}$ recording the total number of available beds after the departures, we have,
\begin{eqnarray}
\mathbb{P}(Q = q \ | \ 
B=b)
&=&\dfrac{(\lambda)^{q}e^{-\lambda}}{(q)!},
\       q=0,1,2\ldots ,b-1,
\\
\mathbb{P}(Q = b \ | \ 
B=b)&=& 1- \sum_{q=0}^{b-1}\dfrac{(\lambda)^{q}e^{-\lambda}}{(q)!},
\end{eqnarray}
and, with $q=\sum_{i=1}^I q_i$,
\begin{eqnarray}
\lefteqn{
	\mathbb{P}
	\left(
	(Q_i=q'_i)_{i=1,\ldots,I}\ | \
	(Z_{k,i}=z_{k,i})_{k=1,\ldots,K;i=1,\ldots,I}
	\right)
}
\nonumber\\
&=&
	\mathbb{P}
	\left(
	(Q_i=q'_i)_{i=1,\ldots,I}\ | \
	B=b
	\right)
\nonumber\\
&=&
\mathbb{P}
\left(
(Q_i=q'_i)_{i=1,\ldots,I}\ | \
Q=q,B=b
\right)
\times \mathbb{P}(Q = q \ | \ 
B=b)
\nonumber\\
&=&
\mathbb{P}
\left(
(Q_i=q'_i)_{i=1,\ldots,I}\ | \
Q=q
\right)
\times \mathbb{P}(Q = q \ | \ 
B=b)
\nonumber\\
&=&
\frac{q!}{\prod_{i=1}^I q_i!}
\times
\prod_{i=1}^I 
\left( \frac{\lambda_i}{\lambda} \right)^{q_i}
\times \mathbb{P}(Q = q \ | \ 
B=b)
,
\end{eqnarray}
since the conditional distribution of $(Q_i=q'_i)_{i=1,\ldots,I}$ given $Q=q$ is multinomial, where $\lambda_i/\lambda$ is the probability that an arriving patient is of type $i$. In this approach, patients arriving to a system with no free beds, are redirected to other health systems. Such approach may be more realistic and can be used to determine the rate of patients that are redirected to help decide a suitable size of the hospital system.

\end{itemize}

To model random variables $Z_{k,i}$, we assume that departures are independent of one another and so $Z_{k,i}$ follows Binomial distribution, $Z_{k,i}\sim Bin(n_{k,i},p_{k,i})$, with	
\begin{eqnarray}
\mathbb{P}(Z_{k,i}=z)&=& \binom{n_{k,i}}{z} 
(p_{k,i})^z (1-p_{k,i})^{n_{k,i}-z},\ z=0,1,\ldots,n_{k,i},\\
	p_{k,i}&=&\mathbb{P}(RLOS_{k,i}\leq 1 ),
\end{eqnarray}
where $RLOS_{k,i}$ is a random variable recording the remaining length of stay of type-$i$ patient that is in ward $k$ at the start of the day. That is, for each type-$i$ patient in ward $k$, we perform a Bernoulli trial with probability of success $p_{k,i}$ to determine if the patient leaves the system on the following day or not.

\subsection{Costs}	
There is an immediate cost $C(s,a)$ associated with decision $a$ given current state $s \in \mathcal{S}$, which may include costs of assignment, transfer, and patients being in nonprimary wards, defined as follows.	
\begin{itemize}
	\item Assignment cost: $\sum_{k=1}^{K} \sum_{i=1}^{I} x_{k,i} \times  c^{(\sigma)}_{k,i}$.

	\item Transfer cost: $\sum_{k=1}^{K} \sum_{\ell=1}^{K} \sum_{i=1}^{I}  y_{k,\ell,i}\times  c^{(t)}_{k,\ell,i} $.
		
	\item Penalty cost for being in a nonprimary ward:  $\sum_{k=1}^{K} \sum_{i=1}^{I} n_{k,i} \times  c^{(p)}_{k,i}  $. 
\end{itemize}
The total immediate cost of decision $a=(x_{k,i},y_{k,l,i})_{i\in\mathcal{I};k,\ell=1,\ldots K}$ given state $s=\left( [n_{k,i}]_{\mathcal{K} \times \mathcal{I}}, [q_i]_{1 \times \mathcal{I}} \right)$ is then,
\begin{eqnarray} \label{DailyCost_I}
	C(s,a) &=&
	\sum_{k=1}^{K} \sum_{i=1}^{I}  \bigg\{ x_{k,i} \times c^{(\sigma)}_{k,i}  + \sum_{\substack{\ell=1 \\ \ell \ne k}}^{K} y_{k,\ell,i}\times c^{(t)}_{k,\ell,i} + n^{(a)}_{k,i}\times c^{(p)}_{k,i}  \bigg\}
	\nonumber\\
	&=&
\sum_{k=1}^{K} \sum_{i=1}^{I}  \bigg\{ x_{k,i} \times c^{(\sigma)}_{k,i}  
+ \sum_{\substack{\ell=1 \\ \ell \ne k}}^{K} 
	y_{k,\ell,i}
	\times c^{(t)}_{k,\ell,i}
	+ \Big( n_{k,i} +  
	x_{k,i} + \sum_{\substack{\ell=1 \\ \ell \ne k}}^{K} (y_{\ell,k,i} - y_{k,\ell,i}) 
	\Big)
	\times c^{(p)}_{k,i} \bigg\}	.
\end{eqnarray}

By above, for all $s=\left( [n_{k,i}]_{\mathcal{K} \times \mathcal{I}}, [q_i]_{1 \times \mathcal{I}}\right)$, Eq.~\eqref{Bellman} can be written as,
\begin{eqnarray}\label{Bellman_A}
E^*+V\left( s \right)  
&=&\min_{a \in \mathcal{A}(s)} \Bigg\{ C(s,a) 
+  \sum_{s'=\left( [n^{'}_{k,i}]_{\mathcal{K} \times \mathcal{I}}, [q'_i]_{1 \times \mathcal{I}}\right) \in \mathcal{S}}\mathbb{P} \left( s^{\prime}\ \vert\  (s,a) \right) V\left( s' \right) \Bigg\}, \nonumber\\
\end{eqnarray}
where $\mathbb{P} \left( s^{\prime}\ \vert\  (s,a) \right)$ and $C\left( s, a \right) $ are given by Eqs.~\eqref{transition_A} and~\eqref{DailyCost_I}, respectively.

\section{Approximate Dynamic Programming approach}\label{esc:ADP}

Policy Iteration is a standard method of solving Equation~\eqref{Bellman} when the size of the state space $\mathcal{S}$ is not too large, see Algorithm~\ref{Al_pol_it} in Appendix~\ref{sec:algorithms}, also refer to Puterman~\cite{1994_Put} for further details.

However, Policy Iteration is not suitable here, and so we apply Approximate Dynamic Programming methods introduced by Powell in~\cite{Powell_WhatYouShldKnow,Powell_ADPText,Powell_2016_persp} to address the curses of dimensionality in Eq.~\eqref{Bellman}, as follows. 
\begin{enumerate}
	\item First, similar to Heydar et al.~\cite{Heydar2021}, we apply basis functions $\phi_f(s)$, $f\in\mathcal{F}$, to record some suitable information about states $s$, referred to as state {\em features}, and then the approximation,
	\begin{eqnarray}
	V(s)\approx \sum_{ f \in \mathcal{F}} \phi_f(s) \theta^{(f)}_n 
	=
	\bphi(s)^T \btheta_n,
	\end{eqnarray}
	where $\bphi(s)=[\phi_f(s)]_{f\in \mathcal{F}}$ is the vector of features, and $\btheta_n=[\theta^{(f)}_n]_{f\in \mathcal{F}}$ is the vector of corresponding weights $\theta^{(f)}_n$, evaluated at the $n$-th iteration of an algorithm. 
	
	\item Next, we apply the Approximate Policy Iteration presented by Dai and Shi in~\cite{Dai2019}, in which the weights $\theta^{(f)}_n$ are recursively updated in an iteration that is repeated $N$ times. Each iteration $n$ involves a simulation of $M$ states for some large $M$. We summarise this approach in Algorithms~\ref{Al1}--\ref{Al1E} in Appendix~\ref{sec:algorithms}. We note that the Markov chains applied in our models are irreducible and positive recurrent, and so the algorithms are guaranteed to converge as $M\to\infty$ ~\cite{Dai2019}.
\end{enumerate}

Consider our Model in Section~\ref{MDP_1}. Suppose that the vector of features of state $s=\left( [n_{k,i}]_{\mathcal{K} \times \mathcal{I}}, [q_i]_{1 \times \mathcal{I}} \right)$ is the vector $\vec{s}$, and so a vector containing all information about state $s$. We note that to compute the expression in Line 6 of Algorithms~\ref{Al1}  in Appendix~\ref{sec:algorithms} and in Line 4 in Algorithm~\ref{Al1E} in Appendix~\ref{sec:algorithms}, it is then convenient to apply the following equivalence. Assuming $s=\left( [n_{k,i}]_{\mathcal{K} \times \mathcal{I}}, [q_i]_{1 \times \mathcal{I}} \right)$, and with $b=\sum_{k=1}^K m_k-\sum_{k=1}^{K}\sum_{i=1}^{I}
n^{(a)}_{k,i}$, we have
\begin{eqnarray}
\sum_{s^{'}\in\mathcal{S}} \mathbb{P} \left( s^{\prime}\ \vert\  (s,a) \right) \bphi(s^{'}) \btheta
&=&
\mathbb{E}(\bphi(s^{'}) \btheta \ \vert\  (s,a))\\
&=&
\sum_{k=1}^{K} \sum_{i=1}^{I}\mathbb{E}(N_{k,i}^{'})\theta_{k,i}
+
 \sum_{i=1}^{I}\mathbb{E}(Q_i^{'})\theta_i
\\
&=&
 \sum_{k=1}^{K}\sum_{i=1}^{I}
 (n^{(a)}_{k,i}-\mathbb{E}(Z_{k,i}))\theta_{k,i}
 +
\sum_{i=1}^{I}\mathbb{E}\left(\mathbb{E}(Q_i^{'}
\ |\ Q\ )\right)\theta_i
\\
&=&
\sum_{k=1}^{K}\sum_{i=1}^{I}
n^{(a)}_{k,i}(1-p_{k,i})\theta_{k,i}
+
\sum_{i=1}^{I}
\frac{\lambda_i}{\lambda}
\times  \mathbb{E}(Q)\theta_i
,
\label{eq:expected_features}
\end{eqnarray}
where $n^{(a)}_{k,i}$ is given by Eq.~\eqref{postA}, $\mathbb{E}(Z_{k,i})=n^{(a)}_{k,i}p_{k,i}$ (binomial mean),  $N_{k,i}^{'}$ is a random variable recording the number of type-$i$ patients in ward $k$ and $Q_i^{'}$ a random variable recording the number of type-$i$ patients waiting to be assigned when state $s^{'}$ is observed, and $Q=\sum_{i=1}^{I} Q_i^{'}$ is the total number of arrivals.

If the number of arrivals is unrestricted with $m_{K+1}=\infty$, then $\mathbb{E}(Q)=\lambda$. Alternatively, if the number of arrivals is restricted so that we only allow arrivals that can be assigned, which can be written as,
	\begin{equation}
	\sum_{k=1}^K\sum_{i=1}^I n_{k,i}^{'}+\sum_{i=1}^I q_i^{'} \leq \sum_{k=1}^K m_k =m, 
	\end{equation}
then we have the following approach for evaluating $\mathbb{E}(Q)$.

Let $Z$ be the random variable recording the total number of departures,
\begin{eqnarray}
Z&=&\sum_{k=1}^{K}\sum_{i=1}^{I} Z_{k,i},
\end{eqnarray}
taking values $z=0,\ldots,z_{max}$, with $z_{max}=\sum_{k=1}^K \sum_{i=1}^I n_{k,i}^{(a)}$. Let $N$ be the number of available beds for the new arrivals, given by
\begin{eqnarray}
N&=&m-\sum_{k=1}^{K}\sum_{i=1}^{I} N_{k,i}^{'}=m-z_{max}+Z,
\end{eqnarray}
taking values $n=m-z_{max},\ldots,m$. Therefore,
\begin{eqnarray*}
\mathbb{E}
\left(
Q
\ |\ Z=z
\right)&=&
\mathbb{E}
\left(
Q
\ |\ N=m-z_{max}+z
\right)
\\
&=&\sum_{k=0}^{m-z_{max}+z}
k\times\frac{\lambda^k}{k!}e^{-\lambda}
+
\sum_{k=m-z_{max}+z+1}^{\infty}
(m-z_{max}+z)
\times
\frac{\lambda^k}{k!}e^{-\lambda}
\\
&=&
\sum_{k=0}^{m-z_{max}+z}
k\times\frac{\lambda^k}{k!}e^{-\lambda}
+
(m-z_{max}+z)
\left(
1-
\sum_{k=0}^{m-z_{max}+z}
\frac{\lambda^k}{k!}e^{-\lambda}
\right)
\\
&=&
\sum_{k=0}^{m-z_{max}+z-1}
\lambda
\times
\frac{\lambda^k}{k!}e^{-\lambda}
+
(m-z_{max}+z)
\left(
1-
\sum_{k=0}^{m-z_{max}+z}
\frac{\lambda^k}{k!}e^{-\lambda}
\right)
\end{eqnarray*}
and so
\begin{eqnarray}\label{eq:lambdatilde1}
\mathbb{E}(Q)
&=&
\sum_{z=0}^{z_{max}}
\mathbb{E}
\left(
Q
\ |\ Z=z
\right)
\mathbb{P}(Z=z)
\nonumber\\
&=&
\sum_{z=0}^{z_{max}}
\left(
\sum_{k=0}^{m-z_{max}+z}
k
\times
\frac{\lambda^k}{k!}e^{-\lambda}
+
(m-z_{max}+z)
\left(
1-
\sum_{k=0}^{m-z_{max}+z}
\frac{\lambda^k}{k!}e^{-\lambda}
\right)
\right)
\mathbb{P}(Z=z)
\nonumber\\
&=&
\sum_{z=0}^{z_{max}}\mathbb{P}(Z=z)
(m-z_{max}+z)
\nonumber\\
&&
+
\sum_{z=0}^{z_{max}}
\left(
\sum_{k=0}^{m-z_{max}+z-1}
\lambda
\times
\frac{\lambda^k}{k!}e^{-\lambda}
-
(m-z_{max}+z)
\left(
\sum_{k=0}^{m-z_{max}+z}
\frac{\lambda^k}{k!}e^{-\lambda}
\right)
\right)
\mathbb{P}(Z=z)
\nonumber\\
&=&
m-z_{max}+\mathbb{E}(Z)
\nonumber\\
&&
+
\sum_{z=0}^{z_{max}}
\left(
\lambda
-(m-z_{max}+z)
\right)
\left(
\sum_{k=0}^{m-z_{max}+z-1}
\frac{\lambda^k}{k!}
e^{-\lambda}
\right)
\mathbb{P}(Z=z)
\nonumber\\
&&
-
\sum_{z=0}^{z_{max}}
(m-z_{max}+z)
\left(
\frac{\lambda^{m-z_{max}+z}}{(m-z_{max}+z)!}
e^{-\lambda}
\right)
\mathbb{P}(Z=z)
,
\end{eqnarray}
where $\mathbb{E}(Z)=\sum_{k=1}^{K}\sum_{i=1}^{I} n^{(a)}_{k,i}p_{k,i}$.

Further, in order to compute $\mathbb{P}(Z=z)$ in~\eqref{eq:lambdatilde1}, we note that $Z=\sum_{k=1}^K \sum_{i=1}^I Z_{k,i} I(n_{k,i}^{(a)}\not=0)$ is a sum of independent binomial random variables $Z_{k,i}\sim Bin(n_{k,i}^{(a)},p_{k,i})$ such that $n_{k,i}^{(a)}\not=0$.

A method for computing the distribution of a sum of independent binomial random variables was discussed by Butler and Stephens in~\cite{2016BS}. Below, we suggest an alternative method, which relies on simple matrix multiplications and standard theory of Markov chains.

Without loss of generality, suppose that $Z=\sum_{i=1}^M Z_i$ is a random variable such that $Z_i\sim Bin(n_i,p_i)$ are independent Binomial random variables with some parameters $n_i>0$ and $0<p_i<1$, for $i=1,\ldots M$.

To compute $\mathbb{P}(Z=z)$, it is convenient to construct a discrete-time Markov chain corresponding to the Bernoulli trials $n_1,n_2,\ldots,n_M$ such that a success at each trial results in the chain moving one step to the right, and failure results in the chain remaining at the original state. Then, the distribution of the chain after all $n_1+\ldots+n_M$ trials will give the distribution of $Z$.

So, we consider a discrete-time Markov chain $\{(J(t)):t= 0,1,2,\ldots,z_{max}\}$ terminating at time $z_{max}=\sum_{i=1}^M n_i$, with state space $\mathcal{S}=\{0,\ldots,z_{max}\}$ and an initial state $J(0)=0$, which evolves as follows. We perform the first $n_1$ Bernoulli trials at times $t=1,\ldots,n_1$ and let $\mathbb{P}(J(t)=J(t-1)+1)=p_1$, $\mathbb{P}(J(t)=J(t-1))=1-p_1$. Then, the distribution of the chain after the first $n_1$ trials, and so at time $t=n_1$, is given by
\begin{eqnarray}
\balpha(n_1)&=&
\left[
\balpha(0)\quad {\bf 0}_{1\times n_1}
\right]
{\bf P}^{(1)},
\end{eqnarray}
where $\balpha(0)=[\alpha(0)_j]_{j=0}$ with $\alpha(0)_0=1$, and ${\bf P}^{[1]}=[P_{j,j'}^{[1]}]_{j,j'=0,1,\ldots,n_1}$ such that $P_{j,j'}^{[1]}=\mathbb{P}(Z_1=j'-j)$.

Next, we repeat this and perform $n_i$ Bernoulli trials at times $t=n_1+\ldots+n_{i-1}+1,\ldots,n_1+\ldots+n_i$ and let $\mathbb{P}(J(t)=J(t-1)+1)=p_i$, $\mathbb{P}(J(t)=J(t-1))=1-p_i$, for each $i=2,\ldots,M$. Then, the distribution after additional $n_i$ trials is given by the recursion
\begin{eqnarray}
\balpha(n_1+\ldots+n_i)&=&
\left[
\balpha(n_1+\ldots+n_{i-1})
\quad
{\bf 0}_{1\times n_i}
\right]
{\bf P}^{(i)},
\end{eqnarray}
where ${\bf P}^{[i]}=[P_{j,j'}^{[i]}]_{j,j'=0,1,\ldots,n_1+\ldots+n_i}$,  $P_{j,j'}^{[i]}=\mathbb{P}(Z_1=j'-j)$. To compute ${\bf P}^{(i)}$ we apply the formula
\begin{eqnarray}\label{eq:Piformula}
{\bf P}^{(i)}&=&\left({\bf A}^{(i)}\right)^{n_i},
\end{eqnarray}
where
\begin{eqnarray*}
	{\bf A}^{(i)}&=&
	[A^{(i)}]_{j,j'=0,1,\ldots,n_1+\ldots+n_i}	
	=
	\left[
	\begin{array}{llllll}
		1-p_i&p_i&0&\ldots&\ldots&0\\
		0&1-p_i&p_i&\ldots&\ldots&0\\
		\vdots&\vdots&\vdots&\vdots&\vdots&\vdots\\
		0&\ldots&0&\ldots&1-p_i&p_i
	\end{array}
	\right].
\end{eqnarray*}
It follows that, for $j=0,1,\ldots,z_{max}$,
\begin{eqnarray}
\mathbb{P}(Z=j)&=&
[\balpha(n_1+\ldots+n_M)]_j
=
\left[
\balpha(n_1+\ldots+n_{M-1})
\quad
{\bf 0}_{1\times n_M}
\right]
{\bf P}^{(M)}.
\end{eqnarray}

This method for computing the probabilities $\mathbb{P}(Z=z)$ for a sum $Z=\sum_{i=1}^M Z_i$ of independent nonnegative random variables $Z_i$ that take values in finite sets, can be applied for any desired discrete distributions of $Z_i$, by replacing~\eqref{eq:Piformula} with a suitable formula for ${\bf P}^{(i)}$.

Alternatively, $\mathbb{P}(Z=z)$ can be computed by numerically inverting the probability generating function $G_Z(s)$ of the random variable $Z$, which here is given by
\begin{eqnarray}
G_Z(s)&=&\prod_{i=1}^M
G_{Z_i}(s)
=
\prod_{i=1}^M (1-p_i+sp_i)^{n_i},
\end{eqnarray}
using a suitable inversion algorithm, for example see Abate and Whitt~\cite{ABATE1992245}.

\section{Numerical examples}\label{sec:numex}

Here, we construct examples to illustrate the theory and the application of algorithms discussed above. First, in Example~\ref{ex1}, we construct a Markov model with a small state space, and demonstrate that the approximate solution converges to the exact optimal solution (which we obtained by applying standard dynamic programming methods). Next, in Example~\ref{ex2}, we construct a realistically-sized Markov model with assumptions driven by practical considerations driven by the conditions of real-world hospitals, and demonstrate the application potential of our methodology.

\subsection{Small-sized example}\label{sec:small_size}

We consider the following simple example to illustrate the application of our Model in Section~\ref{MDP_1}. 	As the size of the state space $\mathcal{S}$ in this example is small enough to allow for the application of standard dynamic programming methods, we apply Algorithm~\ref{Al_pol_it} in Appendix~\ref{sec:algorithms} to find the exact solution and then compare it with the approximation obtained using Algorithms~\ref{Al1}--\ref{Al1E}.

\begin{exa}\label{ex1}

	Consider a healthcare facility with $K = 2$ wards, with capacity $m_k = 1$ for each ward $k$, and $I = 2$ patient types. Assuming that the maximum capacity of the waiting area is $m_3=2$, that is $\sum_{i=1}^I q_i\leq m_3$, and that the arrivals are only accepted if there is available capacity in the system according to 
	\begin{equation}
	\sum_{k=1}^K\sum_{i=1}^I n_{k,i}+\sum_{i=1}^I q_i \leq \sum_{k=1}^K m_k =m, 
	\end{equation}we define the state space of the system as follows:
	\begin{eqnarray*}
		\mathcal{S} &=& 
		\left\{
		\left(
		\begin{bmatrix} n_{1,1}\ n_{1,2}\\n_{2,1}\ n_{2,2} \end{bmatrix}, 
		\begin{bmatrix} q_1\ q_2 \end{bmatrix}
		\right)
		:
		\sum_{j=1}^2 n_{1,j} \leq 1,
		\sum_{j=1}^2 n_{2,j} \leq 1,
 \sum_{i=1}^2 + \sum_{j=1}^2 n_{i,j}+\sum_{i=1}^2 q_i\leq 2
		\right\}
		\\
		&= &
					\left\{
					1,2,3,\ldots ,22
					\right\},
			\end{eqnarray*}
		with
	\begin{eqnarray*}	
		&&
		1\equiv 
		\left(
		\begin{bmatrix} 0\ 0 \\ 0\ 0 \end{bmatrix},
		\begin{bmatrix} 0\ 0 \end{bmatrix} 
		\right),
		2\equiv \left(
		\begin{bmatrix} 0\ 0 \\ 0\ 0 \end{bmatrix},
		\begin{bmatrix} 1\ 0 \end{bmatrix} 
		\right),
		3\equiv 
		\left(
		\begin{bmatrix} 0\ 0 \\ 0\ 0 \end{bmatrix},
		\begin{bmatrix} 0\ 1
		\end{bmatrix} 
		\right),
		4\equiv 
		\left(
		\begin{bmatrix} 0\ 0 \\ 0\ 0 \end{bmatrix},
		\begin{bmatrix} 2\ 0 \end{bmatrix} 
		\right),
		\\
		&&
		5\equiv 
		\left(
		\begin{bmatrix} 0\ 0 \\ 0\ 0 \end{bmatrix},
		\begin{bmatrix} 1\ 1 \end{bmatrix} 
		\right),
		6\equiv 
		\left(
		\begin{bmatrix} 0\ 0 \\ 0\ 0 \end{bmatrix},
		\begin{bmatrix} 0\ 2 \end{bmatrix} 
		\right),
		7\equiv 
		\left(
		\begin{bmatrix} 1\ 0 \\ 0\ 0 \end{bmatrix},
		\begin{bmatrix} 0\ 0 \end{bmatrix} 
		\right),
		8\equiv 
		\left(
		\begin{bmatrix} 1\ 0 \\ 0\ 0 \end{bmatrix},
		\begin{bmatrix} 1\ 0 \end{bmatrix} 
		\right),
		\\
		&&
		9\equiv 
		\left(
		\begin{bmatrix} 1\ 0 \\ 0\ 0 \end{bmatrix},
		\begin{bmatrix} 0\ 1 \end{bmatrix} 
		\right),
		10\equiv 
		\left(
		\begin{bmatrix} 0\ 1 \\ 0\ 0 \end{bmatrix},
		\begin{bmatrix} 0\ 0 \end{bmatrix} 
		\right),
		11\equiv 
		\left(
		\begin{bmatrix} 0\ 1 \\ 0\ 0 \end{bmatrix},
		\begin{bmatrix} 1\ 0 \end{bmatrix} 
		\right),
		12\equiv 
		\left(
		\begin{bmatrix} 0\ 1 \\ 0\ 0 \end{bmatrix},
		\begin{bmatrix} 0\ 1 \end{bmatrix} 
		\right),
		\\
		&&
		13\equiv 
		\left(
		\begin{bmatrix} 0\ 0 \\ 1\ 0 \end{bmatrix},
		\begin{bmatrix} 0\ 0 \end{bmatrix} 
		\right),
		14\equiv 
		\left(
		\begin{bmatrix} 0\ 0 \\ 1\ 0 \end{bmatrix},
		\begin{bmatrix} 1\ 0 \end{bmatrix} 
		\right),
		15\equiv 
		\left(
		\begin{bmatrix} 0\ 0 \\ 1\ 0 \end{bmatrix},
		\begin{bmatrix} 0\ 1 \end{bmatrix} 
		\right),
		16\equiv 
		\left(
		\begin{bmatrix} 0\ 0 \\ 0\ 1 \end{bmatrix},
		\begin{bmatrix} 0\ 0 \end{bmatrix} 
		\right),
		\\
		&&
		17\equiv 
		\left(
		\begin{bmatrix} 0\ 0 \\ 0\ 1 \end{bmatrix},
		\begin{bmatrix} 1\ 0 \end{bmatrix} 
		\right),
		18\equiv 
		\left(
		\begin{bmatrix} 0\ 0 \\ 0\ 1 \end{bmatrix},
		\begin{bmatrix} 0\ 1 \end{bmatrix} 
		\right),
		19\equiv 
		\left(
		\begin{bmatrix} 1\ 0 \\ 1\ 0 \end{bmatrix},
		\begin{bmatrix} 0\ 0 \end{bmatrix} 
		\right),
		20\equiv 
		\left(
		\begin{bmatrix} 0\ 1 \\ 1\ 0 \end{bmatrix},
		\begin{bmatrix} 0\ 0 \end{bmatrix} 
		\right),
		\\
		&&
	21\equiv 
		\left(
		\begin{bmatrix} 1\ 0 \\ 0\ 1 \end{bmatrix},
		\begin{bmatrix} 0\ 0 \end{bmatrix} 
		\right),
		22\equiv 
		\left(
		\begin{bmatrix} 0\ 1 \\ 0\ 1 \end{bmatrix},
		\begin{bmatrix} 0\ 0 \end{bmatrix} 
		\right).
	\end{eqnarray*}

Further, assume that ward $i$ is the preferred ward for patient type $i$, for $i=1,2$, which is reflected by the expected values of the LoS, which are lower when patient type $i$ stays in ward $i$ for the whole duration of their stay. To model this, we assume that probabilities $p_{k,i}=\mathbb{P}(RLOS_{k,i}\leq 1 )$ that patient type $i$ that is in ward $k$ today, leaves the hospital the next day, are given by,
\begin{eqnarray}
{\bf p}=[p_{k,i}]_{\mathcal{K} \times \mathcal{I}}&=&
\left[
\begin{array}{cc}
1/5 & 1/4 \\ 1/10& 1/3 
\end{array}
\right],
\end{eqnarray}
and so the expected LoS for patient type $i=1$ is $E(LoS_{k,i})=5$ days if they stay in ward $k=1$, or $E(LoS_{k,i})=10$ if they stay in ward $k=2$, and some number between these two if the patient is transferred between the wards during their stay at the hospital. For patient type $i=2$ we have $E(LoS_{k,i})=3$ days if they stay in ward $k=2$, or $E(LoS_{k,i})=4$ if they stay in ward $k=1$, and some number between these two if the patient is transferred between the wards during their stay at the hospital.

\end{exa}

We consider two decisions, 
\begin{itemize}
	\item $a=1$, assign arrived patients to their best available wards without transferring the patients between the wards, and
	\item $a=2$, assign arrived patients to their best available wards and allow transferring patients if required.
\end{itemize}
As example, given state
$
s=15\equiv \left(
\begin{bmatrix} 0\ 0 \\ 1\ 0 \end{bmatrix},
\begin{bmatrix} 0\ 1 \end{bmatrix} 
\right),
$
decision $a=1$ with transform it into post decision state 
$
s^{(a)}= \left(
\begin{bmatrix} 0\ 1 \\ 1\ 0 \end{bmatrix}
\right),
$
 and then the process will move to state 
 $
 s=10\equiv \left(
 \begin{bmatrix} 0\ 1 \\ 0\ 0 \end{bmatrix},
 \begin{bmatrix} 0\ 0 \end{bmatrix} 
 \right)
 $ if patient type $1$ departs from ward $2$, or to state 
 $
 s=13\equiv \left(
 \begin{bmatrix} 0\ 0 \\ 1\ 0 \end{bmatrix},
 \begin{bmatrix} 0\ 0 \end{bmatrix} 
 \right)
 $ if patient type $2$ departs from ward $1$. Alternatively, decision $a=2$ will transform $s=15$ into post-decision state 
 $
 s^{(a)}= \left(
 \begin{bmatrix} 1\ 0 \\ 0\ 1 \end{bmatrix}
 \right),
 $
 and then the process will move to state 
 $
 s=16\equiv \left(
 \begin{bmatrix} 0\ 0 \\ 0\ 1 \end{bmatrix},
 \begin{bmatrix} 0\ 0 \end{bmatrix} 
 \right),
 $ if patient type $1$ departs from ward $1$, or to state 
  $
 s=7\equiv \left(
 \begin{bmatrix} 1\ 0 \\ 0\ 0 \end{bmatrix},
 \begin{bmatrix} 0\ 0 \end{bmatrix} 
 \right),
 $ if patient type $2$ departs from ward $2$. No arrivals to waiting area are permitted if they cannot be assigned, that is, when the system is already full.

We assume the following cost parameters:
\begin{itemize}
	\item Assignment cost to ward $k$ per type-$i$ patient: $c^{(\sigma)}_{k,i}=c^{(\sigma)}=1$ for all $k, i$.

	\item Transfer cost from ward $k$ to ward $\ell$ per type-$i$ patient: $ c^{(t)}_{k,\ell,i}=c^{(t)}=1.1$ for all $k,\ell, i$.
	
	\item Penalty cost for being in a nonprimary ward $k$ per type-$i$ patient:  $c^{(p)}_{k,i}=c^{(p)}=0.2$ for all $k, i$. 
\end{itemize}

We apply policy iteration summarised in Algorithm~\ref{Al_pol_it} in Appendix~\ref{sec:algorithms} to find the exact solution. To do so, we first write explicit expressions for the probability matrices ${\bf P}^{(a)}$ and cost vectors ${\bf C}^{(a)}$ for $a=1,2$. The details of this derivation are given in Appendix~\ref{sec:ex1matrices}. 

From policy iteration we identify the optimal policy in our states of interest to be $a=2$, `assign the waiting patients to their most preferred ward by transferring patients currently in the ward'. 
We determine that this optimal policy invokes a long run average cost per stage of $E^*=0.4098$ units.

We now present the results of the approximate policy iteration algorithm, summarised in Algorithms~\ref{Al1} and~\ref{Al1E}. For simulations, we assumed the starting value of each weight $\btheta_0=\boldsymbol{1}\times 10^{-4}$.

\begin{figure}[H]
	\begin{tabular}{cccc}
		{\includegraphics[scale=0.75]{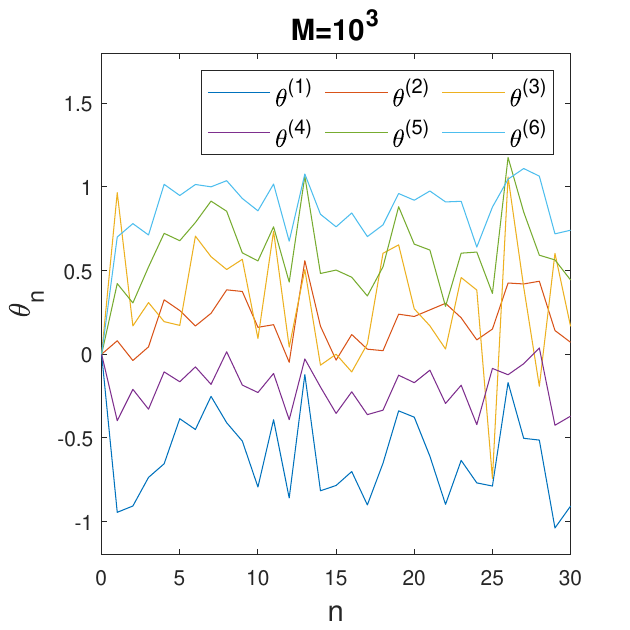}} &
		{\includegraphics[scale=0.75]{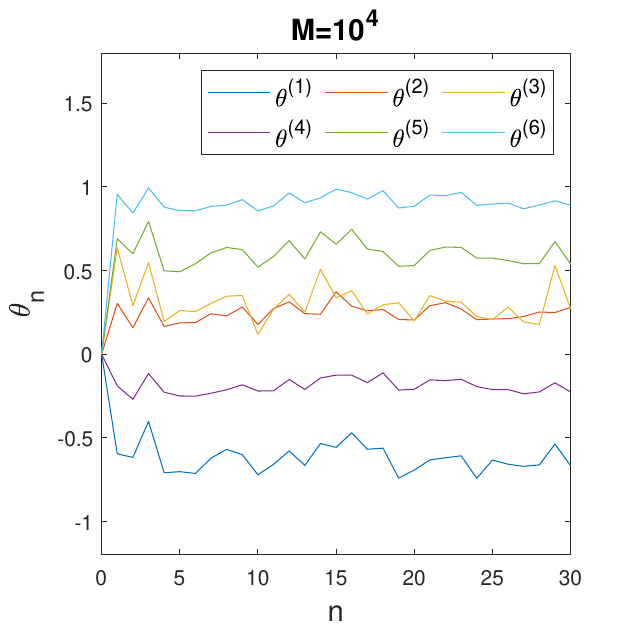}}\\
		{\includegraphics[scale=0.75]{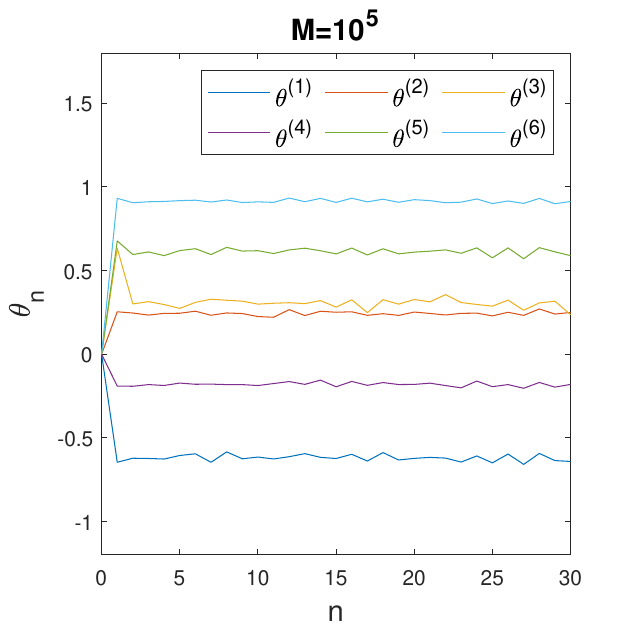}} &
		{\includegraphics[scale=0.75]{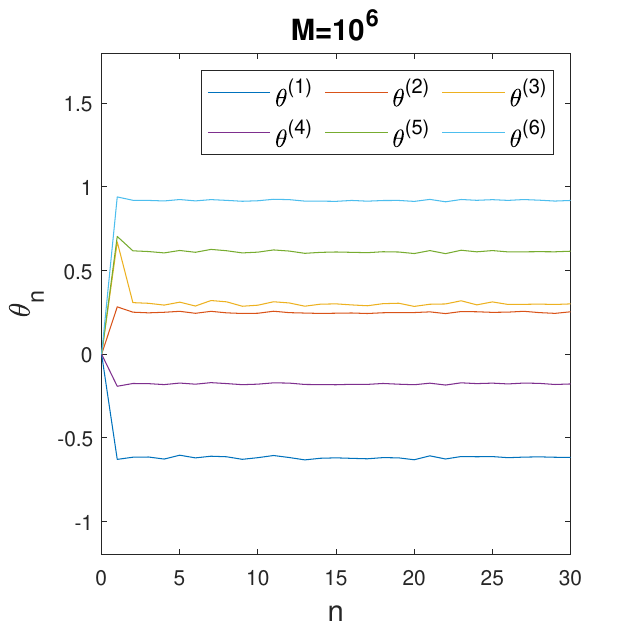}}
	\end{tabular}
	\caption[Values of weights $\btheta$ for Example 3.]{Values of $\theta^{(f)}_n$ in iteration $n$ of Algorithm~\ref{Al1} in  in Example~\ref{ex1}. Simulations were run for $M= 10^3, 10^4, 10^5$ and $10^6$ states, respectively.}
	\label{ex4_fig_theta}
\end{figure}

\begin{figure}[H]
	\begin{tabular}{cccc}
		{\includegraphics[scale=0.75]{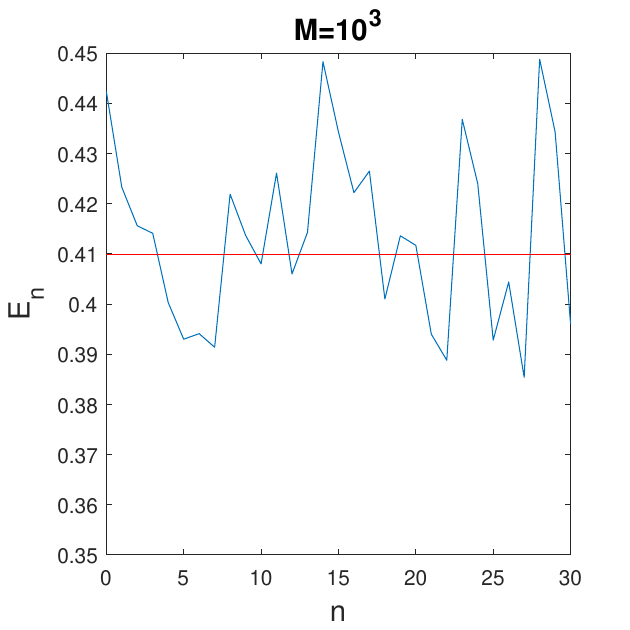}} &
		{\includegraphics[scale=0.75]{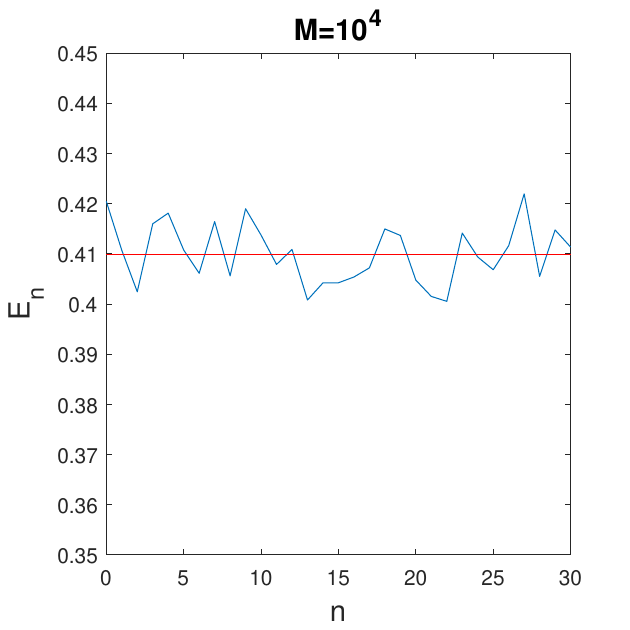}}\\
		{\includegraphics[scale=0.75]{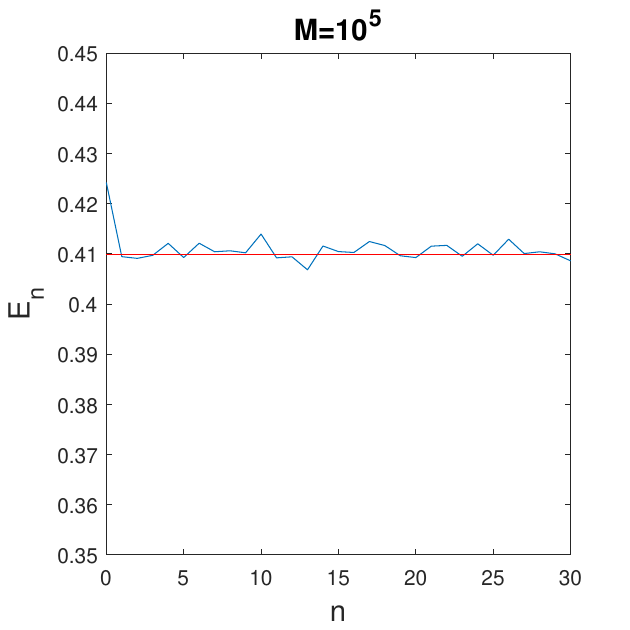}} &
		{\includegraphics[scale=0.75]{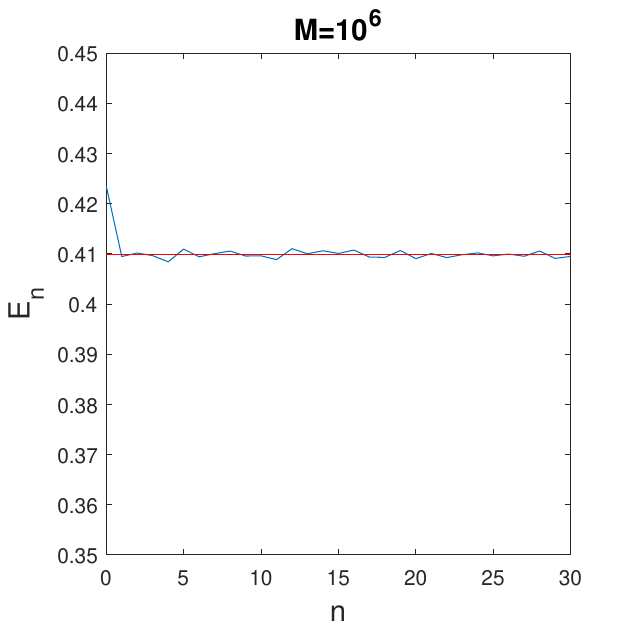}}
	\end{tabular}
	\caption[Values of $E$ for Example 3 compared to the optimal value $E^*=0.4098$.]{Values of $E_n$ (blue) compared to the optimal value $E^*=0.4098$ (red) in iteration $n$ of Algorithm~\ref{Al1} in Example~\ref{ex1}. Simulations were run for $M= 10^3, 10^4, 10^5$ and $10^6$ states, respectively.}
	\label{ex4_fig_E}
\end{figure}

We find that the estimated value of $E$ converges to the optimal value $E^*=0.4098$ for large values of $M$, after a small number of iterations. As expected, the values of $\theta$ and $E$ converge more tightly as $M$, the number of states simulated, increases, as shown in Figure~\ref{ex4_fig_E}. Second, we observe that in this example, each value of $\btheta$ converges to a distinct value, as shown in Figure~\ref{ex4_fig_theta}. This is because the parameters $\boldsymbol{p}$ and $\boldsymbol{\lambda}$ are asymmetric, so we do not expect patients of different types $i$ to be weighted equally.


\newpage
\subsection{Realistically-sized example}

We consider a realistically-sized example of our Model in Section~\ref{MDP_1}. The structure of the model here is influenced by the practical considerations within real-world hospitals, in which patients are allocated to their primary (preferred) wards based on their types (medical needs), and when these wards are full, suitable policies are applied to decide which nonprimary ward a patient should be allocated to instead. Furthermore, we assume that patients arriving to a hospital when the system is full (all wards are full), are redirected to another facility. Therefore, the number of daily arrivals is unrestricted, however the number of arrivals admitted to the hospital is bounded by the system capacity. We fitted the parameters of the arrival and length of stay distributions to data, and focused on the application of our model to the analysis of how decision making affects the number of patients in nonprimary wards.

\begin{exa}\label{ex2} Similary to Dai and Shi \cite{Dai2019}, assume five types of patients, $i=1,\ldots,5$ ($I=5$): Orthopedic (Orth), Cardio (Card), Surgery (Surg), General Medicine (GenMed), and Other Medicine (OthMed); and corresponding wards, $k=i=1,\ldots,5$ ($K=5$), most suitable to these patient types, respectively.

To estimate the key parameters of our model, we applied Matlab and performed the analysis of five-year patient flow data obtained from an Australian tertiary referral hospital, see Table~\ref{tab:parametersModel}.

\begin{table}[h]
	\centering
	
	\setlength{\tabcolsep}{1pt}	
	\begin{tabular} {l@{\hspace{1cm}}r@{\hspace{1cm}}r@{\hspace{1cm}}r@{\hspace{1cm}}r@{\hspace{1cm}}r}
		
		\hline
		$i$ &	$1$ (Ortho)&	$2$ (Card)&	$3$ (Surg)&	$4$ (GenMed)& $5$ (OthMed)\\
		\hline
		$\lambda_i$ &	$2.0252$&	$3.3565$&
	$10.0159$&	$11.7442$ &
			$38.7853$\\
		$\mathbb{E}(LoS_i)$ &	$5.1473$&	$4.1414$&$3.9373$&	
		$4.5209$ &		$2.8505$\\
		$p_{k,i}$, $k=i$ &$0.1943$&	$0.2415$&	
		$0.2540$&		
		$0.2212$ &		$0.3508$\\  	
$m_k$, $k=i$ &	$12$&	$15$&
$38$&	$50$ &	$99$\\		     
		\hline
	\end{tabular}
	\caption{
		Model parameters in Example~\ref{ex2} ($\lambda_i$ and $\mathbb{E}(LoS_i)$ were fitted to data, and the remaining parameters were estimated): $\lambda_i$ (daily arrival rate of patients type $i$), $p_{k,i}$ (probability of patient of type~$i$ departing from ward $k$ within a day), and $m_k$ (capacity of ward $k$). The total capacity of the system is $m=\sum_k m_k = 214$. We apply $p_{i,i}=1/\mathbb{E}(LoS_i)$ (probability of patient of type $i$ departing within a day, assuming they are in the most suitable ward). When $k\not=i$, we let $p_{k,i}=\beta\times p_{i,i}$, for 
		some $\beta<1$ (so that the LoS of patients that are in less suitable wards, is increased). Here, we set 
		$\beta=1/1.25$ (and so the LoS of patients in less suitable wards increases by $25\%$ on the average). 
		Further, to estimate $m_k$, for each $k=i$, we applied M/M/N/N queueing model with arrival rate $\lambda_i$ and service rate $\mu_i=1/\mathbb{E}(LoS_i)$ to find $N=m_k$ such that the proportion of time the queue is full is $\pi_N<0.15$.%
		}
	\label{tab:parametersModel}
\end{table}

\noindent In order to focus on the impact of decision making on the number of patients observed in nonprimary wards, we applied the following cost parameters:
\begin{itemize}

	\item Transfer cost from ward $k$ to ward $\ell$ per type-$i$ patient: $ c^{(t)}_{k,\ell,i}=c^{(t)}=1.1$ for all $k,\ell\not=k, i$ and assume large value at $k=\ell$ (to avoid such decision).

	\item Penalty cost for being in a nonprimary ward $k$ per type-$i$ patient:  $c^{(p)}_{k,i}=c^{(p)}=0.2$ for all $k, i$.
\end{itemize}
(Here, we assumed the assignment costs  $c^{(\sigma)}=1$ for each admitted or redirected patient, but have not included these costs in the optimisation, since all patients need to be admitted to some facility, and our focus is on the number of patients in nonprimary wards.)

We assume that the capacity of the waiting area is unlimited with $m_{k+1}=\infty$ and so the total number of arrivals $\sum_{i=1}^I q_i\leq m_{k+1}$ is unrestricted. However, the arrivals are only accepted if there is available capacity in the system according to $\sum_{i=1}^I n'_{k,i}\leq m_k$, where $\left( [n'_{k,i}]_{\mathcal{K} \times \mathcal{I}} \right)$ is post-decision state. Patients that may not be admitted, are redirected to another facility. The state space of the system is then
\begin{eqnarray*}
	\mathcal{S}&=&
	\{
	\left( [n_{k,i}]_{\mathcal{K} \times \mathcal{I}}, [q_i]_{1 \times \mathcal{I}} \right):
	n_{k,i}\geq 0,q_i\geq 0,  \sum_{i=1}^{I} n_{k,i}\leq m_k
	\},
\end{eqnarray*} 
while the set of post-decision states is
\begin{eqnarray*}
	\mathcal{S}' &=& 
	\left\{
	\left( [n'_{k,i}]_{\mathcal{K} \times \mathcal{I}} \right)
	: 
	\sum_{i=1}^I n'_{k,i}\leq m_k
	\right\}\ ,
\end{eqnarray*}
with $|\mathcal{S}'|=\prod_{k=1}^K\sum_{g=0}^{m_k}\binom{g+I-1}{I-1} =\prod_{k=1}^K\binom{m_k+I}{I}$ by Heydar et al.~\cite{Heydar2021}. Here, $|\mathcal{S}'|=2.9544e+28$.

Next, we assume the following order of assignment, represented as a matrix ${\bf O}=[O_{i,k}]$, where $O_{i,k}=1$ means that ward $k$ is the first choice for type $i$, $O_{i,k}=2$ is the second choice for type $i$, and so on, with
\begin{eqnarray}\label{eq:order}
{\bf O}&=&
\left[
\begin{array}{l||ccccc}
&k=\mbox{Ortho}&k=\mbox{Card}&k=\mbox{Surg}&k=\mbox{GenMed}&k=\mbox{OthMed}\\\hline\hline
i=\mbox{Ortho}&1&     5&     2&     3&     4\\\hline
i=\mbox{Card}&3&     1&     2&     4&     5\\\hline
i=\mbox{Surg}&2&     3&     1&     5&     4\\\hline
i=\mbox{GenMed}&3&     5&     4&     1&     2\\\hline
i=\mbox{OthMed}&3&     5&     4&     2&     1
\end{array}
\right].
\end{eqnarray}
As example, for the assignment of 'Ortho' patients, the considered order of wards is `Ortho', then `Surg', then `GenMed', then `OthMed', and `Card' is the last choice (used when there are no other choices). This order is similar to the order of assignment in Dai and Shi~\cite[Figure 1]{Dai2019}.
\end{exa}

Suppose that patient type represents the order of assignment, and so we assign patients in the order of their type, assigning type $1$ patients first, then type $2$ patients, and so on. We consider two decisions, 
\begin{itemize}
	\item $a=1$, assign arrived patients to their best available wards, in the order of their types, \underline{without} transferring patients between the wards, and
	\item $a=2, 3$, assign arrived patients to their best available wards, in the order of their types, and \underline{allow} transferring patients if required, with the constraint of no more than $y=4, 10$ transfers in total, respectively.
\end{itemize}

As example, if 
\begin{itemize}
	\item $q_1=7$, $q_2=3$, $q_i=0$ for $i=3,4,5$ (there are $7$ type-$1$ patients and $3$ type-$2$ patients in the queue); 
	\item and $n_{1,1}=11$, $n_{1,2}=4$, $n_{1,i}=0$ for $i=3,4,5$ (there are $22$ patients in ward $1$, so it is full), 
\end{itemize} 
then applying $a=2$ means that
\begin{itemize}
	\item $4$ type-$1$ patients will be assigned to ward $1$, while $3$ type-$1$ patients will be assigned to their best available ward other than ward $1$, 
	\item while $4$ type-$2$ patients will need to be transferred from ward $1$ to their best available ward (to make space for type-$1$ patients), 
	\item and finally, $3$ type-$2$ patients will be assigned from the queue to their best available ward.
\end{itemize} 	
More precisely, the assignments of patients corresponding to decisions $a=1$ and $a=2, 3$, are described in Algorithms~\ref{a4}-\ref{a5} in Appendix~\ref{sec:algorithms}, respectively.

\begin{figure}[H]
	\centering	{\includegraphics[scale=0.4]{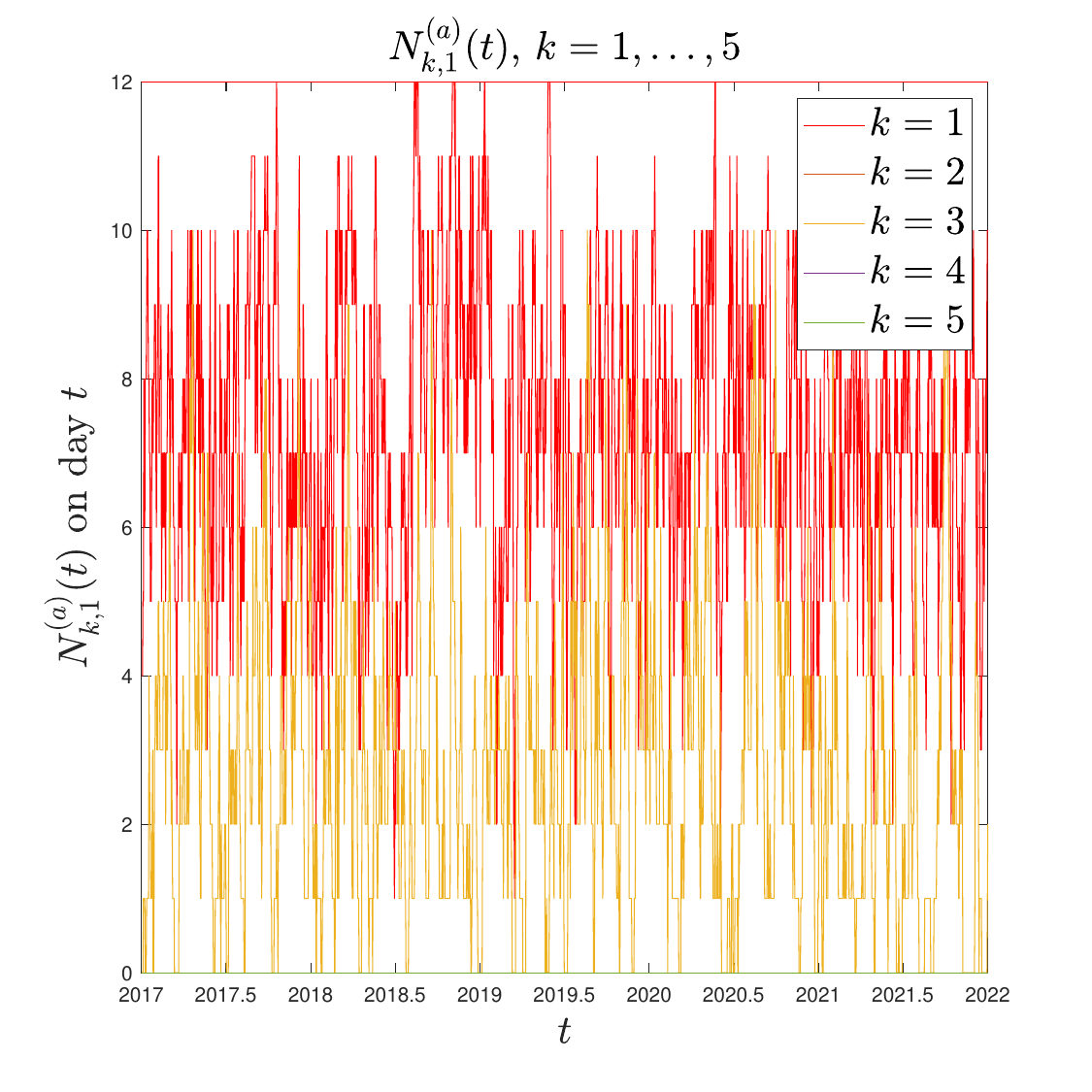}}
	\ 
	{\includegraphics[scale=0.4]{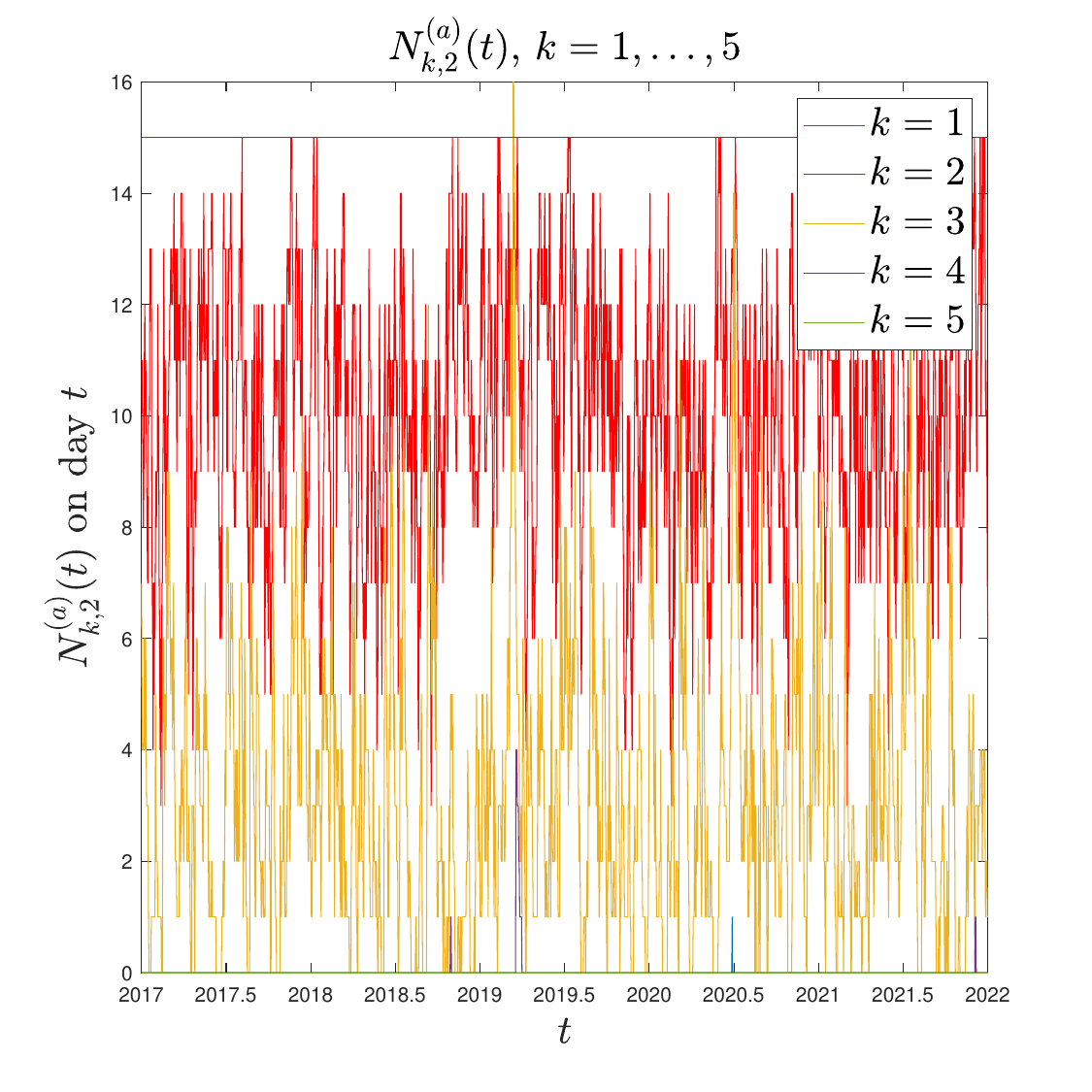}}\\
	{\includegraphics[scale=0.4]{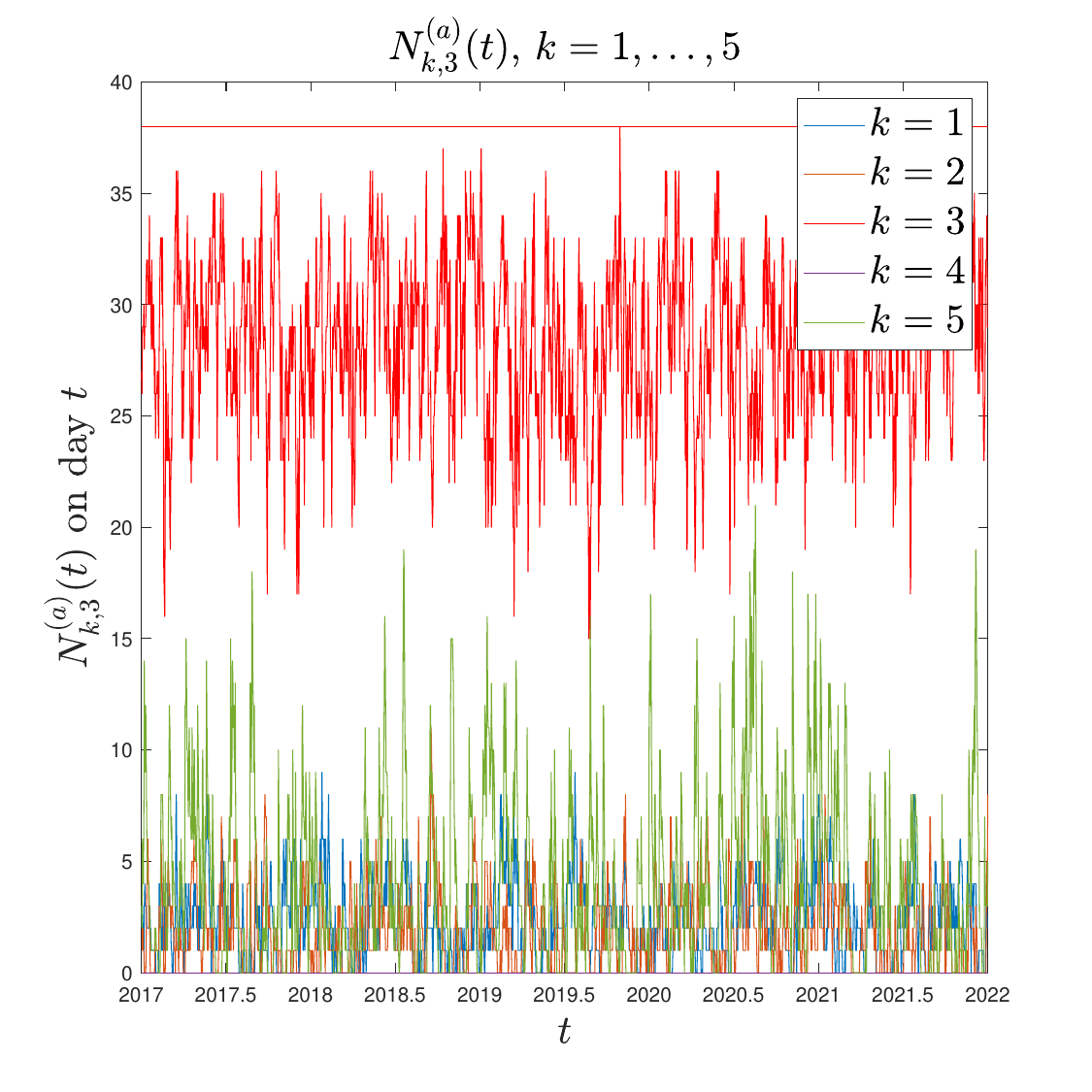}}
	\ 
	{\includegraphics[scale=0.4]{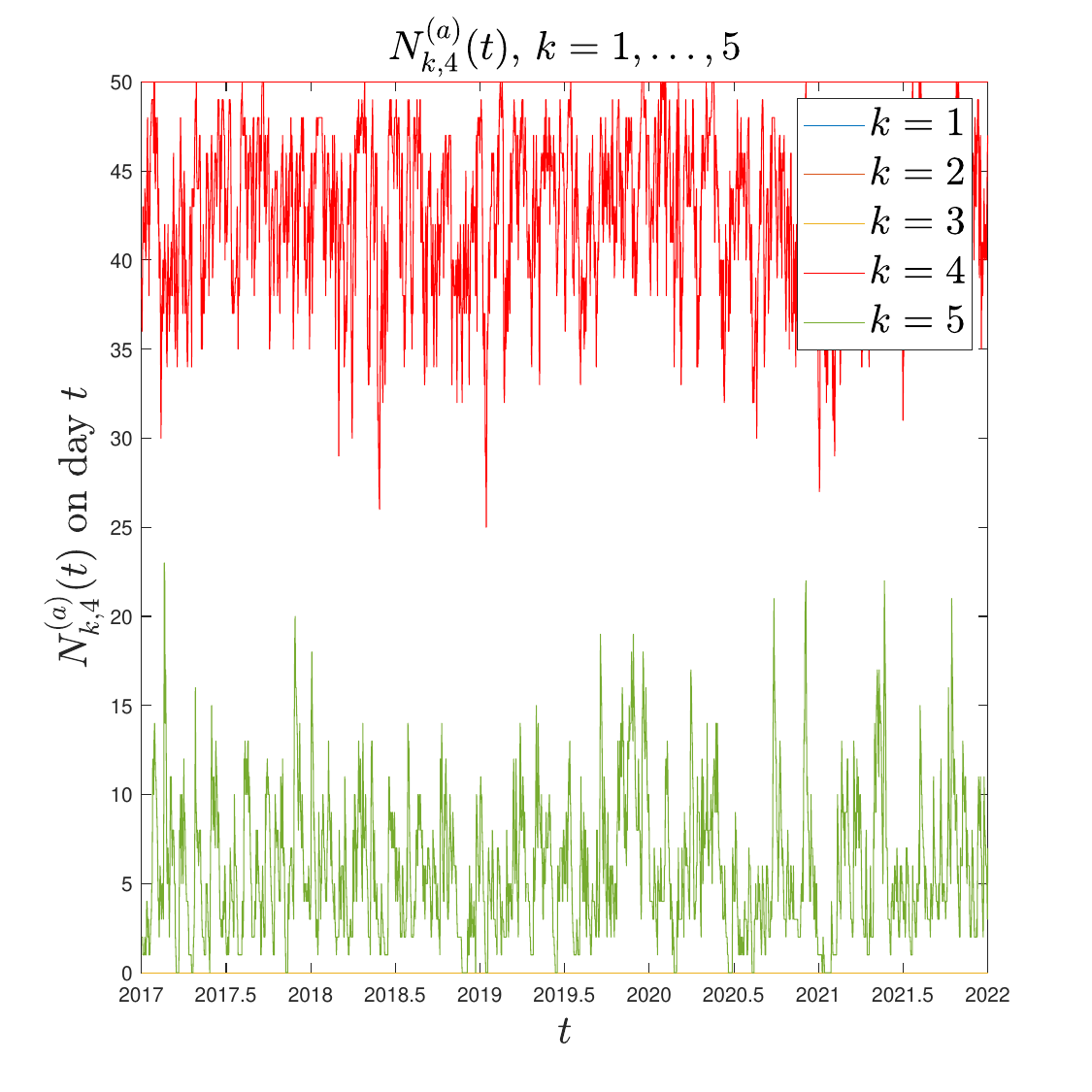}}\\
	{\includegraphics[scale=0.4]{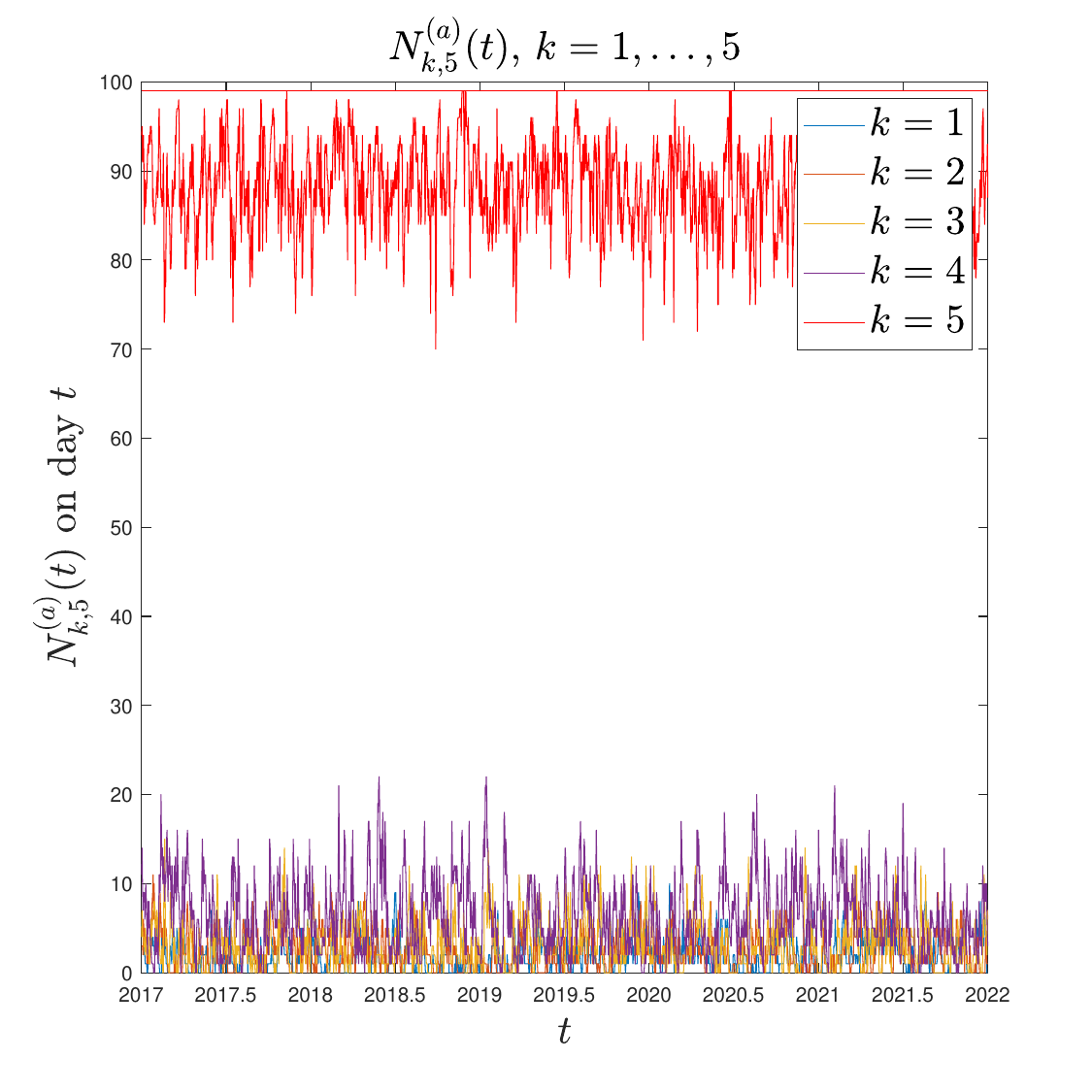}}
	\ 
	{\includegraphics[scale=0.4]{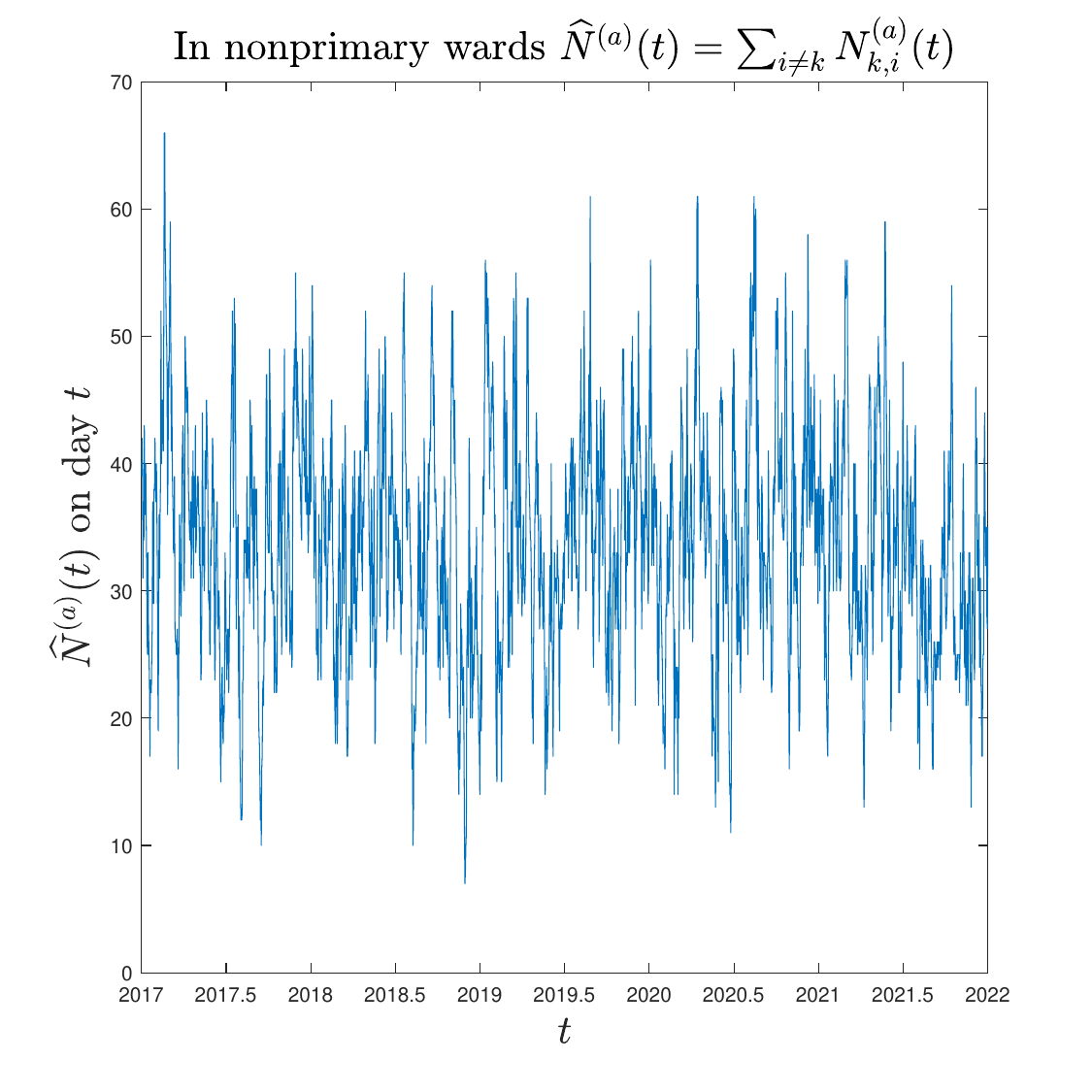}}\\
	\caption{Simulation of a system with limited capacity under policy that applies $a=1$ (no transfers): $N^{(a)}_{k,i}(t)$ is the number of type $i$ patients in ward $k$, $\widehat N^{(a)}(t)=\sum_{k\not=i}N^{(a)}_{k,i}(t)$ is the number of patients in nonprimary wards on day $t$, post decision. We apply~\eqref{eq:order} and the parameters from Table~\ref{tab:parametersModel}.}
	\label{ex2_bounded_a1}
\end{figure}

\begin{figure}[H]
	\centering	{\includegraphics[scale=0.4]{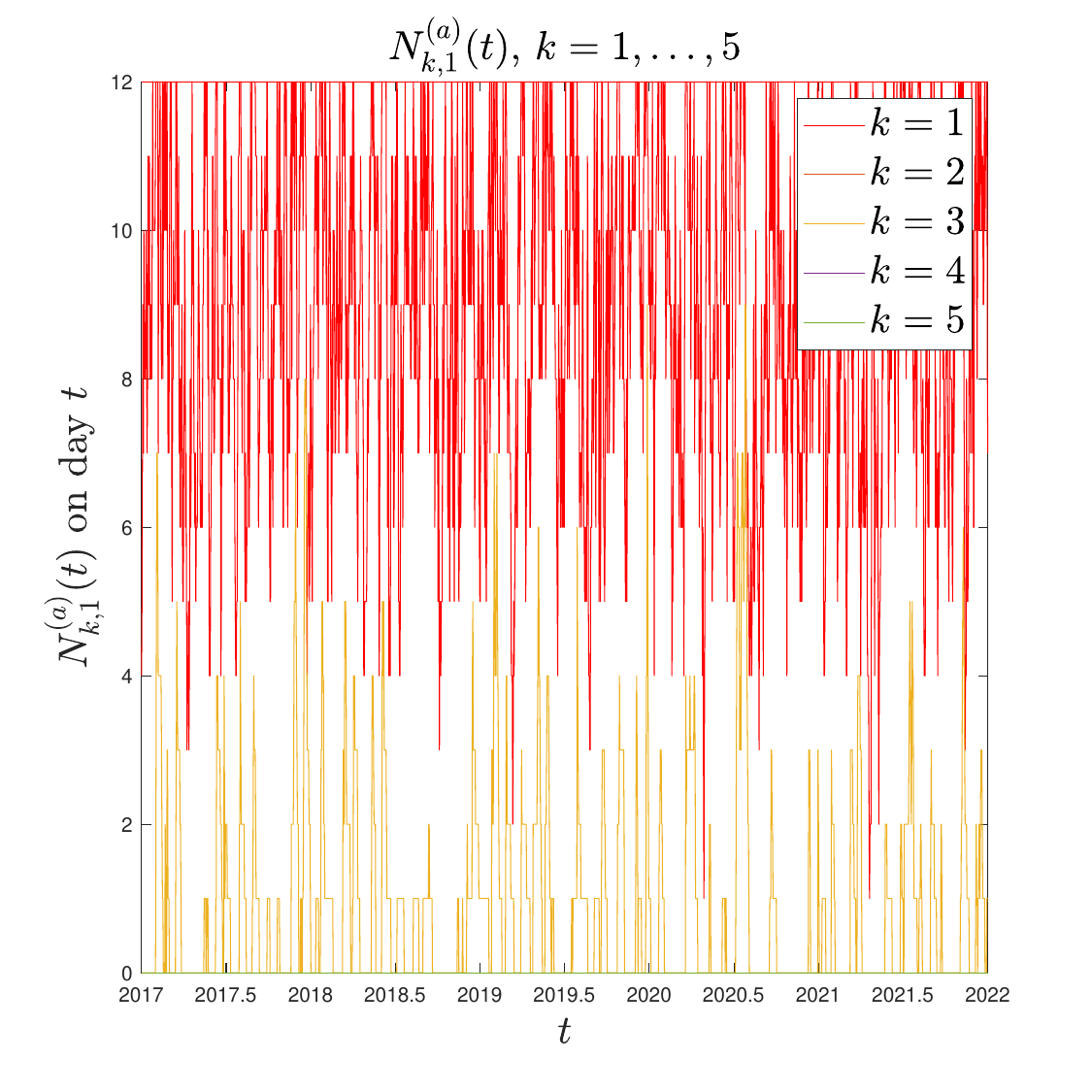}}
	\ 
	{\includegraphics[scale=0.4]{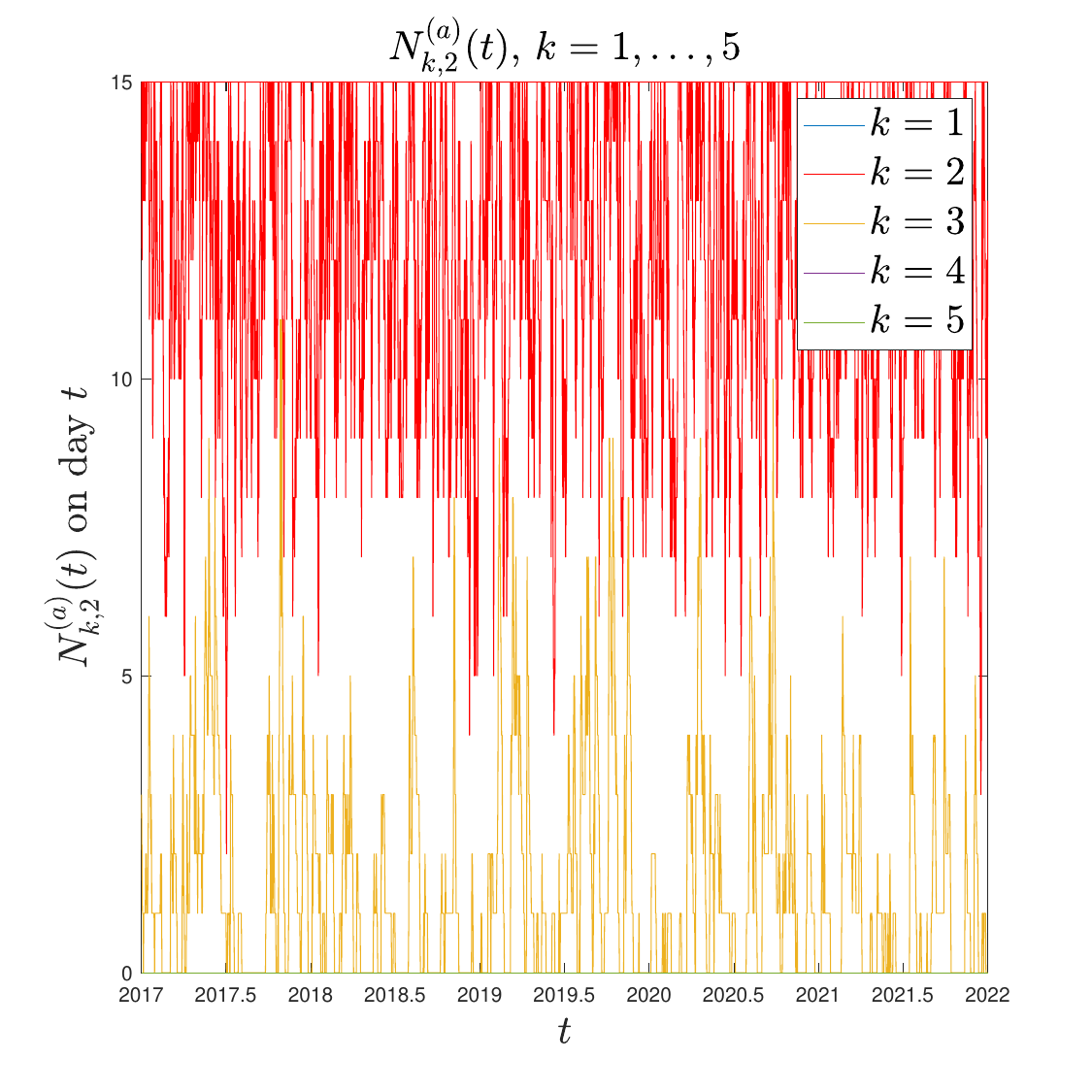}}\\
	{\includegraphics[scale=0.4]{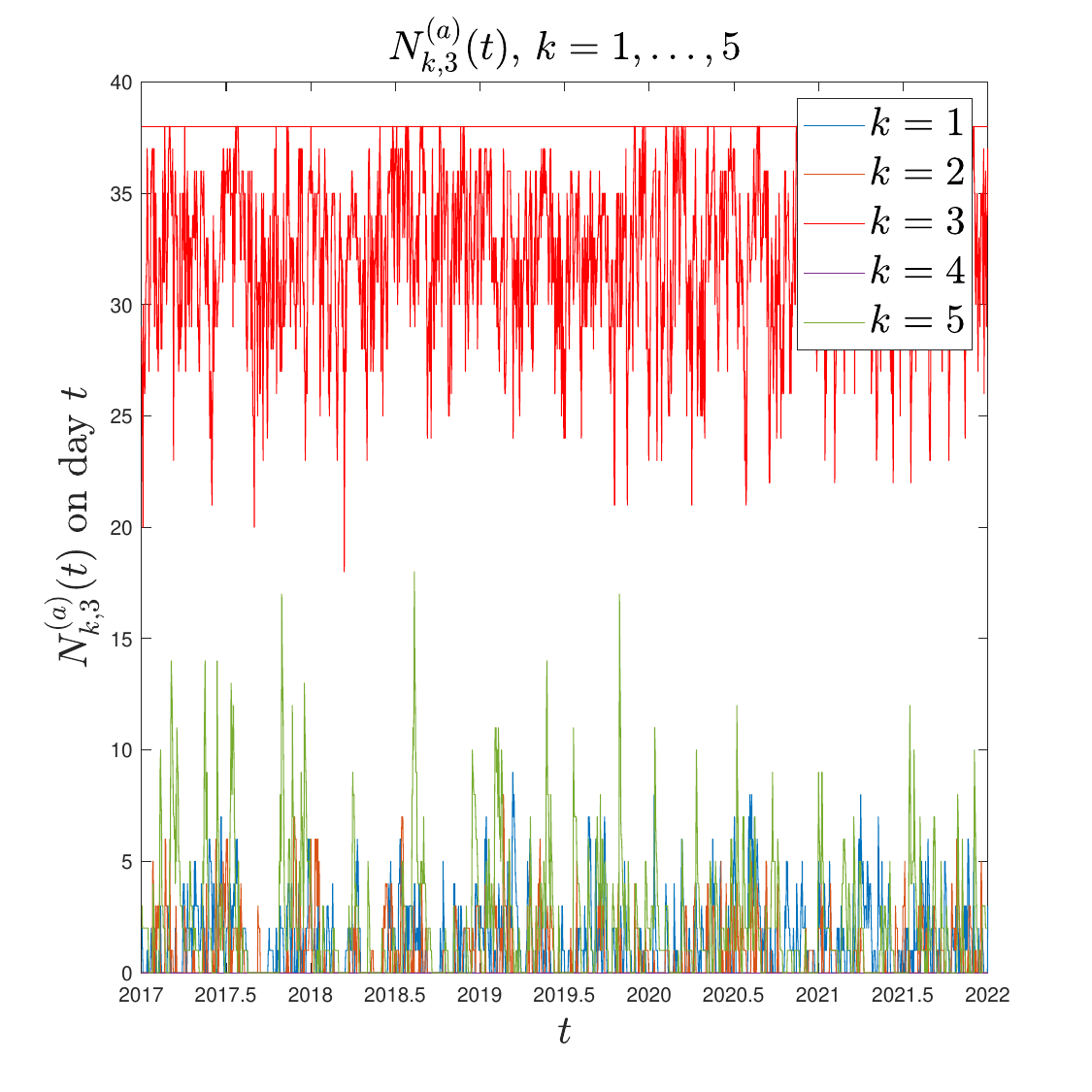}}
	\ 
	{\includegraphics[scale=0.4]{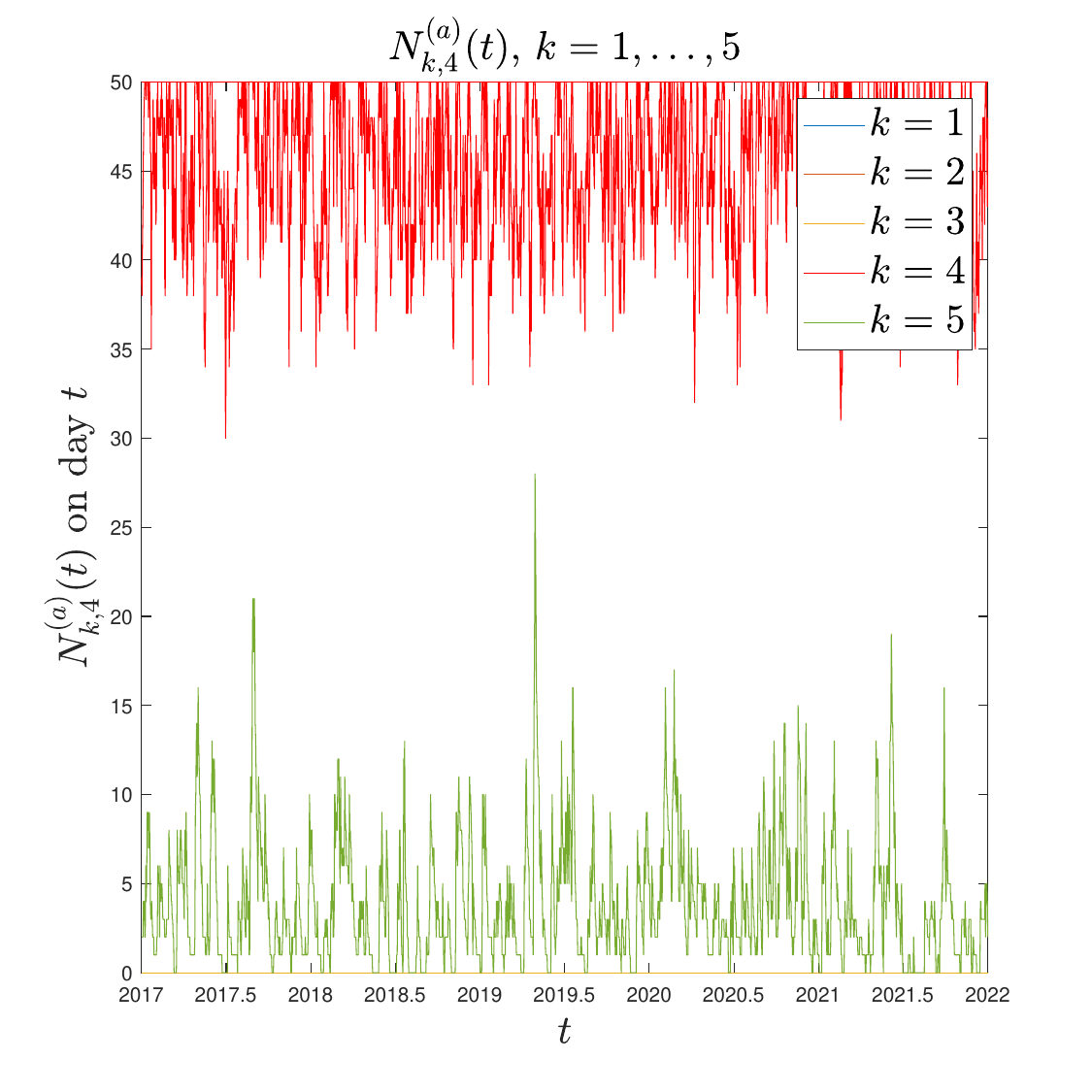}}\\
	{\includegraphics[scale=0.4]{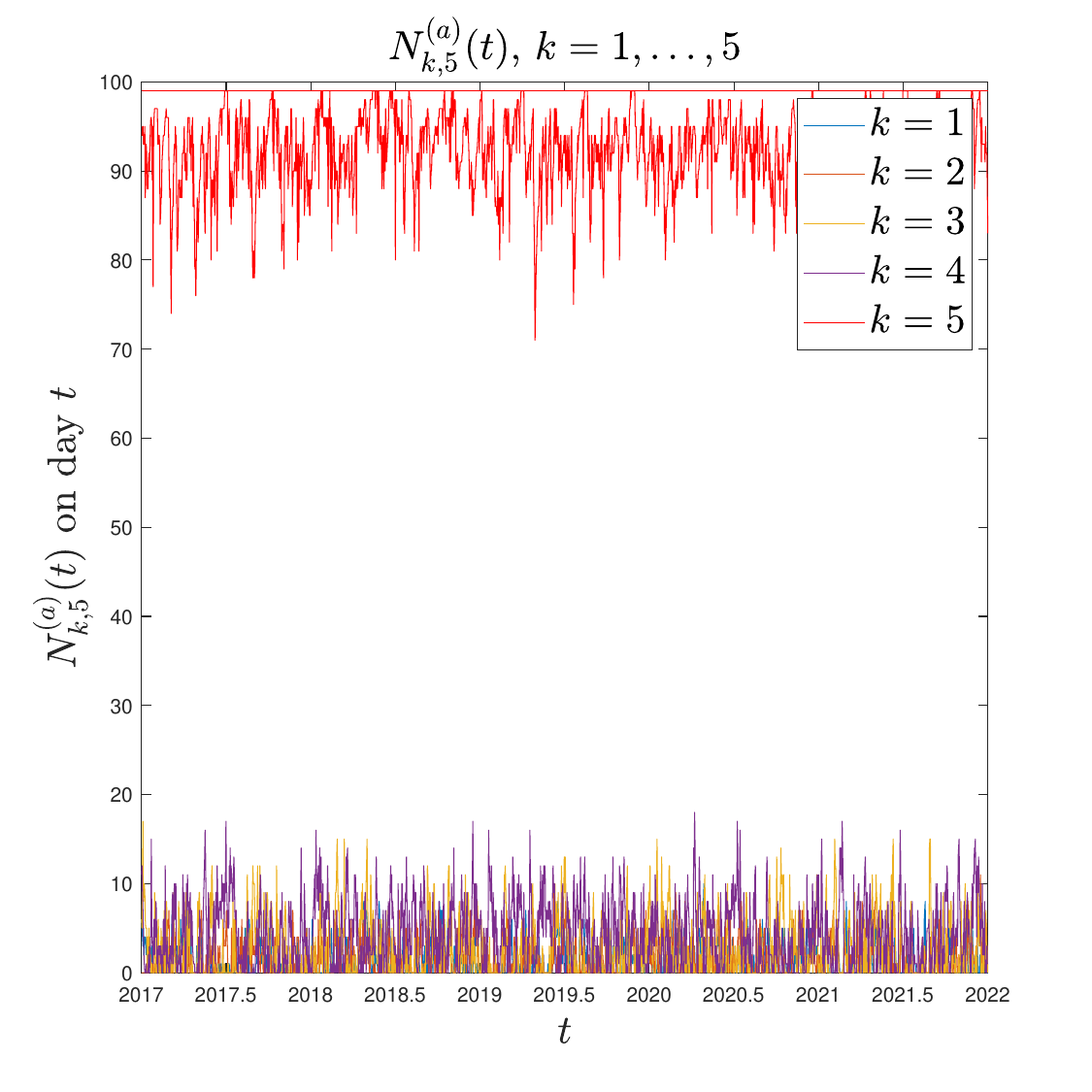}}
	\ 
	{\includegraphics[scale=0.4]{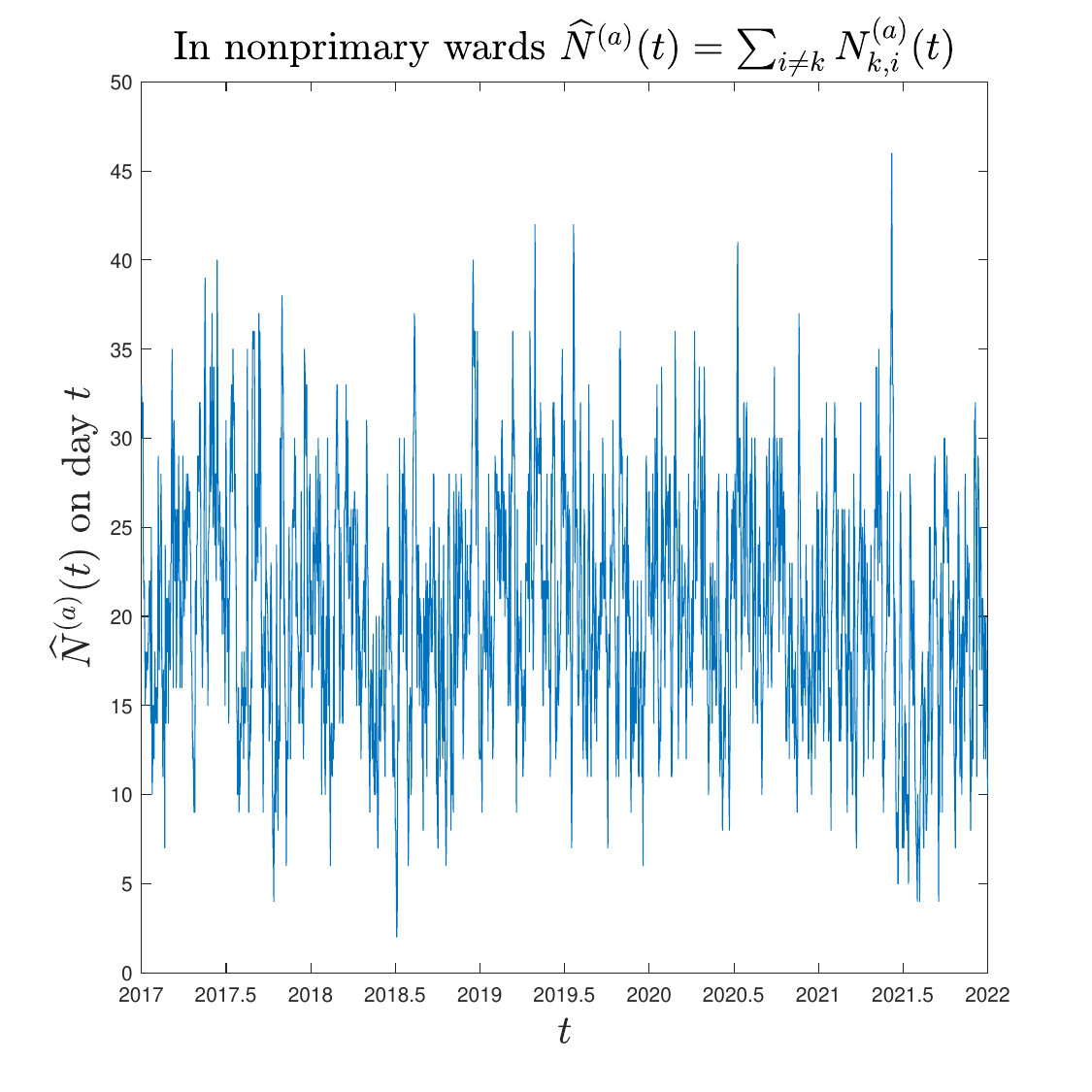}}\\
	\caption{Simulation of a system with limited capacity under policy that applies $a=2$ (with no more than $y=4$ transfers): 
		We note the reduction of the total number of patients in nonprimary wards in comparison to the output for $a=1$ in Figure~\ref{ex2_bounded_a1}.}
	\label{ex2_bounded_a2}
\end{figure}

\begin{figure}[H]
	\centering	{\includegraphics[scale=0.4]{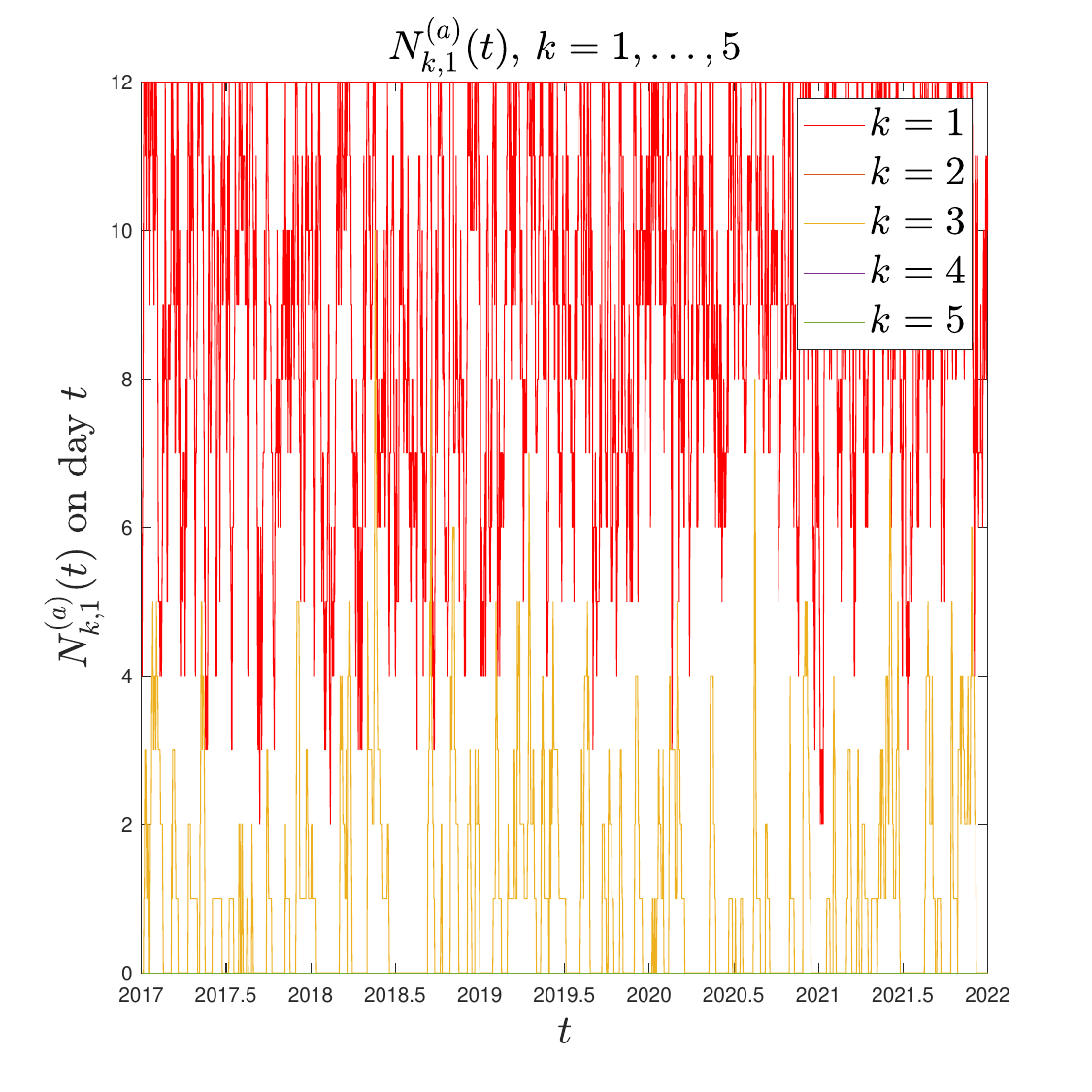}}
	\ 
	{\includegraphics[scale=0.4]{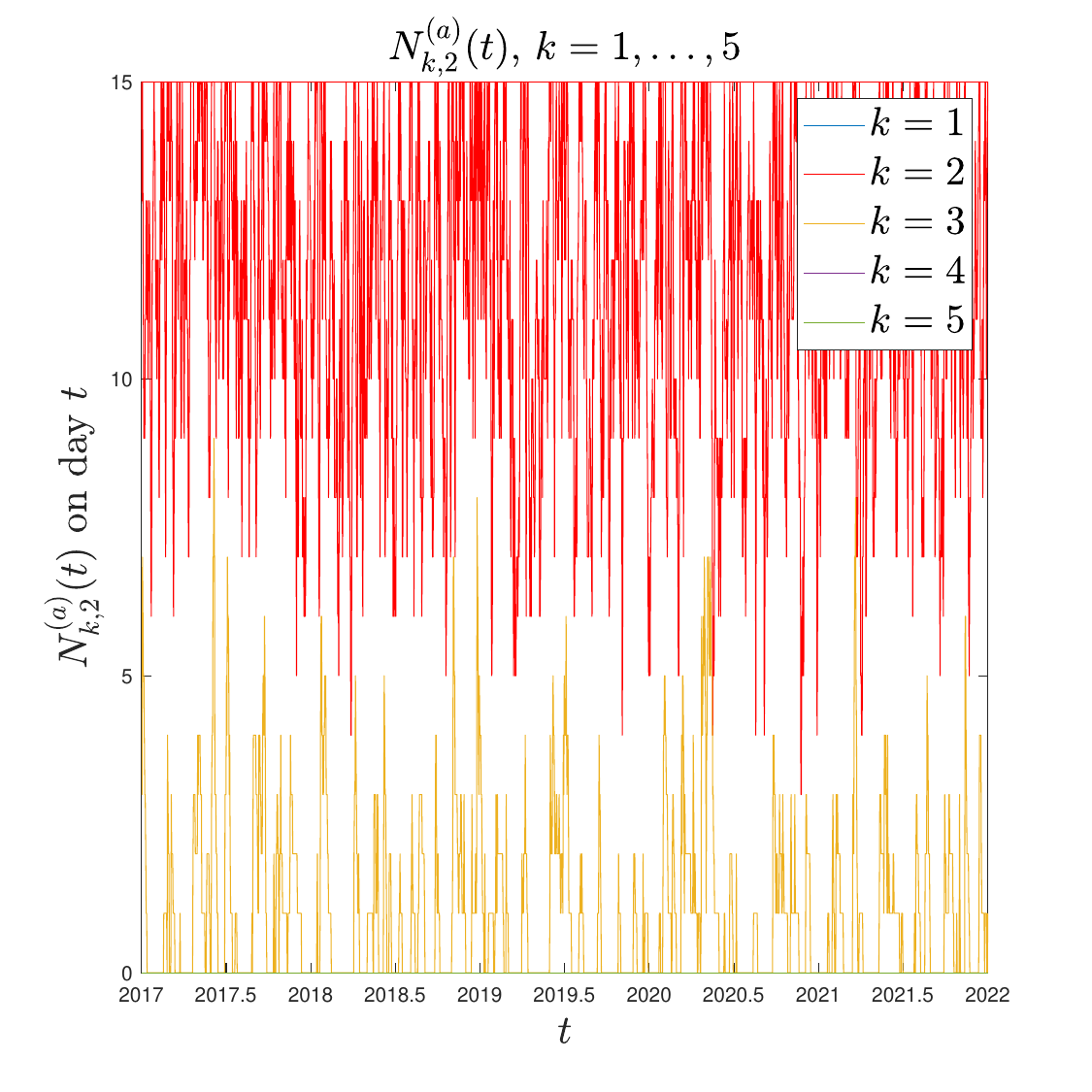}}\\
	{\includegraphics[scale=0.4]{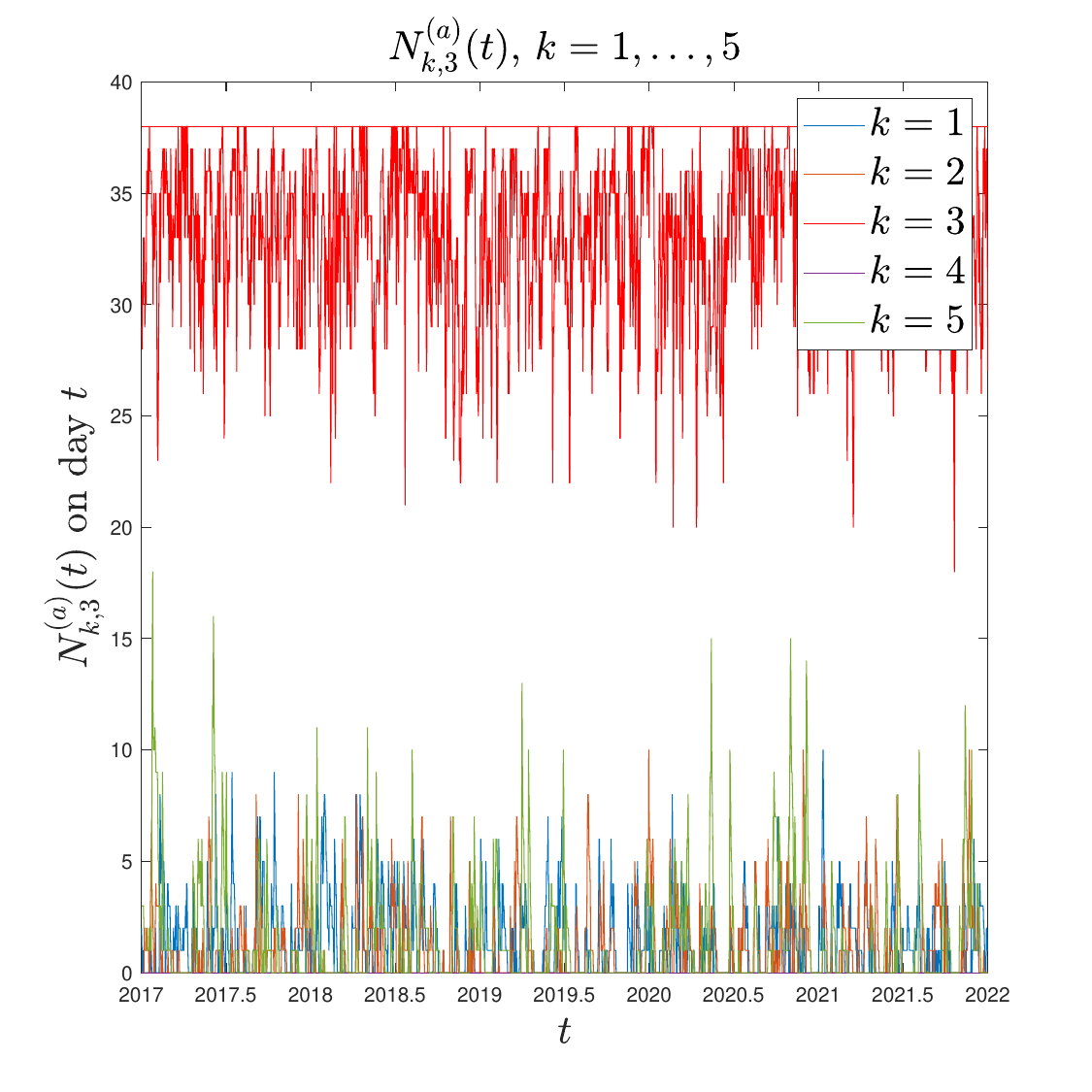}}
	\ 
	{\includegraphics[scale=0.4]{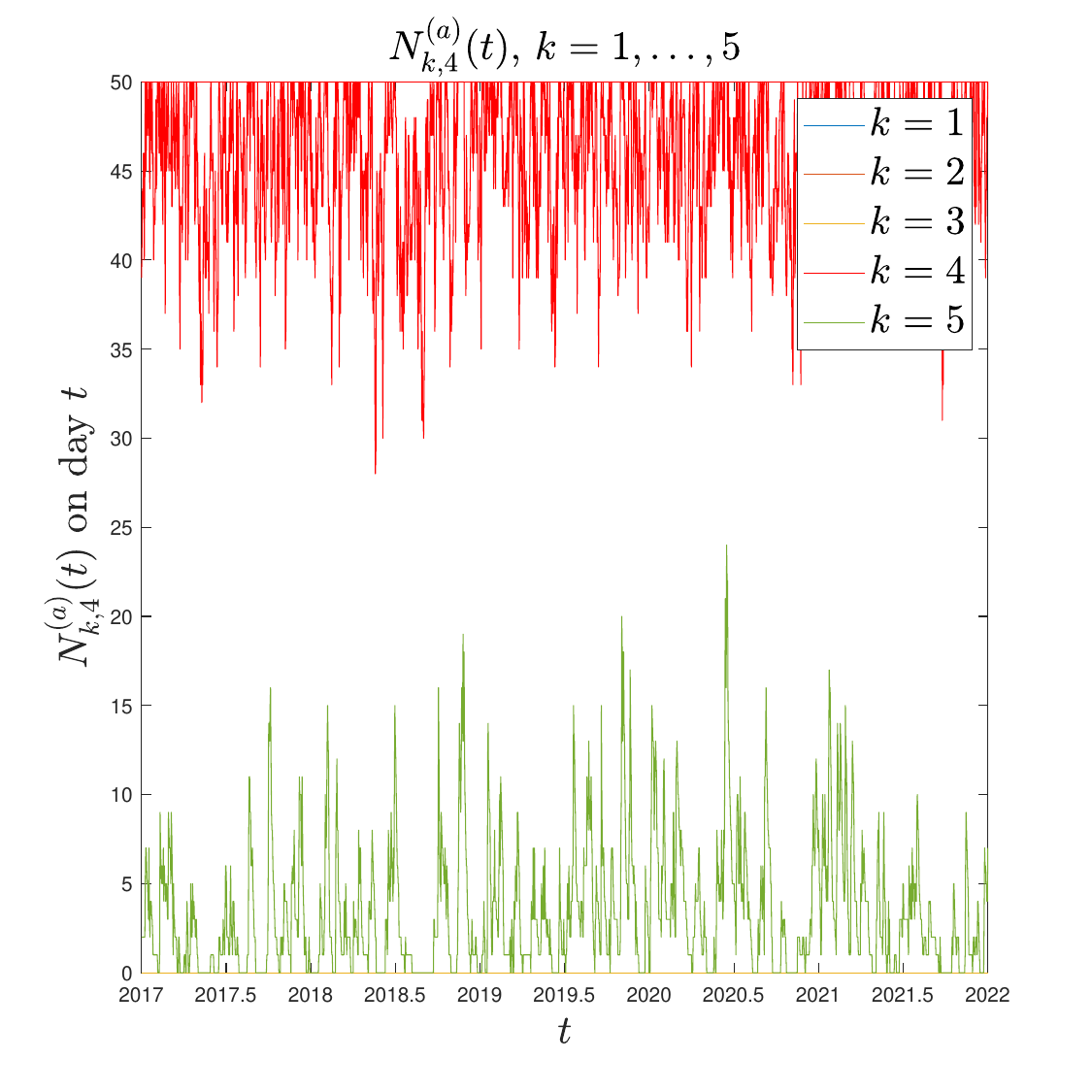}}\\
	{\includegraphics[scale=0.4]{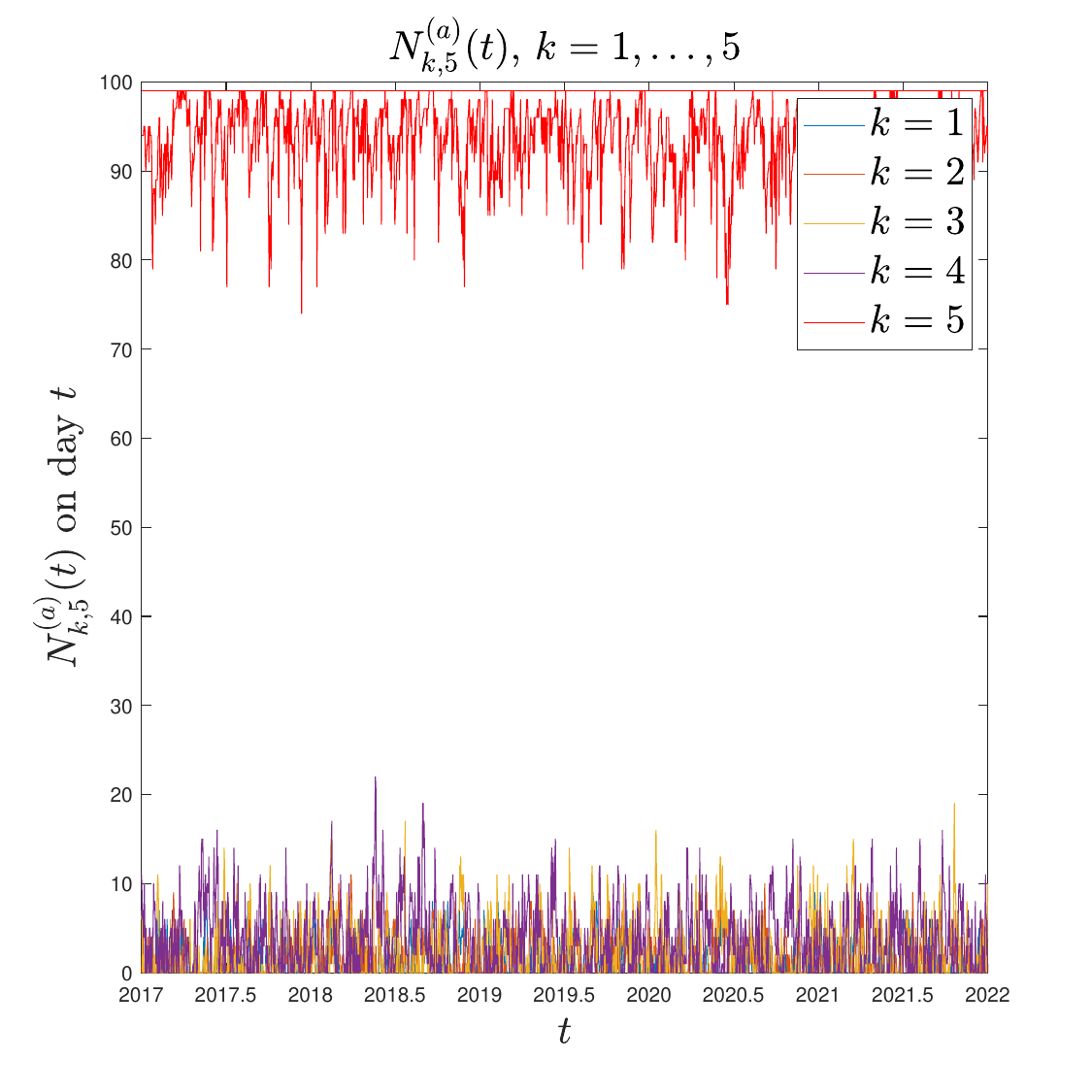}}
	\ 
	{\includegraphics[scale=0.4]{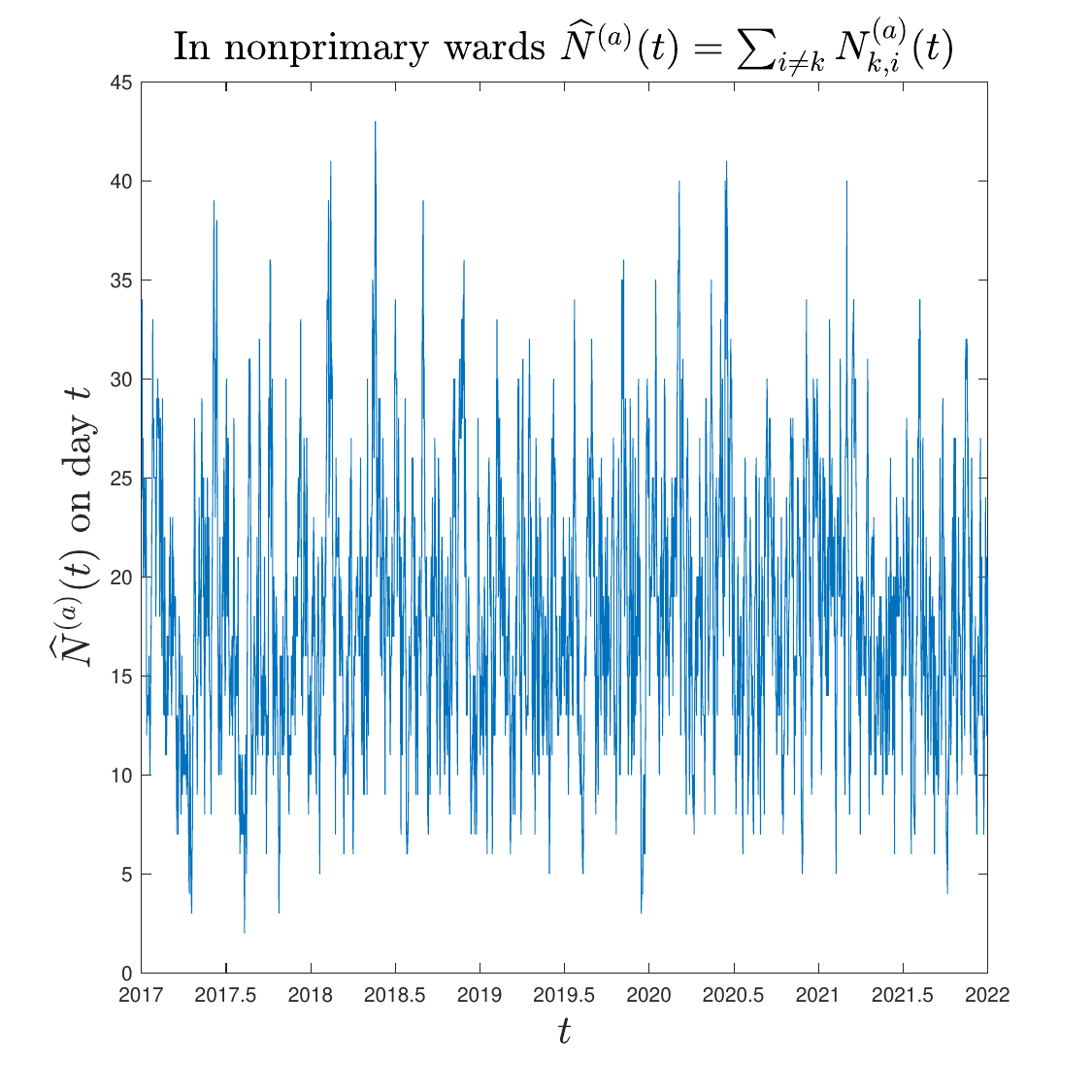}}\\
	\caption{Simulation of a system with limited capacity under policy that applies $a=3$ (with no more than $y=10$ transfers): 
		We note the further reduction of the total number of patients in nonprimary wards in comparison to the output for $a=2$ in Figure~\ref{ex2_bounded_a2}.}
	\label{ex2_bounded_a3}
\end{figure}

Next, we present the simulation output in Figures~\ref{ex2_bounded_a1}-\ref{ex2_bounded_a3}, focusing on the impact of these decisions on the number of patients observed that are in nonprimary wards (and ignoring the costs of these decisions for now). Simulation output suggests that under the model parameters assumed in Table~\ref{tab:parametersModel}, a policy in which we apply decision $a=3$ all the time, results in a better performance of the system in the sense that it reduces the number of patients in nonprimary wards when compared with applying decision $a=1$ or $a=2$ all the time. The performance of the system under decision $a=3$ is only slightly better than when compared with decision $a=2$ however.

In order to evalute the performance of the system under decisions $a=1,2,3$, we require to also consider the costs. Since the size of the state space is very large, we apply the Approximate Policy Iteration summarised in Algorithms~\ref{Al1}--\ref{Al1E} in Appendix~\ref{sec:algorithms}.

We consider the vector of features $\bphi(s)$ of state $s=\left( [n_{k,i}]_{\mathcal{K} \times \mathcal{I}}, [q_i]_{1 \times \mathcal{I}}\right) $ defined as
\begin{eqnarray}
\bphi(s)
&=&
\left[
\begin{array}{cccccccccc}
n_{1,1},&\sum_{i\not= 1}n_{1,i},&\ldots ,&
n_{5,5},&\sum_{i\not= 5}n_{5,i},
&q_1, &\ldots ,& q_5
\end{array}
\right],
\end{eqnarray}
that is, for each ward $k=1,\ldots,5$, $n_{k,k}=$ the total number of type-$i$ patients in the primary ward $i=k$; $\sum_{i\not= k}n_{k,i}$, the total number of other-type patients in that ward; and for each patient type $i=1,\ldots,5$, $q_i=$ the number of type-$i$ patients in the queue waiting to be assigned.

\begin{figure}[H]
	\begin{tabular}{cc}
{\includegraphics[scale=0.5]{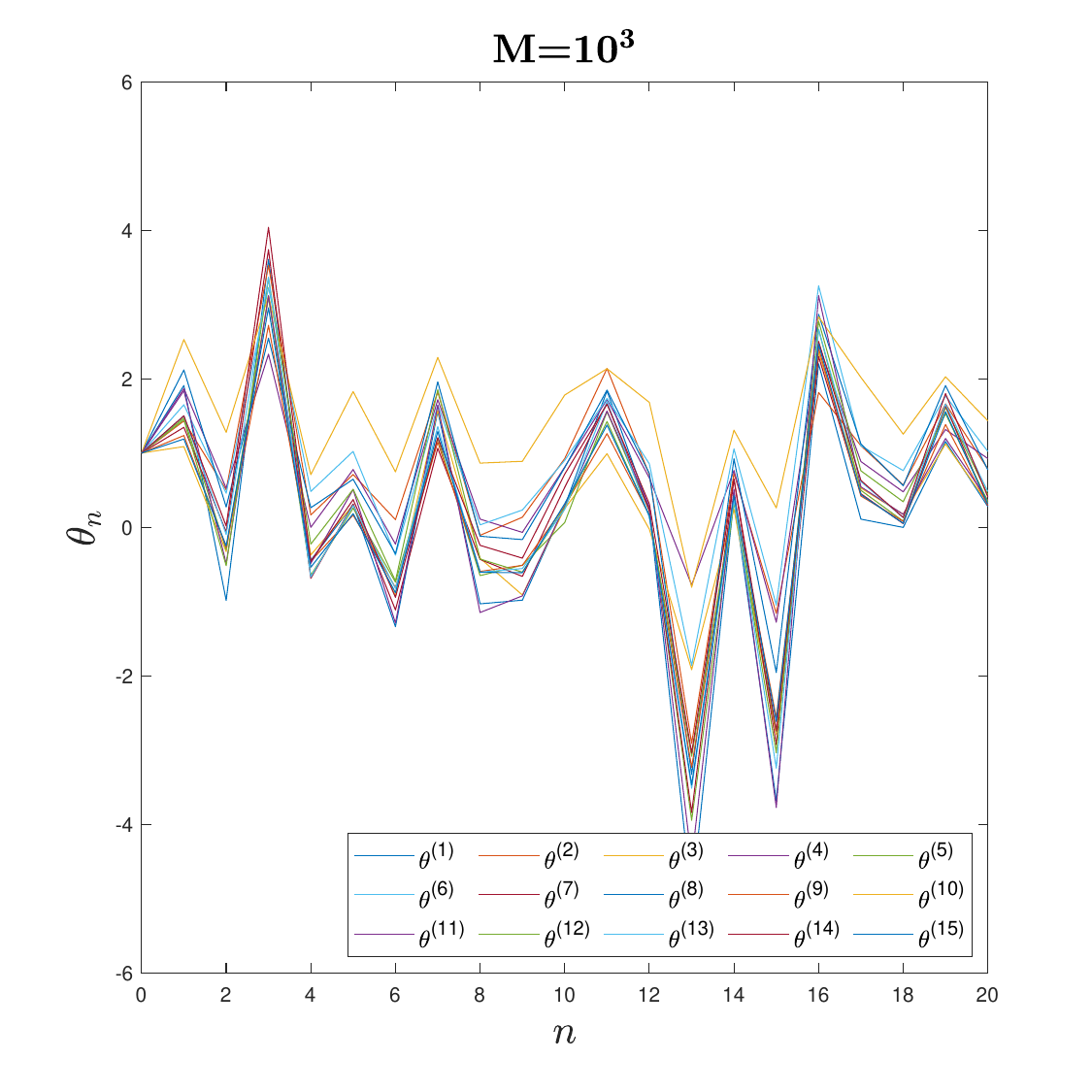}}
 &
{\includegraphics[scale=0.5]{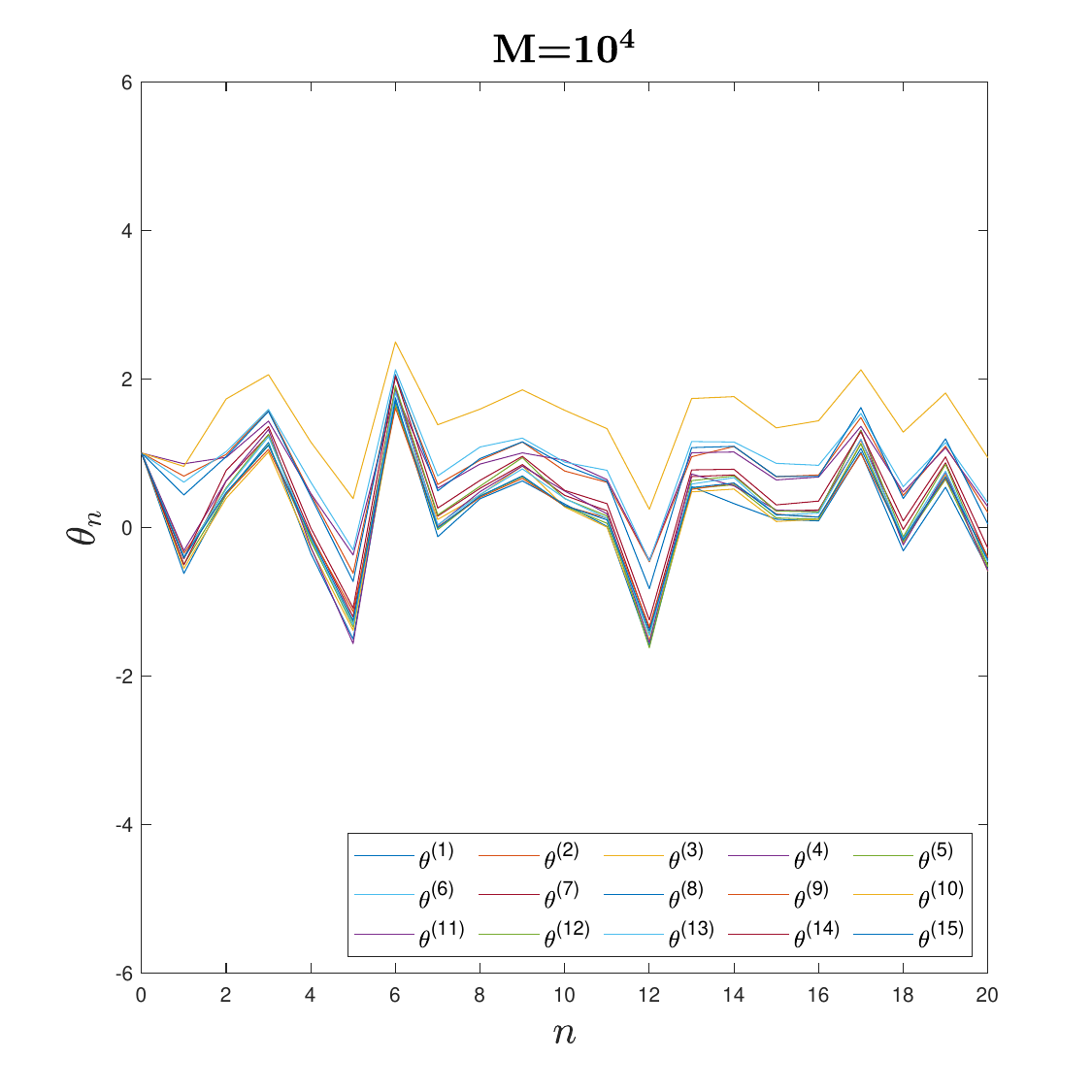}}\\
{\includegraphics[scale=0.5]{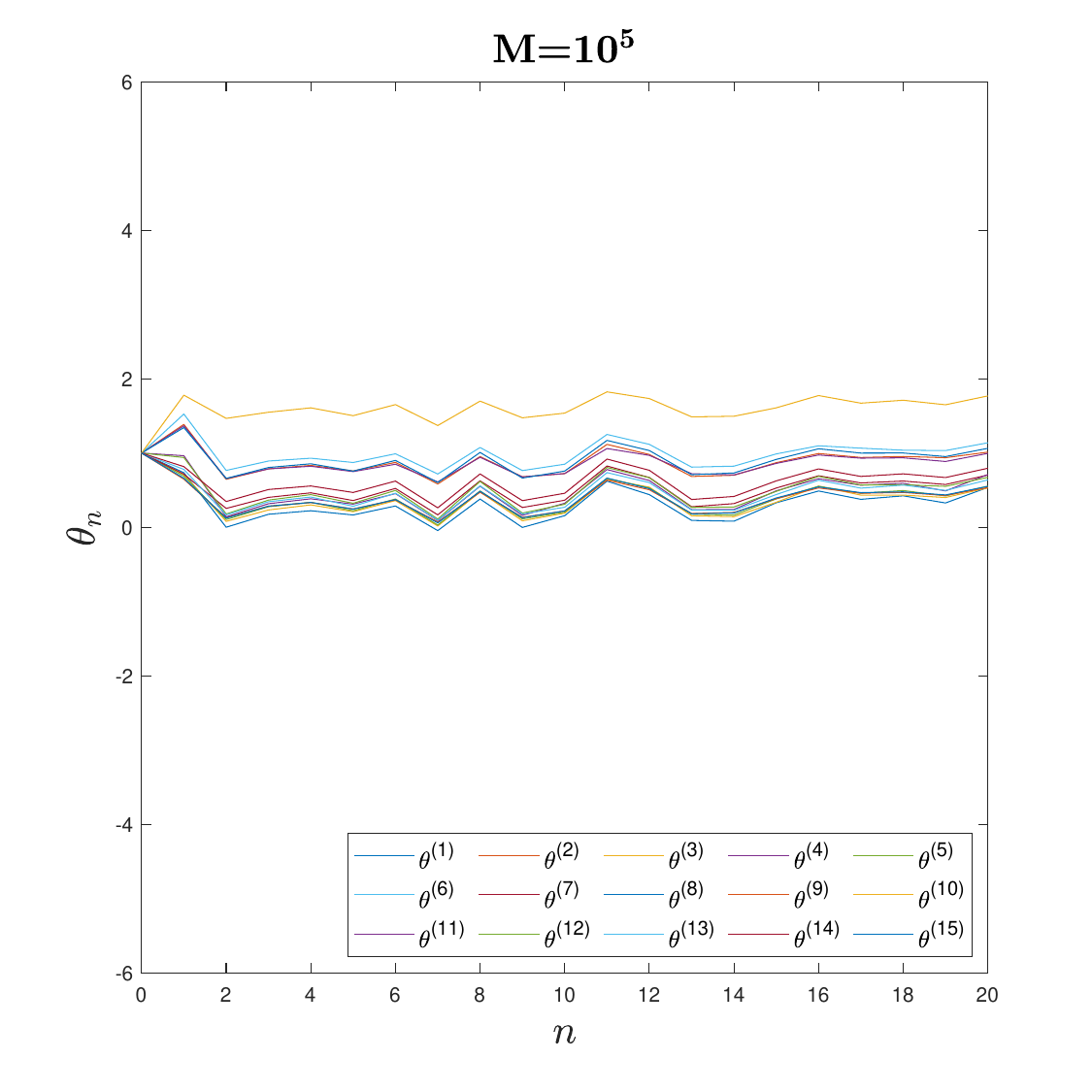}}
&
{\includegraphics[scale=0.5]{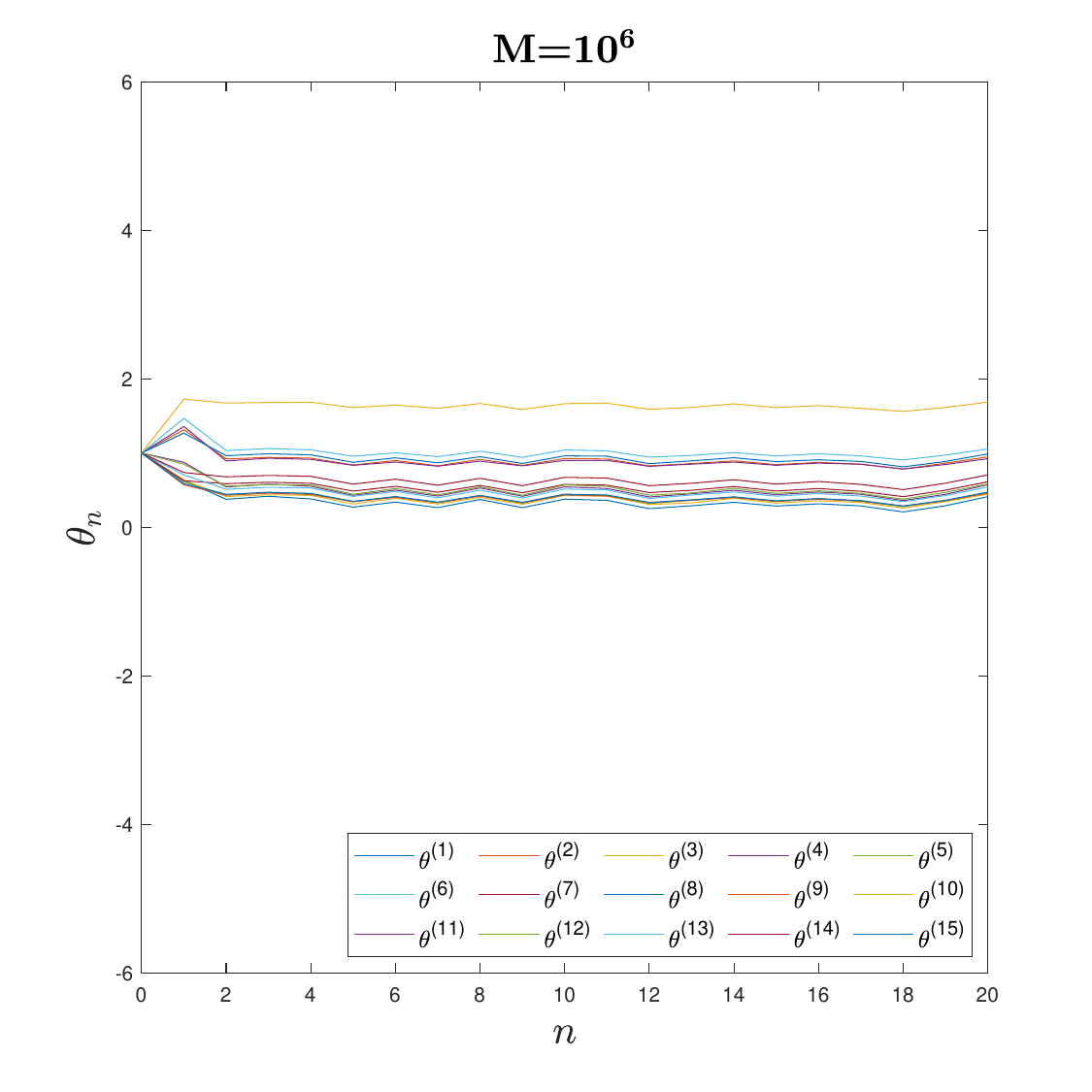}}
	\end{tabular}
	\caption[Values of weights $\btheta$ for Example 2.]{Values of $\theta^{(f)}_n$ in iteration $n$ of Algorithm~\ref{Al1} in Example~\ref{ex2}, for $M= 10^3, 10^4, 10^5$ and $10^6$.
		}
	\label{ex2_fig_theta}
\end{figure}

\begin{figure}[H]
	\begin{tabular}{cc}
{\includegraphics[scale=0.5]{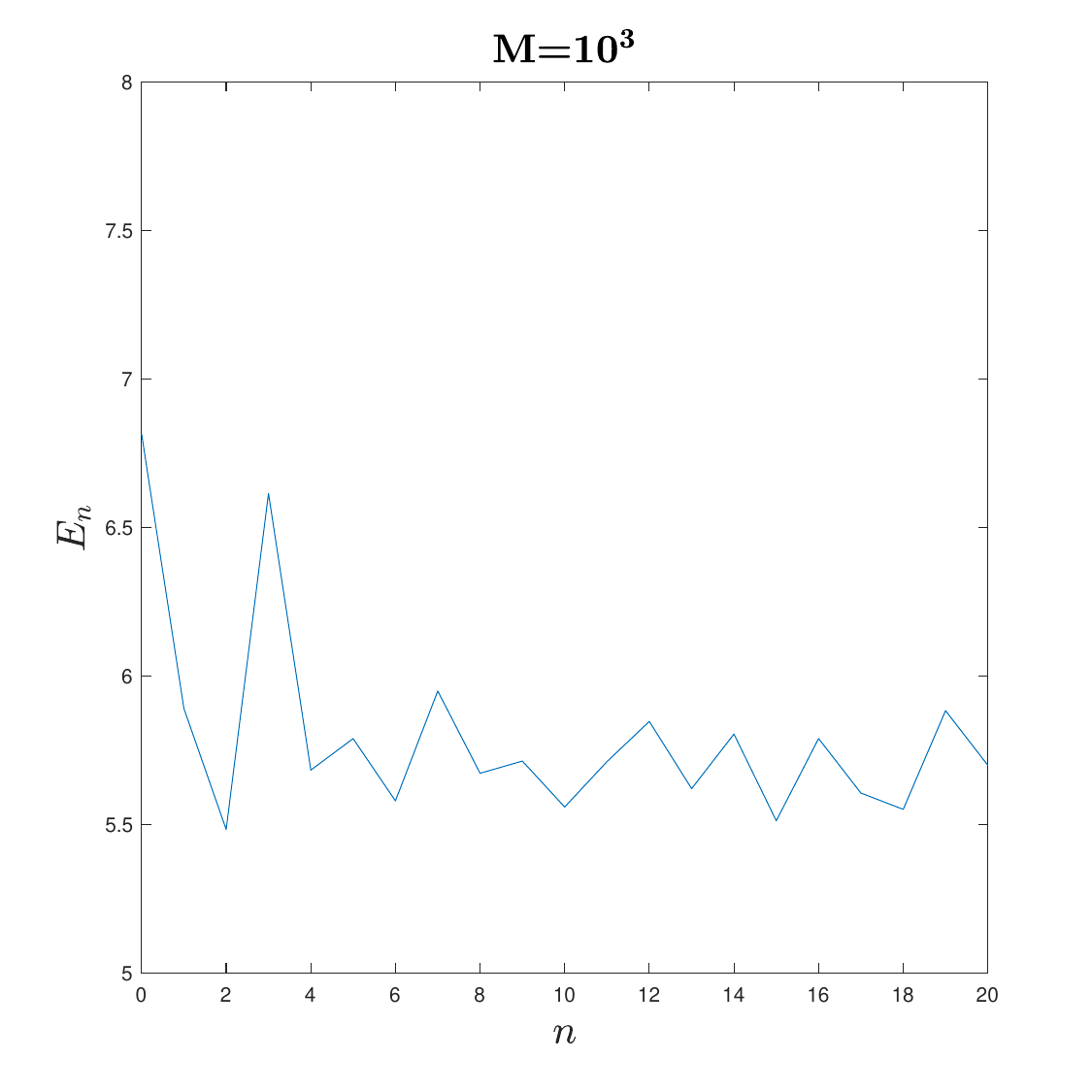}} &
{\includegraphics[scale=0.5]{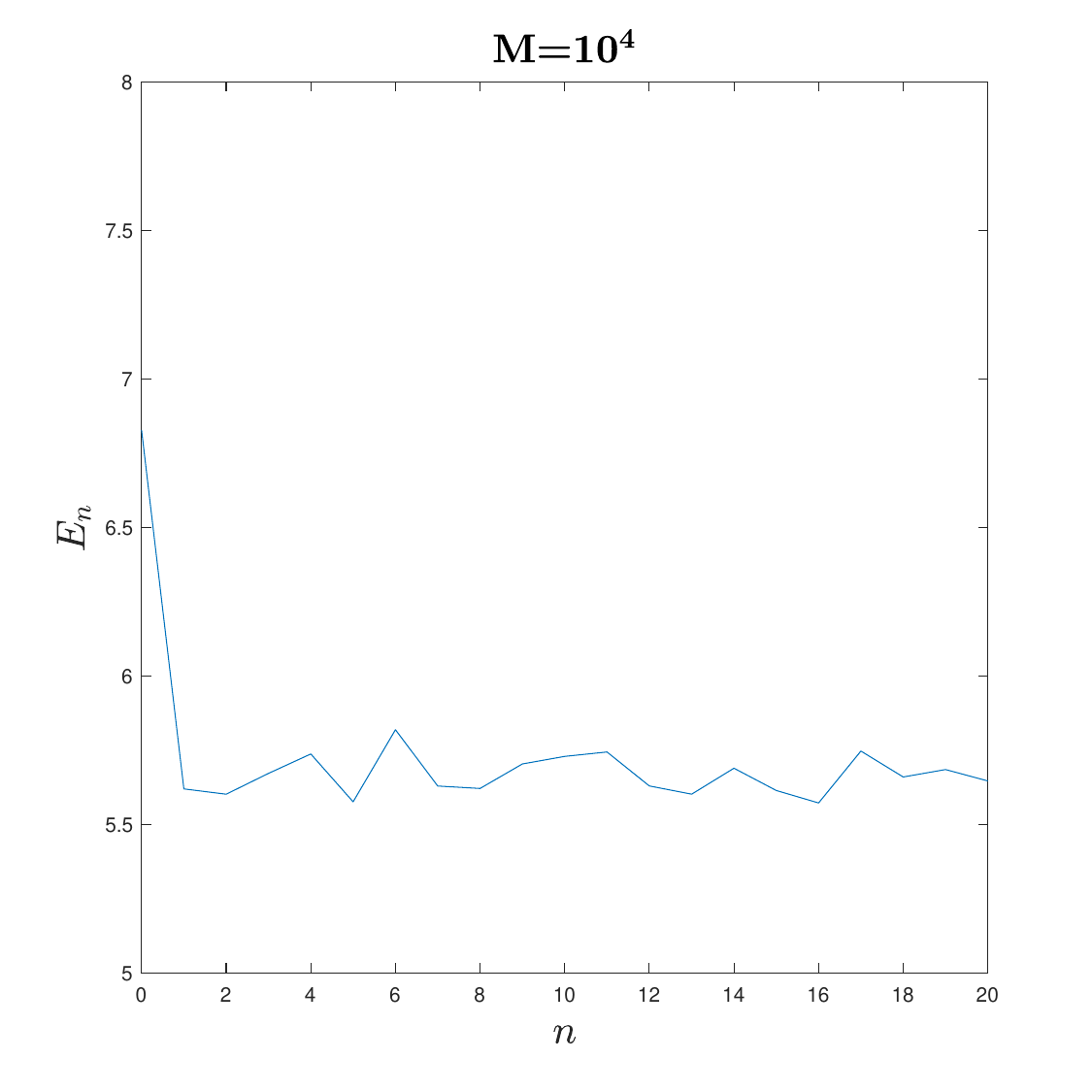}}\\
{\includegraphics[scale=0.5]{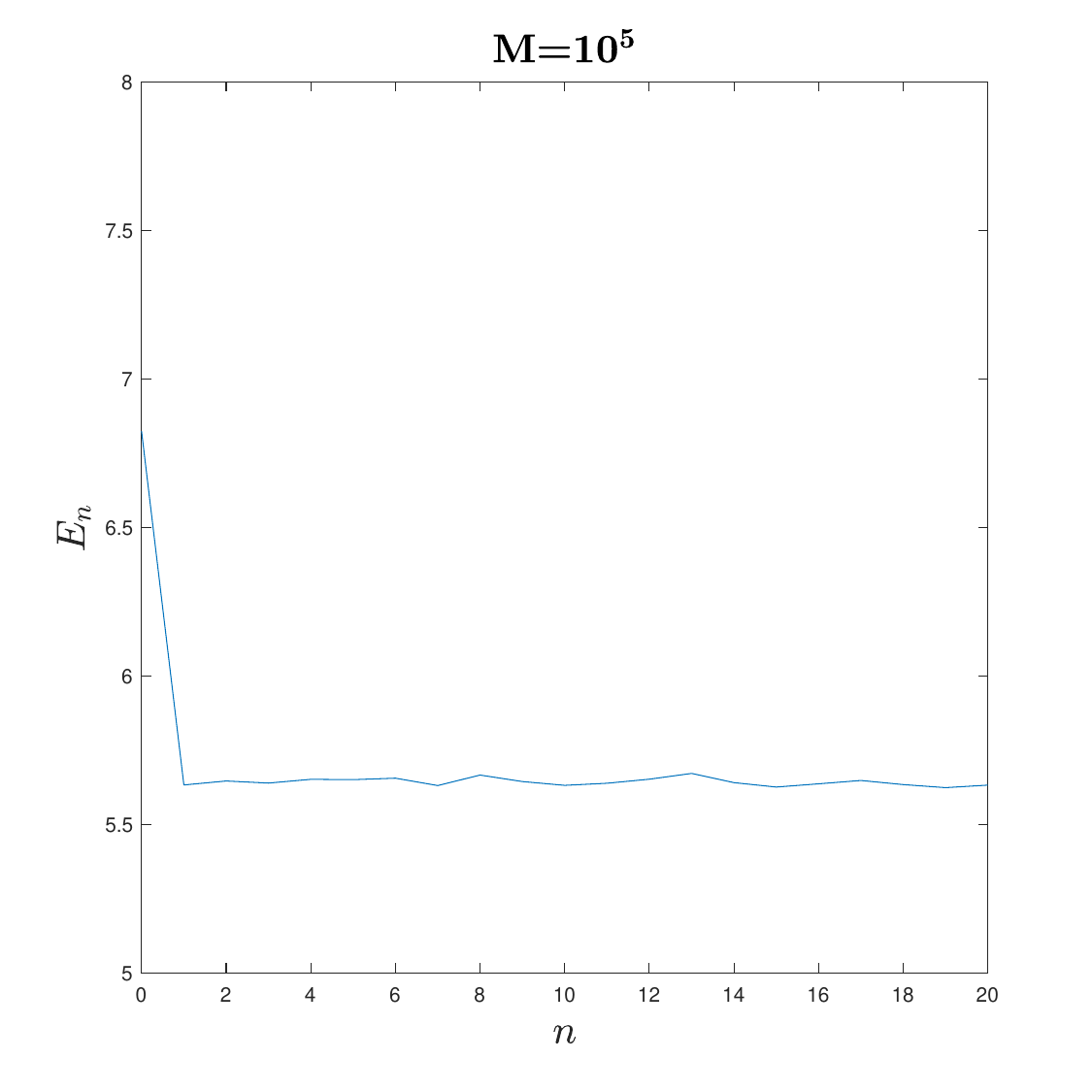}} &
{\includegraphics[scale=0.5]{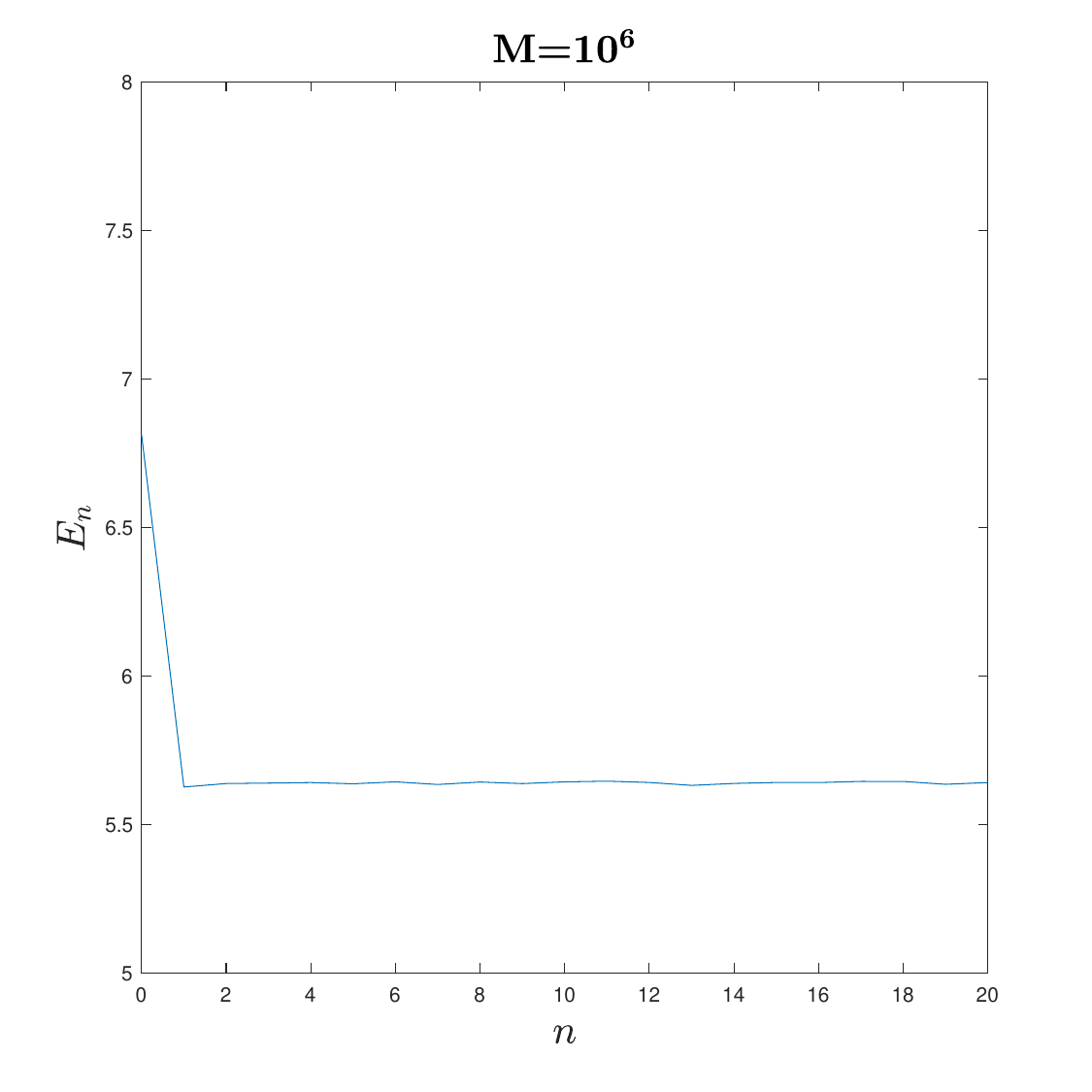}}
	\end{tabular}
	\caption[Values of $E$ for Example 2 compared to the mean value $E^*=xxx$.]{Values of $E_n$ in iteration $n$ of Algorithm~\ref{Al1}  in Example~\ref{ex2}, for $M= 10^3, 10^4, 10^5$ and $10^6$. 
		}
	\label{ex2_fig_E}
\end{figure}

Then, by~\eqref{eq:expected_features}, $\mathbb{E}(Q)=\lambda$, and given state $s$ and decision $a$, we have
\begin{eqnarray*}
	\sum_{s^{'}\in\mathcal{S}} \mathbb{P} \left( s^{\prime}\ \vert\  (s,a) \right) \bphi(s^{'}) \btheta
	&=&
	\mathbb{E}(\bphi(s^{'}) \btheta \ \vert\  (s,a))
	=
	\sum_{k=1}^{K}\sum_{i=1}^{I}
	n^{(a)}_{k,i}(1-p_{k,i})\theta_{k,i}
	+
	\sum_{i=1}^{I}
	\lambda_i\theta_i
	,
\end{eqnarray*}
which is the expression we use in our computatations in Algorithms~\ref{Al1}--\ref{Al1E}. The output presented in Figures~\ref{ex2_fig_theta}-\ref{ex2_fig_E} demonstrates that the algorithm converged as expected.

To analyse the impact of the algorithm on the performance of the system, we have then also simulated the evolution of a system in which the decision $a(s)$ given current state $s$ is chosen according to
\begin{eqnarray}\label{eq:approx}
a(s) = \argmin_{a\in\mathcal{A}(s)} \{ C(s,a)
+
\sum_{k=1}^{K}\sum_{i=1}^{I}
n^{(a)}_{k,i}(1-p_{k,i})\theta_{k,i}
+
\sum_{i=1}^{I}
\lambda_i\theta_i
,
\end{eqnarray}
using the parameters $\btheta$ obtained from the output of the algorithm, and compared the output with the performance of the system under decisions $a=1,2,3$. Each simulation $M=1,\ldots,1000$, was executed over a 5-year period. The output is presented in Figures~\ref{ex2_bounded_algo_compare_many}-\ref{ex2_bounded_algo_compare_2_withzeros} and Table~\ref{tab:ex2_output} below.

\begin{table}[H]
\centering
\footnotesize
\begin{tabular}{l|r r r }
$\textbf{Policy}$&
$\textbf{Mean cost}$&
$\textbf{Mean nonprimary}$&
$\textbf{Mean redirected}$
\\\hline
$\mbox{Apply}$ a=1 $\mbox{always}$&
6.8177&34.0887&6.6719\\
$\mbox{Apply}$ a=2 $\mbox{always}$&
5.7643&20.1370&5.9918\\
$\mbox{Apply}$ a=3 $\mbox{always}$&
5.6983&18.9862&5.9408\\\hline
$\textbf{Near optimal solution}$&
{\bf 5.6449}&{\bf 20.2566}&
{\bf 5.9720}\\
\end{tabular}
\caption{\textit{The output from simulations in Example~\ref{ex2}.}}
\label{tab:ex2_output}
\end{table}

We observe (Figure~\ref{ex2_bounded_algo_compare_many}, Table~\ref{tab:ex2_output}) that under the near-optimal policy obtained by the application of this approach, the costs are minimised (costs are the lowest under near-optimal policy, and are the highest under $a=1$), with mean costs: $6.8177$ ($a=1$), $5.7643$ ($a=2$), $5.6983$ ($a=3$), and ${\bf 5.6449}$ (our near optimal solution). Simultaneously, the performance of the system in the sense of reducing the number of patients observed in the nonprimary wards and the number of redirected patients per day, is also improved. The number of patients in nonprimary wards is significantly lower than under $a=1$, and similar to $a=2$ but with lower costs than under $a=2,3$.

Next (Figure~\ref{ex2_bounded_algo_compare_2_withzeros}), under near-optimal solution, $82.1918\%$ of the time no additional transfers are required. Also, under near-optimal solution, decisions $a=1,2,3$ are chosen approximately $50\%$, $40\%$ and $10\%$ of the time, respectively.

\begin{figure}[H]
	\centering
	$a=1$:
	\\
	{\includegraphics[scale=0.4]{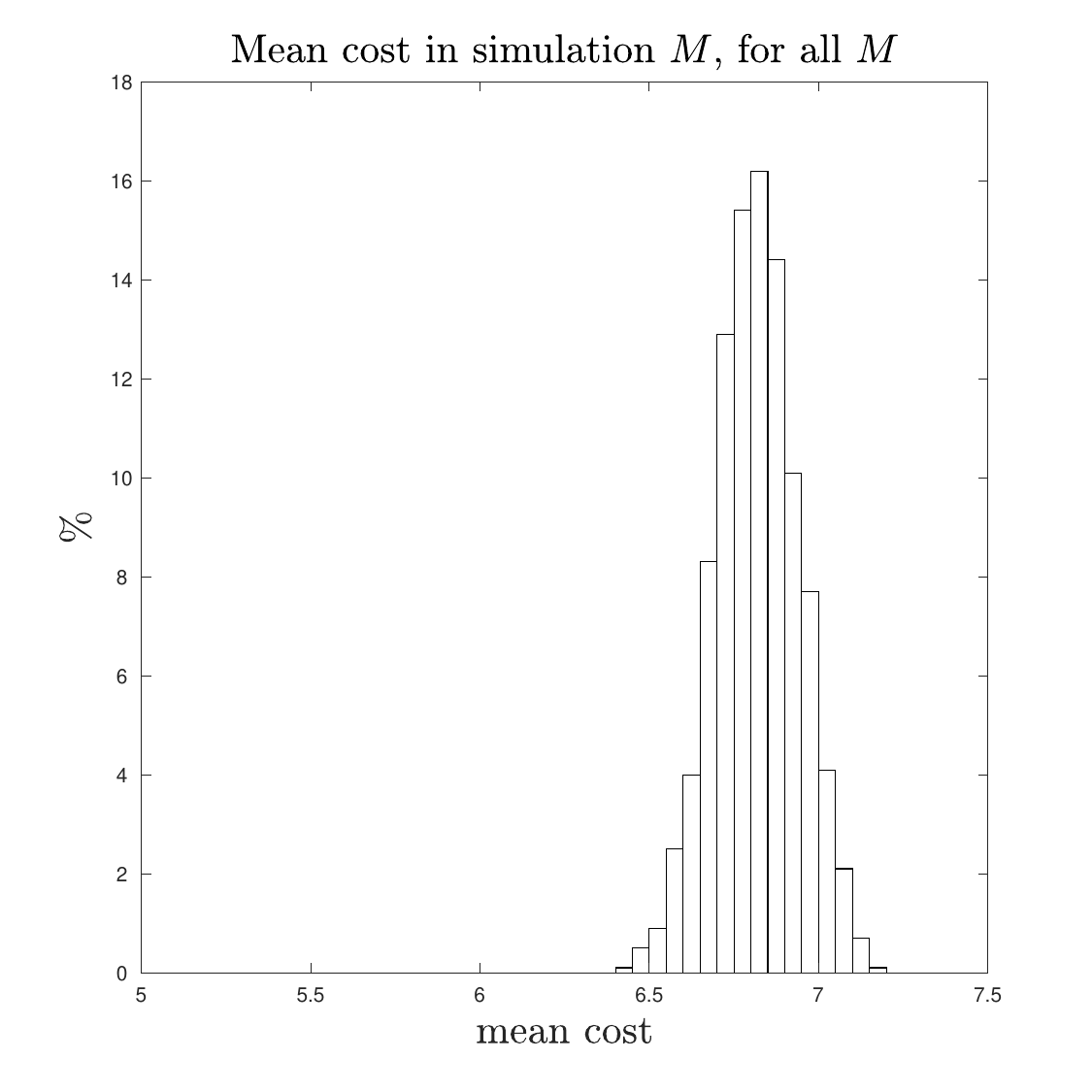}}
	\
	{\includegraphics[scale=0.4]{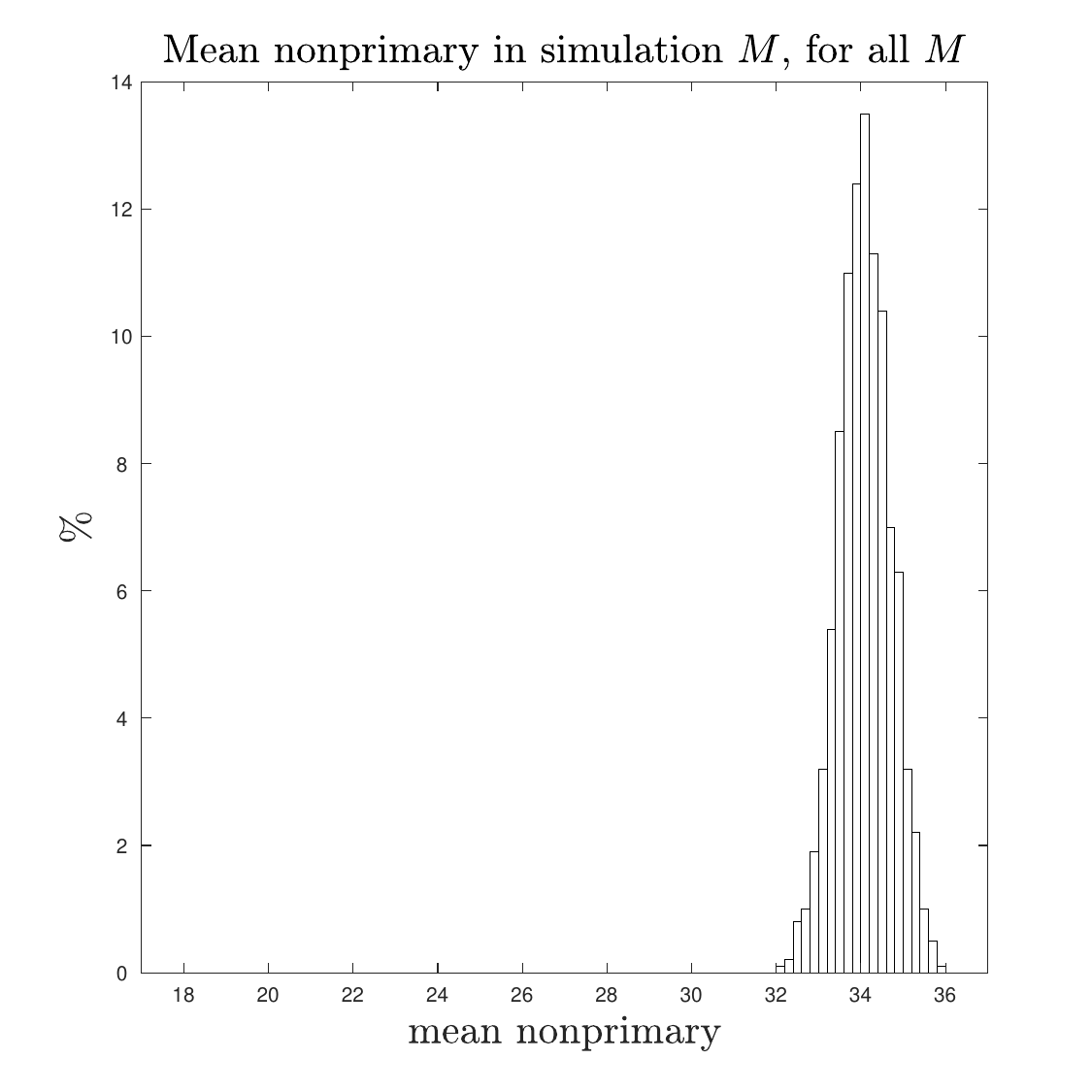}}	
	\\
	$a=3$:
	\\
	{\includegraphics[scale=0.4]{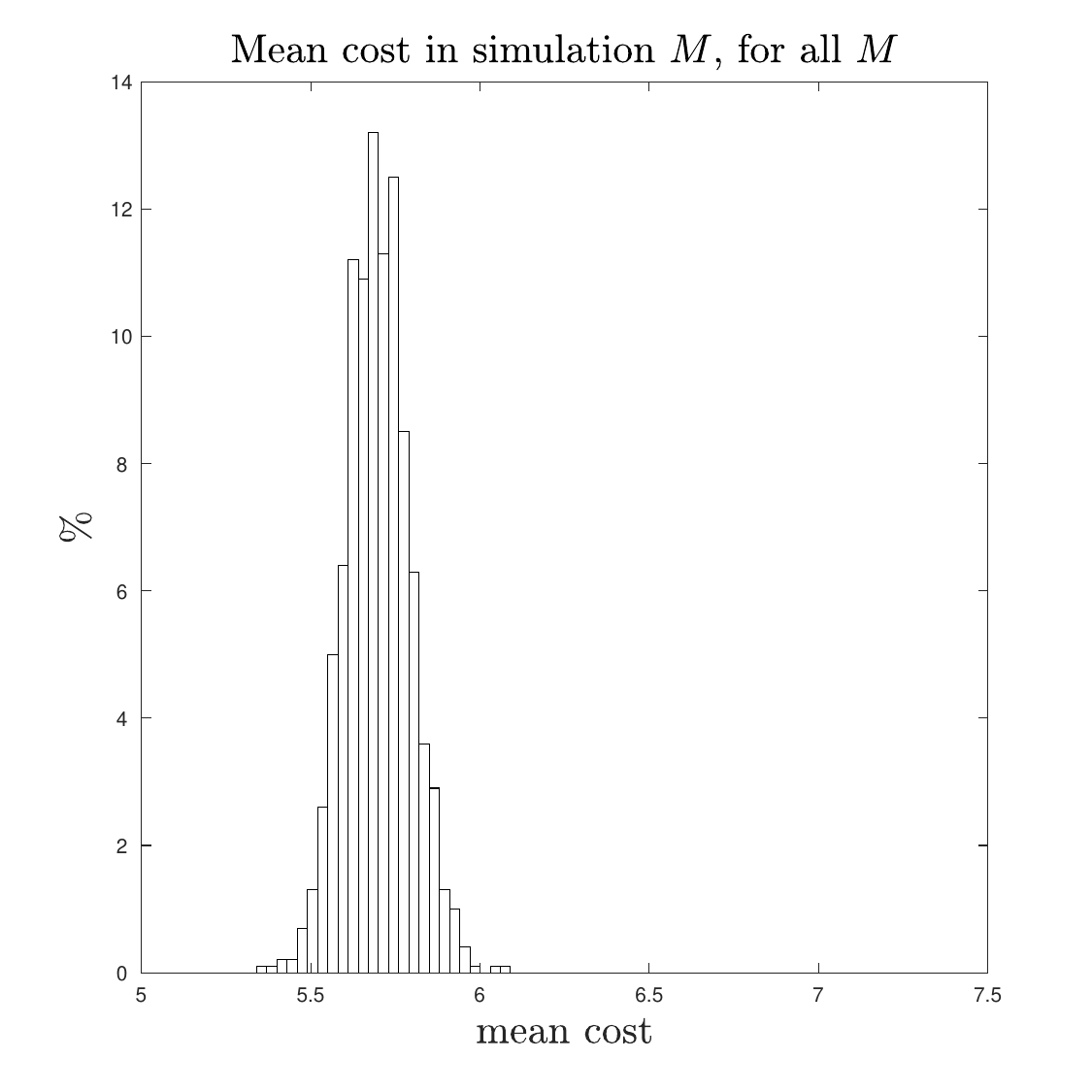}}
	\
	{\includegraphics[scale=0.4]{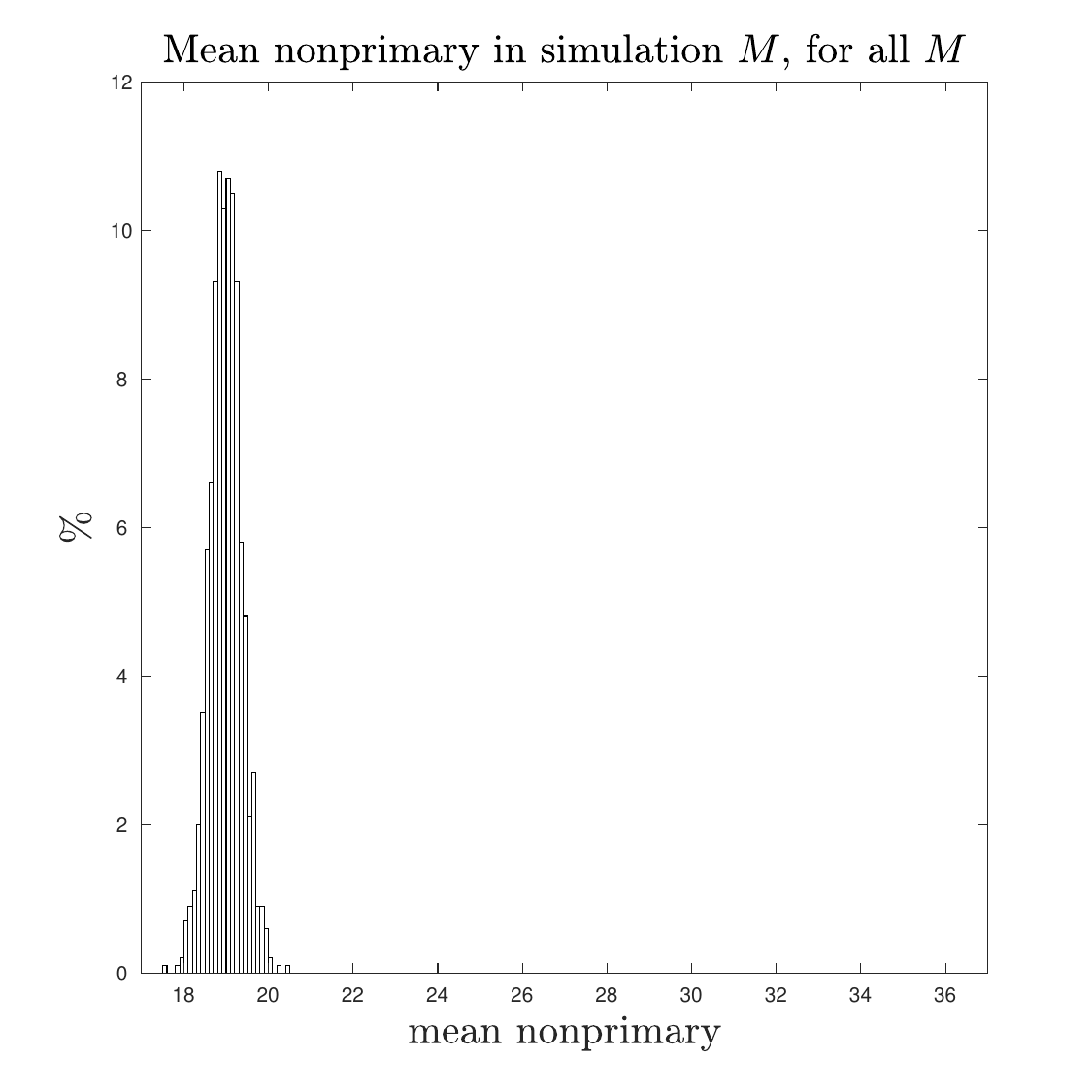}}	
	\\
	near-optimal solution obtained using~\eqref{eq:approx}:
	\\
	{\includegraphics[scale=0.4]{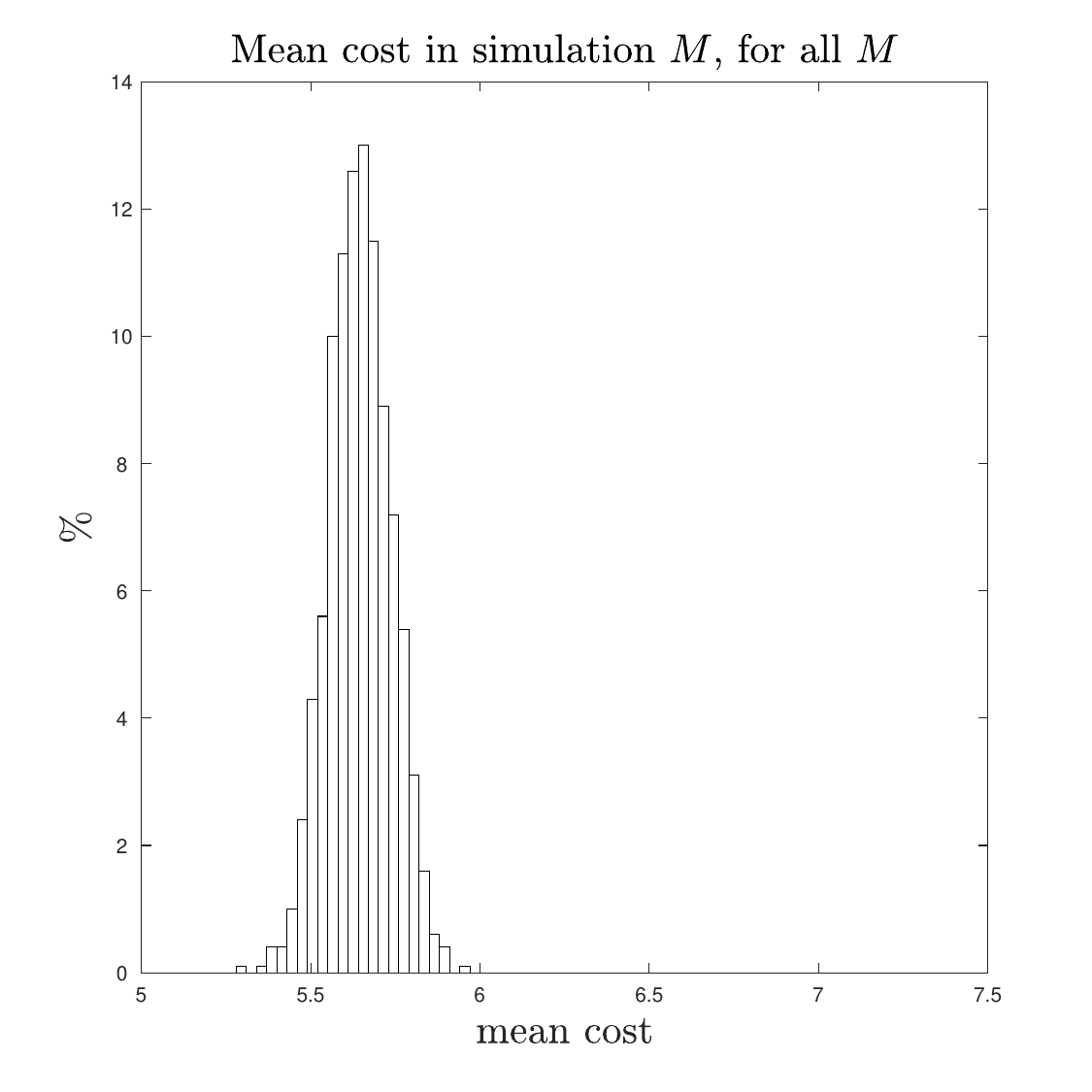}}
	\
	{\includegraphics[scale=0.4]{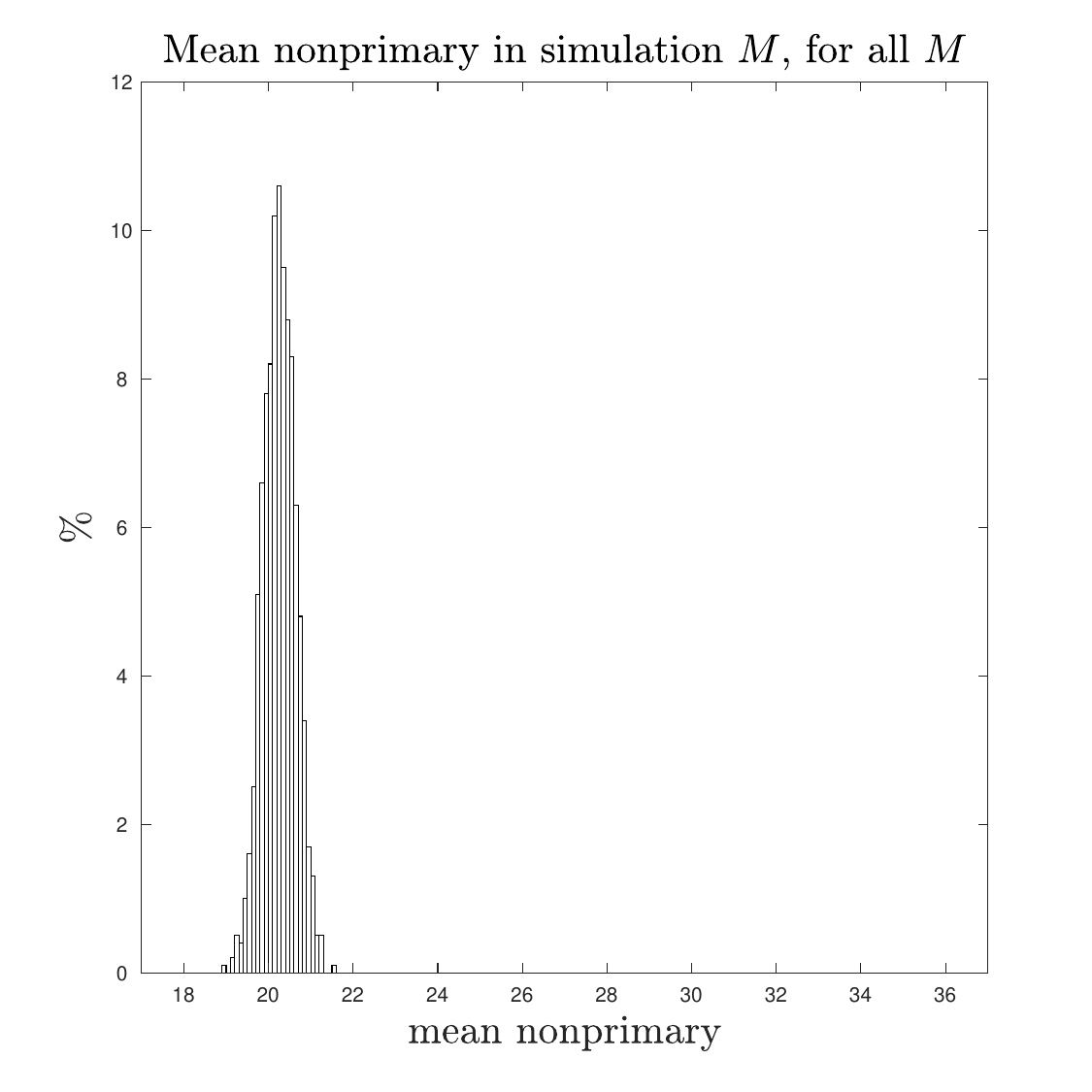}}		
	\\
	\caption{Simulations under decisions $a=1$ (top row), $a=3$ (middle row), and near-optimal solution.}
	\label{ex2_bounded_algo_compare_many}
\end{figure}

\begin{figure}[H]
	\centering
	{\includegraphics[scale=0.4]{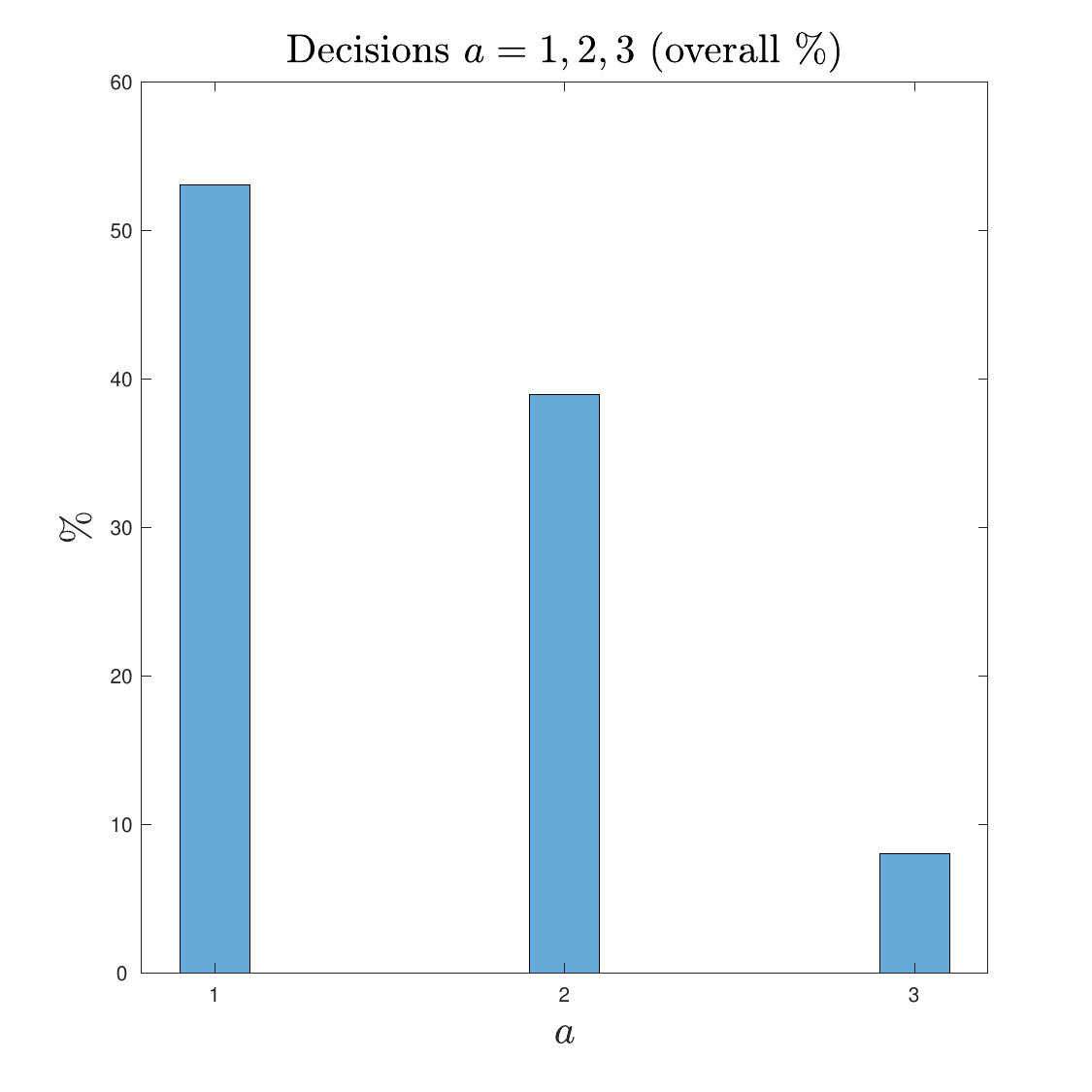}}	
	\\
	{\includegraphics[scale=0.4]{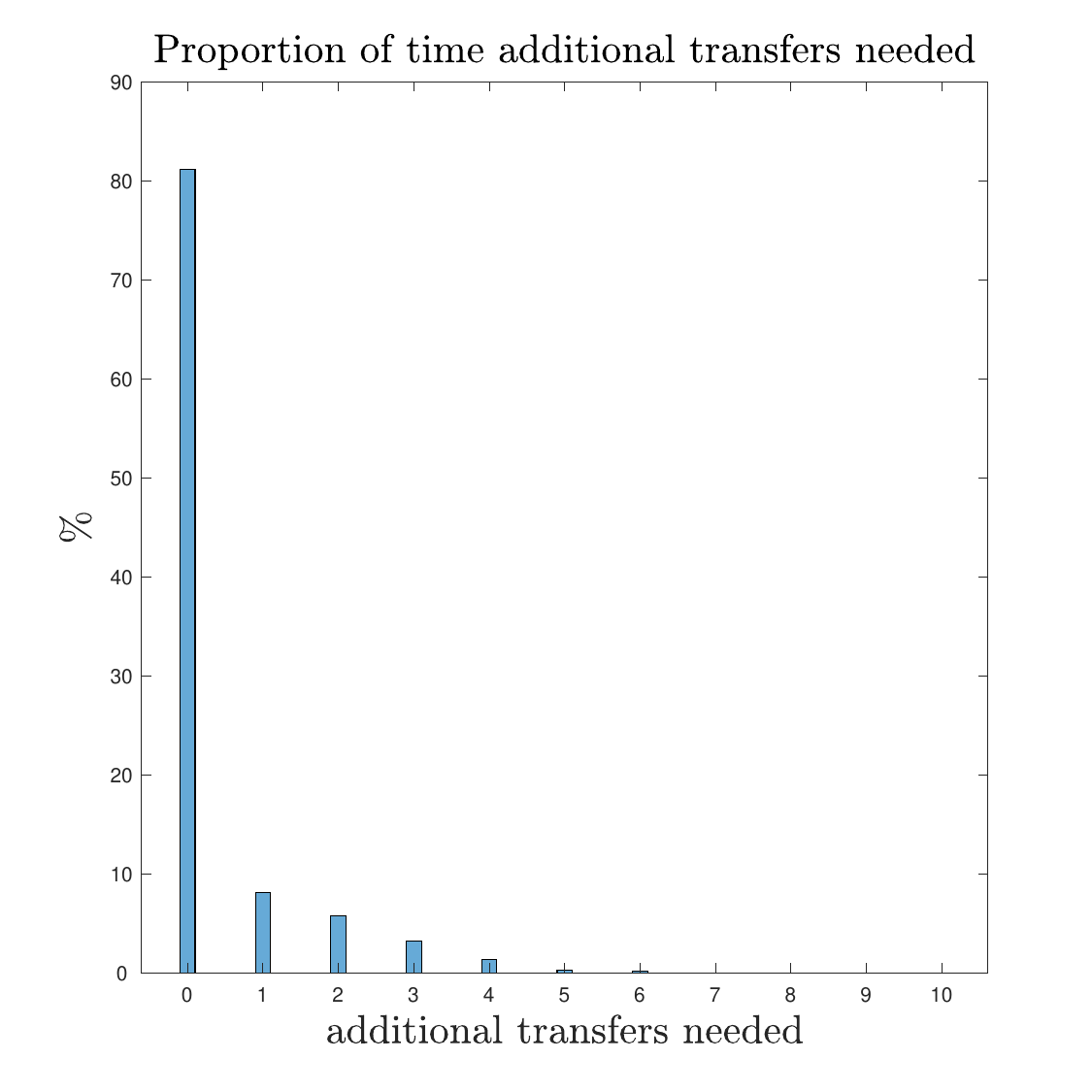}}	
	\
	{\includegraphics[scale=0.4]{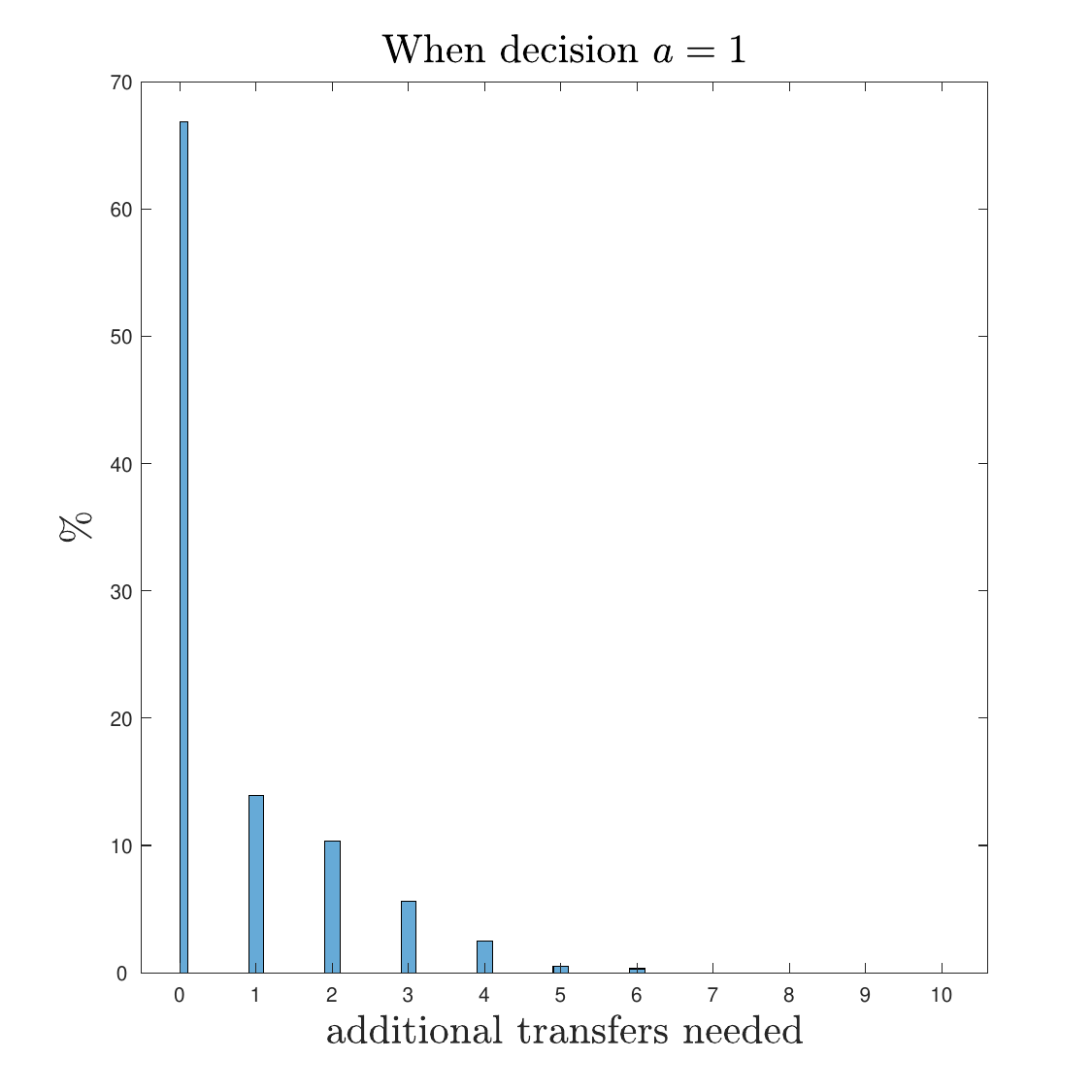}}
	\\	
	{\includegraphics[scale=0.4]{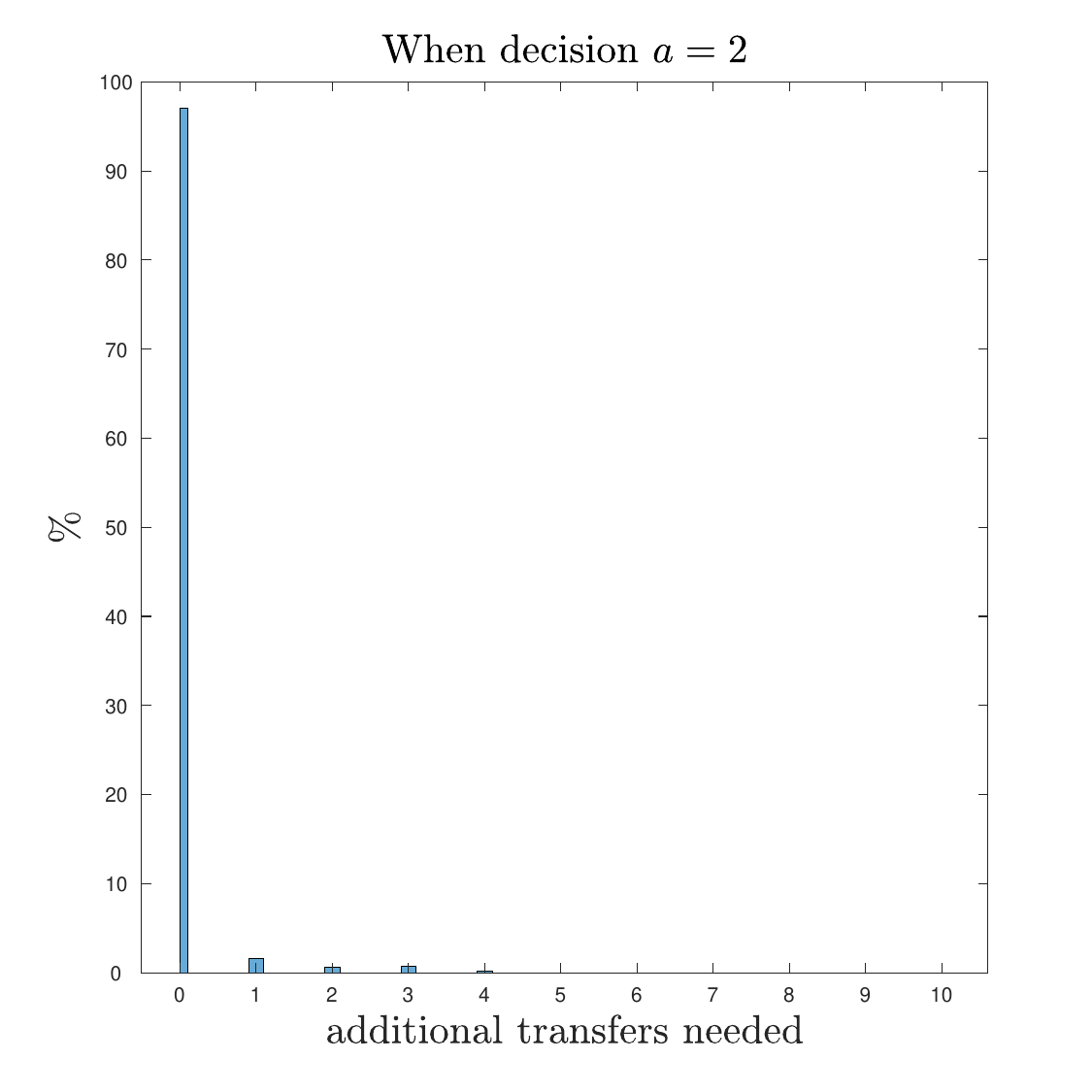}}
	\
	{\includegraphics[scale=0.4]{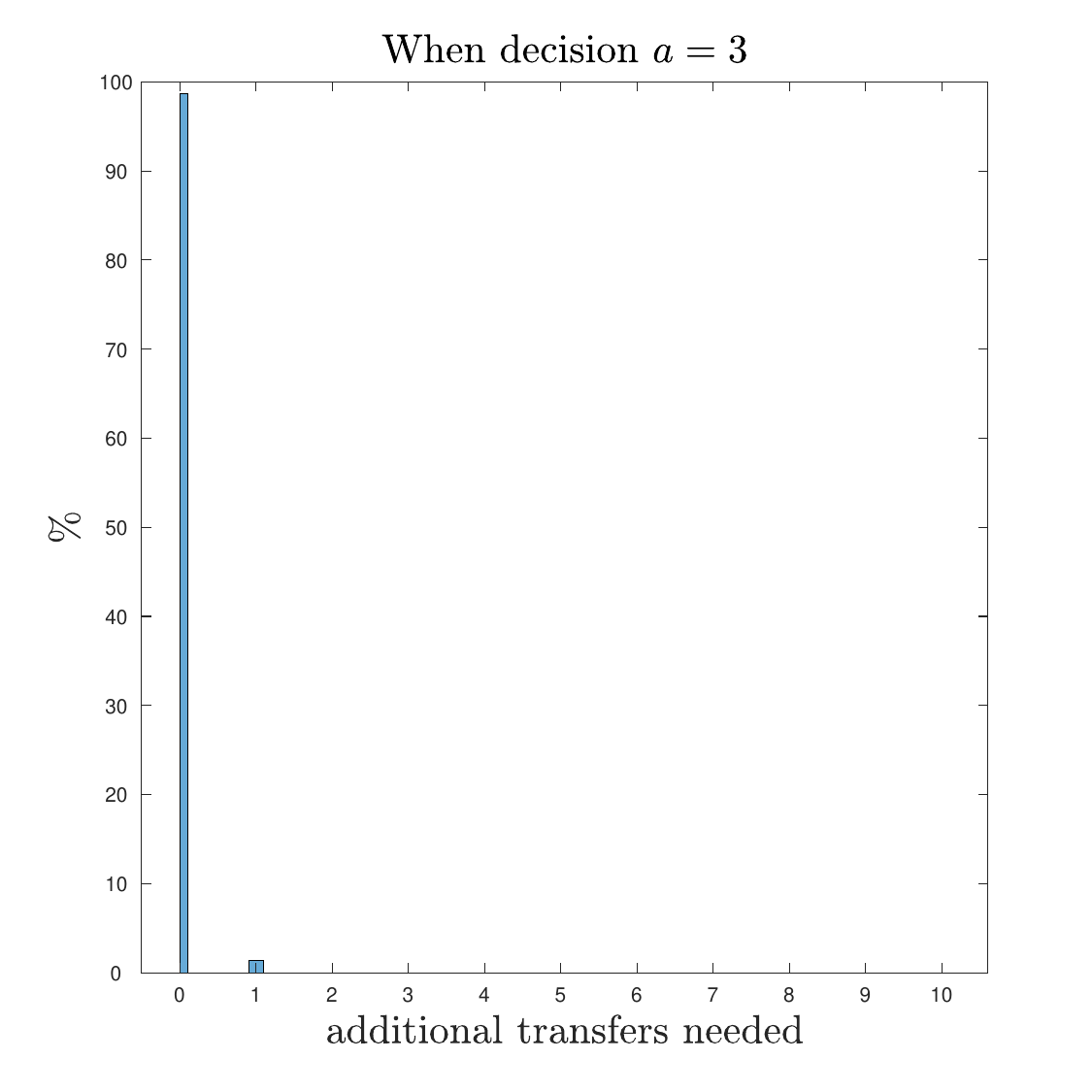}}
	\\	
	\caption{Percentage of decisions $a=1,2,3$ made under near-optimal solution (top left); and the number of additional transfers needed (if they were permitted without constraints).}
	\label{ex2_bounded_algo_compare_2_withzeros}
\end{figure}

\section{Conclusion}\label{sec:conclusion}

We proposed a Markov Decision Process for a Patient Assignment Scheduling problem where, at the start of each time period (e.g. day, 8-hour block etc), patients are allocated to suitable wards, following relevant hospital policy. Due to the large size of the state space required to record the information about the hospital system, we developed methodology based on Approximate Dynamic Programming, for the computation of near-optimal solutions. We demonstrated the application potential of our model through examples with parameters fitted to data from a tertiary referral hospital in Australia.

Our model can be applied to a range of scenarios that are of interest in practical contexts in hospitals. We will analyse scheduling problems in data-driven examples, and consider additional costs, such as waiting costs or redirecting costs, and decisions such as transfers of inpatients within the hospitals that are not triggerred by arrivals, and take into account the duration of stay at the time of decision making. We will also consider scenarios with time varying parameters, in particular with seasonal components. These analyses will be reported in our forthcoming and future work.

\section{Statements and declarations}

%
%

\bigskip\noindent{\bf Data}\\

Data used in the paper was obtained following ethical approval from the Tasmanian Health and Medical Human Research Ethics Committee (HREC No 23633) and site-specific approval from the Research Governance Office of the Tasmanian Health Service.

\bigskip\noindent{\bf Authorship contribution statement}\\

This research has contributed to the Honours thesis by Krasnicki~\cite{Krasnicki}. The following are the contributions of the authors, Ma\l gorzata M. O'Reilly (MMO), Sebastian Krasnicki (SK), James Montgomery (JM), Mojtaba Heydar (MH), Richard Turner (RT), Pieter Van Dam (PVD), Peter Maree (PM): 
\begin{itemize}
	\item	Conceptualisation, Mathematical background: MMO, SK, MH, and JM;
	\item  
	Conceptualisation, Clinical background: RT, PVD, PM;
	\item Sections~\ref{sec:MDPmodel}-\ref{sec:numex}, Appendix~\ref{sec:algorithms}-\ref{sec:ex1matrices}, Methodology development: MMO, SK, and JM;
	\item Sections~\ref{sec:MDPmodel}-\ref{sec:numex}, Derivation of the probability formulae: MMO;
	\item Section~\ref{sec:numex}, Example~\ref{ex1}, Coding and analysis: SK, MMO, and MH;
	\item Section~\ref{sec:numex}, Example~\ref{ex2}, Coding and analysis: MMO and JM.
\end{itemize}

{\bf Declaration of competing interests}\\

The authors have no competing interests to declare that are relevant to the content of this article.

\bibliographystyle{abbrv}
\bibliography{PWA_model}

\appendix
\section{Appendix: Algorithms}\label{sec:algorithms}

\begin{algorithm}
	\caption{Policy Iteration (Howard~\cite{1960_Howard})}
	\label{Al_pol_it}
	\begin{algorithmic}[1] %
		\State Initialise an aribitrary policy $\widehat\bpi=(a^{\widehat\bpi}(s))_{s\in\mathcal{S}}$.
		\State Let $v(m)=0$, where $\mathcal{S}=\{1,\ldots,m\}$.
		\State Solve $E+v(s)-\sum_{s' \in \mathcal{S}} \mathbb{P} \left( s^{\prime}\ \vert\  (s,a^{\widehat\bpi}(s)) \right) v(s')=C(s,a^{\widehat\bpi}(s))$ for all $s=1,\ldots,m-1$.
		\State Solve $\bpi=\arg\min_{a \in \mathcal{A}(s)} \Bigg\{ C(s,a) +  
		\sum_{s' \in \mathcal{S}} \mathbb{P} \left( s^{\prime}\ \vert\  (s,a) \right) v(s') \Bigg\}$.
		\State If $\bpi= \widehat\bpi$, then go to Step 6. Else, let $\widehat\bpi= \bpi$ and go to Step 2.
		\State The optimal policy is $\bpi^*=\bpi$.
	\end{algorithmic}
\end{algorithm}

\begin{algorithm}
	\caption{(Adapted from Dai and Shi~\cite{Dai2019})}
	\label{Al1}
	\begin{algorithmic}[1] %
		\State Initialise $\btheta = \btheta_0$, $n=0$. Choose large $N$.
		\State Initialise $E_0$ using Algorithm~\ref{Al1E} for $E$ with $\btheta = \btheta_0$.
		\For{$n=0,1,2,\ldots,N-1$}
		\State Randomly initialise starting state $s_0$. Let $\boldsymbol{A}={\bf O}$, $\boldsymbol{b}=\boldsymbol{0}$.
		\For{$m=0,1,2,\ldots,M-1$}
		\State Given $s=s_m$, take decision $a(s) = \argmin_{a\in\mathcal{A}(s)} \{ C(s,a)+\sum_{s^{'}\in\mathcal{S}} \mathbb{P} \left( s^{\prime}\ \vert\  (s,a) \right) \bphi(s^{'}) \btheta \}$.
		\State Randomly generate state $s_{m+1}$ using model parameters and post-decision state $s_m^{(a)}$.
		\State Compute $\boldsymbol{A}_m =  \frac{1}{M} \bphi(s_m)^T(\bphi(s_m)-\bphi(s_{m+1}))$.
		\State Compute $\boldsymbol{b}_m = \frac{1}{M} \bphi(s_m)^T(C(s_m , a(s_m))-E_n)$.
		\State Let $\boldsymbol{A}=\boldsymbol{A}+\boldsymbol{A}_m$.
		\State Let $\boldsymbol{b}=\boldsymbol{b}+\boldsymbol{b}_m$.
		\EndFor
		\State Solve the set of linear equations 
		$\boldsymbol{A} \btheta_{n+1} = \boldsymbol{b}$ for $\btheta_{n+1}$.
		\State Update $\btheta=\btheta_{n+1}$.
		\State Let $n=n+1$.
		\State Compute $E_n$ using Algorithm~\ref{Al1E} for $E$ with $\btheta = \btheta_n$.
		\EndFor
		\State $E^{\bpi^*}\approx E_N$.
	\end{algorithmic}
\end{algorithm}

\begin{algorithm}
	\caption{(Adapted from Dai and Shi~\cite{Dai2019})}
	\label{Al1E}
	\begin{algorithmic}[1] %
		\State Input $\btheta$. Let $E=0$.
		\State Randomly initialise starting state $s_0$.
		\For{$m=0,1,2,\ldots,M-1$}
		\State Given state $s=s_m$, take decision $a(s) = \argmin_{a\in\mathcal{A}(s)} \{ C(s,a)+\sum_{s^{'}\in\mathcal{S}} \mathbb{P} \left( s^{\prime}\ \vert\  (s,a) \right) \bphi(s^{'}) \btheta \}$.
		\State Randomly generate state $s_{m+1}$ using model parameters and post-decision state $s_m^{a(s)}$.
		\State Let $E=E+C(s,a(s))$.
		\EndFor 
		\State Compute $E=E/M$.
	\end{algorithmic}
\end{algorithm}

\begin{algorithm}
	\caption{Decision $a=1$. We make a decision based on priorities so as to guarantee the best possible bed for higher priority patients, without allowing transfers. Recursively, we find the highest priority of patients to allocate yet and transfer them to their best available wards.}
	\label{a4}
	\begin{algorithmic}[1] %
		\State Input $s=\left( [n_{k,i}]_{\mathcal{K} \times \mathcal{I}}, [q_i]_{1 \times \mathcal{I}} \right)$.
		\State Initialise $s^{(a)}= [n_{k,i}^{(a)}]_{\mathcal{K} \times \mathcal{I}} = [n_{k,i}]_{\mathcal{K} \times \mathcal{I}}$ and $[x_{k,i}]_{\mathcal{K} \times \mathcal{I}}={\bf 0}$.
		\While{$\sum_{i=1}^I q_i>0$}
		\State Let $i=\min\{j:q_j>0\}$, current highest priority to allocate.
		\While{$q_i>0$} (allocate type-$i$ patients)
		\State Find $\ell^* = \min \{\ell: \sum_{j=1}^I n_{w(\ell,i),j}^{(a)}< m_{w(\ell,i)} \}$. 
		\State Let $k^*=w(\ell^*,i)$, best ward for type-$i$ patient with an available bed.
		\State Let $x_{k^*,i}=x_{k^*,i}+1$, and $n_{k^*,i}^{(a)}=n_{k^*,i}^{(a)}+1$.
		\State Let $q_i=q_i -1$.
		\EndWhile 
		\EndWhile
	\end{algorithmic}
\end{algorithm}

\begin{algorithm}
	\caption{Decisions $a=2,3$.  We make a decision based on priorities so as to guarantee the best possible bed for higher priority patients, and allow up to $y$ transfers. Recursively, we find the highest priority of patients to allocate/transfer yet (Line 6). Next, we transfer type-$i$ patients (Lines 7-13) before allocating newly arrived type-$i$ patients (Lines 14-30). When transferring patients, we find a type-$i$ patient that is in the worst ward and transfer them to a ward with the lowest transfer cost.}
	\label{a5}
	\begin{algorithmic}[1] %
		\State Input $s=\left( [n_{k,i}]_{\mathcal{K} \times \mathcal{I}}, [q_i]_{1 \times \mathcal{I}} \right)$ and the maximum number of transfers $y$.
		\State Initialise $y^*=y$, current number of available transfers.
		\State Initialise $[y_{k,i}]_{\mathcal{K} \times \mathcal{I}}={\bf 0}$, patients to be transferred.
		\State Initialise $s^{(a)}= [n_{k,i}^{(a)}]_{\mathcal{K} \times \mathcal{I}} = [n_{k,i}]_{\mathcal{K} \times \mathcal{I}}$ and $(x_{k,i},y_{k,\ell,i})_{i\in\mathcal{I};k,\ell=1,\ldots K}={\bf 0}$.
		\While{$\sum_{i=1}^I q_i>0$ or $\sum_{i=1}^I\sum_{k=1}^K y_{k,i}>0$}
		\State Let $i=\min\{j:q_j>0 \mbox{ or }\sum_{k=1}^K y_{k,j}>0\}$, current highest priority to allocate/transfer.
		\While{$\sum_{k=1}^K y_{k,i}>0$} (transfer type-$i$ patient from ward $k^*$ to $\ell^*$)
		\State  $h^*=\max\{\ell:y_{w(\ell,i),i}>0\}$;
		\State  $k^*=w(h^*,i)$ is the worst ward with a type-$i$ patient;
		\State $\ell^*=\arg\min_{\ell=1,\ldots,K}\{c^{(t)}_{k^*,\ell,i}: \sum_{j=1}^I n_{\ell,j}^{(a)}<m_{\ell}\}$ has the lowest transfer cost.
		\State Let $n_{\ell^*,i}^{(a)}=n_{\ell^*,i}^{(a)}+1$.
		\State Let $y_{k^*,\ell^*,i}=y_{k^*,\ell^*,i}+1$, $y_{k^*,i}=y_{k^*,i}-1$.
		\EndWhile
		\While{$q_i>0$} (allocate type-$i$ patients)
		\If{$y^*=0$ (no more transfers allowed)} (allocate type-$i$ patient)
		\State $\ell^* = \min \{\ell: \sum_{j=1}^I n_{w(\ell,i),j}^{(a)}< m_{w(\ell,i)} \}$; 
		\State $k^*=w(\ell^*,i)$ is the best ward for patient type $i$ with an available bed.
		\State Let $x_{k^*,i}=x_{k^*,i}+1$, and $n_{k^*,i}^{(a)}=n_{k^*,i}^{(a)}+1$.
		\EndIf
		\If{$y^*>0$ (transfers allowed)} (allocate type-$i$ patient to ward $k^*$)
		\State $\ell^* = \min \{ \ell:\sum_{j=1}^i n_{w(\ell,i),j}^{(a)}< m_{w(\ell,i)} \}$; 
		\State $k^*=w(\ell^*,i)$ is the best ward for patient type $i$ with a potential bed.
		\State Let $x_{k^*,i}=x_{k^*,i}+1$, and $n_{k^*,i}^{(a)}=n_{k^*,i}^{(a)}+1$.
		\If{$\sum_{j=1}^I n_{k^*,j}^{(a)}= m_{k^*}+1$} (must transfer patient $i^*$ out of ward $k^*$)
		\State $i^*=\max\{j:j=i+1,\ldots,I; n_{k^*,j}^{(a)}\not= 0\}$ is the lowest priority patient in ward $k^*$.
		\State Let $y_{k^*,i^*}=y_{k^*,i^*}+1$, $y^*=y^*-1$.
		\EndIf
		\EndIf
		\State Let $q_i=q_i -1$.	
		\EndWhile
		\EndWhile
	\end{algorithmic}
\end{algorithm}

\newpage
\section{Appendix: Example 1 -- matrices ${\bf P}^{(a)}$ and vectors ${\bf C}^{(a)}$}\label{sec:ex1matrices}

We note that the set of all possible post-decision states $s^{(a)}$ is $\{1,7,10,13,16,19,20,21,22\}$ for both $a=1,2$, with
\begin{eqnarray}
s^{(1)}=\left\{
\begin{array}{lll}
1 &\mbox{ for }&s=1\\
7 &\mbox{ for }&s=2,7\\
10 &\mbox{ for }&s=10\\
13 &\mbox{ for }&s=13\\
16 &\mbox{ for }&s=3,16\\
19 &\mbox{ for }&s=4,8,14,19\\
20 &\mbox{ for }&s= 11,15,20\\
21 &\mbox{ for }&s=5,9,17,21\\
22 &\mbox{ for }&s=6,12,18,22
\end{array}	
\right.
\quad
,
\quad
s^{(2)}=
\left\{
\begin{array}{lll}
s^{(1)} &\mbox{ for }&s\not= 11,15\\
21 &\mbox{ for }&s=11,15
\end{array}	
\right.
\quad ,
\end{eqnarray}
where the difference between $s^{(1)}$ and $s^{(2)}$ is in places corresponding to $s=11,15$, which are the only two states where choosing between $a=1$ or $a=2$ can make a difference.

The probability matrices ${\bf P}^{(a)}$ of the discrete-time Markov chains corresponding to decisions $a=1,2$ are evaluated as follows. First, denote $q_{k,i}=1-p_{k,i}$, the probability that type-$i$ patient does not depart from ward $k$. Next, denote $p(k)=\frac{\lambda^k}{k!}e^{-\lambda}$, $k=0,1,2,\ldots$, the probability of observing $k$ arrivals in a Poisson process with rate $\lambda$. Let
\begin{eqnarray}
	p(\widetilde n)&=& 
	1-\sum_{k=0}^{n-1} p(k)
	=1-\sum_{k=0}^{n-1} \frac{\lambda^k}{k!}e^{-\lambda},
	\quad n=1,2,\ldots,
\end{eqnarray} 
with $p(\widetilde 0)=1$, interpreted as the probability that at least $n$ arrivals were observed; and
\begin{eqnarray}
	p(\ell_1,\ldots,\ell_I\ | \ n)&=&
	\frac{n!}{\prod_{i=1}^I \ell_i!}
	\times
	\prod_{i=1}^I 
	\left( \frac{\lambda_i}{\lambda} \right)^{\ell_i}
	=\frac{n!}{\lambda^n}
	\times
	\prod_{i=1}^I 
	\frac{\lambda_i^{\ell_i}}{\ell_i!}
	\quad n=1,2,\ldots,
\end{eqnarray} 
with $\ell_1+\ldots+\ell_I=n$, interpreted as the conditional probability that, given $n$ patients have been admitted, there were $\ell_1,\ldots,\ell_I$ patients of type $1,\ldots,I$, respectively.

Further, let $p^{(n)}_{(\ell_1,\ldots,\ell_I)}$ be the probability that $n$ patients have been admitted and there were $\ell_1,\ldots,\ell_I$ patients of type $1,\ldots,I$, respectively, with $\ell_1+\ldots+\ell_I=n$. Then, if the number of admissions is restricted by $\ell_1+\ldots+\ell_I\leq N$, where $N$ is the number of available beds, we have,
\begin{eqnarray}
	p^{(n)}_{(\ell_1,\ldots,\ell_I)}&=&p(\ell_1,\ldots,\ell_I\ |\ n)p(n)
	=
	e^{-\lambda}
	\times \prod_{i=1}^I 
	\frac{\lambda_i^{\ell_i}}{\ell_i!}
	\quad\mbox{ for } n< N,\\
	p^{(\widetilde N)}_{(\ell_1,\ldots,\ell_I)}&=&p(\ell_1,\ldots,\ell_I\ |\ N)p(\widetilde N)
	=
	\left(
	1-\sum_{k=0}^{N-1} \frac{\lambda^k}{k!}e^{-\lambda}
	\right)
	\times 
	\frac{N!}{\lambda^N}
	\times
	\prod_{i=1}^I 
	\frac{\lambda_i^{\ell_i}}{\ell_i!},
\end{eqnarray}  
and $p^{(n)}_{(\ell_1,\ldots,\ell_I)}=0$ otherwise, with $p^{(\widetilde 0)}_{(0,0)}=1$.

Also, denote,
\begin{eqnarray}
	{\bf p}
	&=&
	\left[
	\begin{array}{cccccc}
		p^{(0)}_{(0,0)}& p^{(1)}_{(1,0)}& p^{(1)}_{(0,1)}& p^{(\widetilde 2)}_{(2,0)} & p^{(\widetilde 2)}_{(1,1)} & p^{(\widetilde 2)}_{(0,2)}
	\end{array}
	\right],\\
	{\bf w}
	&=&
	\left[
	\begin{array}{ccc}
		p^{(0)}_{(0,0)}& p^{(\widetilde 1)}_{(1,0)}& p^{(\widetilde 1)}_{(0,1)}
	\end{array}
	\right]
	.
\end{eqnarray}
Then, for $a=1$, the nonzero entries of ${\bf P}^{(1)}$ are,
\begin{eqnarray}
	\left[
	\begin{array}{ccc}
		{\bf P}^{(1)}_{1,1}& \ldots & {\bf P}^{(1)}_{1,6}
	\end{array}
	\right]
	&=&
	{\bf p},
\end{eqnarray}
for $s=2,7$,
\begin{eqnarray}
	\left[
	\begin{array}{ccc}
		{\bf P}^{(1)}_{s,1}& \ldots & {\bf P}^{(1)}_{s,6}
	\end{array}
	\right]
	&=&
	p_{1,1} \times {\bf p},
	\\
	\left[
	\begin{array}{ccc}
		{\bf P}^{(1)}_{s,7}& \ldots & {\bf P}^{(1)}_{s,9}
	\end{array}
	\right]
	&=& q_{1,1} \times {\bf w},
\end{eqnarray}
for $s=3,16$,
\begin{eqnarray}
	\left[
	\begin{array}{ccc}
		{\bf P}^{(1)}_{s,1}& \ldots & {\bf P}^{(1)}_{s,6}
	\end{array}
	\right]
	&=& p_{2,2} \times {\bf p},
	\\
	\left[
	\begin{array}{ccc}
		{\bf P}^{(1)}_{s,16}& \ldots & {\bf P}^{(1)}_{s,18}
	\end{array}
	\right]
	&=& 
	q_{2,2}\times {\bf w}
	,
\end{eqnarray}
for $s=4,8,14,19$,
\begin{eqnarray}
	\left[
	\begin{array}{ccc}
		{\bf P}^{(1)}_{s,1}& \ldots & {\bf P}^{(1)}_{s,6}
	\end{array}
	\right]
	&=& p_{1,1}p_{2,1} \times {\bf p},
	\\
	\left[
	\begin{array}{ccc}
		{\bf P}^{(1)}_{s,7}& \ldots & {\bf P}^{(1)}_{s,9}
	\end{array}
	\right]
	&=&
	q_{1,1}p_{2,1}\times {\bf w}\\
	\left[
	\begin{array}{ccc}
		{\bf P}^{(1)}_{s,13}& \ldots & {\bf P}^{(1)}_{s,15}
	\end{array}
	\right]
	&=&
	p_{1,1}q_{2,1}\times {\bf w},\\
	{\bf P}^{(1)}_{s,19}&=&q_{1,1}q_{2,1},
\end{eqnarray}
for $s=5,9,17,21$,
\begin{eqnarray}
	\left[
	\begin{array}{ccc}
		{\bf P}^{(1)}_{s,1}& \ldots & {\bf P}^{(1)}_{s,6}
	\end{array}
	\right]
	&=& p_{1,1}p_{2,2} \times {\bf p},\\
	\left[
	\begin{array}{ccc}
		{\bf P}^{(1)}_{s,7}& \ldots & {\bf P}^{(1)}_{s,9}
	\end{array}
	\right]
	&=&
	q_{1,1}p_{2,2}\times {\bf w},\\
	\left[
	\begin{array}{ccc}
		{\bf P}^{(1)}_{s,16}& \ldots & {\bf P}^{(1)}_{s,18}
	\end{array}
	\right]
	&=&
	p_{1,1}q_{2,2}\times {\bf w},\\
	{\bf P}^{(1)}_{s,21}&=&q_{1,1}q_{2,2},\\
\end{eqnarray}
for $s=6,12,18,22$,
\begin{eqnarray}
	\left[
	\begin{array}{ccc}
		{\bf P}^{(1)}_{s,1}& \ldots & {\bf P}^{(1)}_{s,6}
	\end{array}
	\right]
	&=& p_{1,2}p_{2,2} \times {\bf p},\\
	\left[
	\begin{array}{ccc}
		{\bf P}^{(1)}_{s,10}& \ldots & {\bf P}^{(1)}_{s,12}
	\end{array}
	\right]
	&=& 
	q_{1,2}p_{2,2}\times {\bf w},\\
	\left[
	\begin{array}{ccc}
		{\bf P}^{(1)}_{s,16}& \ldots & {\bf P}^{(1)}_{s,18}
	\end{array}
	\right]
	&=&
	p_{1,2}q_{2,2}\times {\bf w},\\
	{\bf P}^{(1)}_{s,22}&=&q_{1,2}q_{2,2},
\end{eqnarray}
for $s=10$,
\begin{eqnarray}
	\left[
	\begin{array}{ccc}
		{\bf P}^{(1)}_{10,1}& \ldots & {\bf P}^{(1)}_{10,6}
	\end{array}
	\right]
	&=& p_{1,2} \times {\bf p},\\
	\left[
	\begin{array}{ccc}
		{\bf P}^{(1)}_{10,10}& \ldots & {\bf P}^{(1)}_{10,12}
	\end{array}
	\right]
	&=&
	q_{1,2}\times {\bf w},\\
\end{eqnarray}
for $s=11,15,20$,
\begin{eqnarray}
	\left[
	\begin{array}{ccc}
		{\bf P}^{(1)}_{s,1}& \ldots & {\bf P}^{(1)}_{s,6}
	\end{array}
	\right]
	&=& p_{1,2}p_{2,1} \times {\bf p},\\
	\left[
	\begin{array}{ccc}
		{\bf P}^{(1)}_{s,10}& \ldots & {\bf P}^{(1)}_{s,12}
	\end{array}
	\right]
	&=&
	q_{1,2}p_{2,1}\times {\bf w},
	\\
	\left[
	\begin{array}{ccc}
		{\bf P}^{(1)}_{s,13}& \ldots & {\bf P}^{(1)}_{s,15}
	\end{array}
	\right]
	&=&
	p_{1,2}q_{2,1}\times {\bf w},
	\\
	{\bf P}^{(1)}_{s,20}&=&
	q_{1,2}q_{2,1},
\end{eqnarray}
and for $s=13$,
\begin{eqnarray}
	\left[
	\begin{array}{ccc}
		{\bf P}^{(1)}_{13,1}& \ldots & {\bf P}^{(1)}_{13,6}
	\end{array}
	\right]
	&=& p_{2,1} \times {\bf p},\\
	\left[
	\begin{array}{ccc}
		{\bf P}^{(1)}_{13,13}& \ldots & {\bf P}^{(1)}_{13,15}
	\end{array}
	\right]
	&=&
	q_{2,1}\times {\bf w}.
\end{eqnarray}

For $a=2$, rows $s=11,15$ are the only rows in ${\bf P}^{(2)}$ that differ from the rows of ${\bf P}^{(1)}$, with their nonzero entries given by,
\begin{eqnarray}
	\left[
	\begin{array}{ccc}
		{\bf P}^{(1)}_{s,1}& \ldots & {\bf P}^{(1)}_{s,6}
	\end{array}
	\right]
	&=& p_{1,1}p_{2,2} \times {\bf p},\\
	\left[
	\begin{array}{ccc}
		{\bf P}^{(1)}_{s,7}& \ldots & {\bf P}^{(1)}_{s,9}
	\end{array}
	\right]
	&=&
	q_{1,1}p_{2,2}\times {\bf w},
	\\
	\left[
	\begin{array}{ccc}
		{\bf P}^{(1)}_{s,16}& \ldots & {\bf P}^{(1)}_{s,18}
	\end{array}
	\right]
	&=&
	p_{1,1}q_{2,2}\times {\bf w},
	\\
	{\bf P}^{(1)}_{s,21}&=&
	q_{1,1}q_{2,2},
\end{eqnarray}
and the remaining rows are the same as in ${\bf P}^{(1)}$.

Finally, costs $C(s,a)$ are recorded as vectors ${\bf C}^{(a)}=[C(s,a)]_{s=1,\ldots,22}$ corresponding to decisions $a=1,2$, such that,
\begin{eqnarray}
	C(s,1)&=& 2c^{(\sigma)}I\{s= 4,5,6 \}
	+c^{(\sigma)}I\{s=2,3,8,9,11,12,14,15,17,18\}
	\nonumber\\
	&&
	+ c^{(p)}I\{s=4,6,8,10,12, 13,14,18,19,22\}
	+ 2 c^{(p)}I\{s=11,15,20\}
	\nonumber\\
	C(s,2)&=& 2c^{(\sigma)}I\{s= 4,5,6 \}
	+c^{(\sigma)}I\{s=2,3,8,9,11,12,14,15,17,18\}
	\nonumber\\
	&&
	+ c^{(p)}I\{s=4,6,8,10,12, 13,14,18,19,22\}
	+ 2 c^{(p)}I\{s=20\}+ c^{(t)}I\{s=11,15\},
	\nonumber\\
\end{eqnarray}
since the costs do not depend on the next observed state.

\end{document}